\documentclass[11pt]{article}
\usepackage{amssymb}
\usepackage{amsmath}
\usepackage{amsthm}
\newcommand{\cA}{{\cal A}}
\newcommand{\cB}{{\cal B}}
\newcommand{\cC}{{\cal C}}
\newcommand{\cD}{{\cal D}}
\newcommand{\cH}{{\cal H}}
\newcommand{\cE}{{\cal E}}
\newcommand{\cG}{{\cal G}}
\newcommand{\cI}{{\cal I}}

\newcommand{\cO}{{\cal O}}
\newcommand{\cL}{{\cal L}}
\newcommand{\cM}{{\cal M}}

\newcommand{\cF}{{\cal F}}
\newcommand{\cK}{{\cal K}}
\newcommand{\cP}{{\cal P}}

\newcommand{\cS}{{\cal S}}
\newcommand{\cT}{{\cal T}}
\newcommand{\cU}{{\cal U}}
\newcommand{\cV}{{\cal V}}
\newcommand{\cW}{{\cal W}}
\newcommand{\cX}{{\cal X}}
\newcommand{\cY}{{\cal Y}}
\newcommand{\cZ}{{\cal Z}}

\renewcommand{\AA}{{\mathbb A}}
\newcommand{\GG}{{\mathbb G}}

\newcommand{\ZZ}{{\mathbb Z}}
\newcommand{\QQ}{{\mathbb Q}}
\newcommand{\CC}{{\mathbb C}}
\newcommand{\PP}{{\mathbb P}}

\newcommand{\SSS}{{\mathbb S}}

\newcommand{\bb}{\mathfrak{b}}
\renewcommand{\gg}{\mathfrak{g}}  
\newcommand{\gh}{\mathfrak{h}}

\newcommand{\gp}{\mathfrak{p}}

\newcommand{\on}{\operatorname}

\newcommand{\Rep}{{\on{Rep}}}

\newcommand{\Qlb}{\mathbb{\bar Q}_\ell}
\newcommand{\Gm}{\mathbb{G}_m}

\newcommand{\A}{\mathbb{A}}

\newcommand{\toup}[1]{\stackrel{#1}{\to}}
\newcommand{\hook}[1]{\stackrel{#1}{\hookrightarrow}}

\newcommand{\getsup}[1]{\stackrel{#1}{\gets}}
\newcommand{\Sp}{\on{\mathbb{S}p}}
\newcommand{\Spin}{\on{\mathbb{S}pin}}
\newcommand{\GSpin}{\on{G\mathbb{S}pin}}

\newcommand{\IC}{\on{IC}}

\newcommand{\Hom}{\on{Hom}}
\newcommand{\Mod}{\on{Mod}}
\newcommand{\Ext}{\on{Ext}}
\newcommand{\End}{\on{End}}

\newcommand{\Sym}{\on{Sym}}
\newcommand{\SO}{\on{S\mathbb{O}}}

\newcommand{\Ker}{\on{Ker}}

\newcommand{\Aut}{\on{Aut}}

\newcommand{\RG}{\on{R\Gamma}}

\newcommand{\Bun}{\on{Bun}}
\newcommand{\uBun}{\on{\underline{Bun}}}
\newcommand{\Bunb}{\on{\overline{Bun}} }
\newcommand{\Bunt}{\on{\widetilde\Bun}}

\newcommand{\Spec}{\on{Spec}}

\newcommand{\HOM}{{{\cal H}om}}

\newcommand{\Gr}{\on{Gr}}

\newcommand{\GL}{\on{GL}}

\newcommand{\pr}{\on{pr}}
\newcommand{\id}{\on{id}}
\newcommand{\tr}{\on{tr}}
\newcommand{\QED}{$\square$} 
\newcommand{\Fq}{\mathbb{F}_q}  
\newcommand{\Fp}{\mathbb{F}_p}  
\newcommand{\iso}{{\widetilde\to}}

\newcommand{\comp}{\circ}
\newcommand{\Four}{\on{Four}}
\renewcommand{\H}{{\on{H}}}   

\newcommand{\DD}{\mathbb{D}}  
\newcommand{\D}{\on{D}}       
\newcommand{\DP}{\on{DP}}
\newcommand{\wt}{\widetilde}
\newcommand{\ov}[1]{\overline{#1}}
\newcommand{\select}[1]{{\it{#1}}}

\renewcommand{\div}{\on{div}}

\renewcommand{\P}{{\on{P}}}

\newcommand{\<}{\langle}
\renewcommand{\>}{\rangle}

\newcommand{\ev}{\mathit{ev}}

\newcommand{\Loc}{\on{Loc}}

\newcommand{\Res}{\on{Res}}

\newcommand{\act}{\on{act}}
\newcommand{\dimrel}{\on{dim.rel}}

\renewcommand{\Im}{\on{Im}}
\newcommand{\codim}{\on{codim}}

\newcommand{\SL}{\on{SL}}
\newcommand{\St}{\on{St}}
\newcommand{\tboxtimes}{\,\tilde\boxtimes\,}

\newcommand{\ra}{\rightarrow}
\newcommand{\la}{\leftarrow}

\newcommand{\diag}{\on{diag}}

\newcommand{\ST}{\on{ST}}

\newcommand{\Div}{\on{Div}}
\newcommand{\ocS}{\overset{\scriptscriptstyle\circ}{\cS}}
\newcommand{\ocY}{\overset{\scriptscriptstyle\circ}{\cY}}
\newcommand{\ocZ}{\overset{\scriptscriptstyle\circ}{\cZ}}
\newcommand{\ocL}{\overset{\scriptscriptstyle\circ}{\cL}}

\newtheorem{Lm}{Lemma}
\newtheorem{Th}{Theorem}
\newtheorem{Pp}{Proposition}
\newtheorem{Cor}{Corollary}
\newtheorem{Con}{Conjecture}

\theoremstyle{remark}
\newtheorem{Rem}{Remark}

\theoremstyle{definition}
\newtheorem{Def}{Definition}

\newenvironment{Prf}{\par\noindent {\it Proof }}{\QED}
\newcommand{\Step}[1]{\par\noindent{\bf Step {#1}}.}

\begin{document}
\author{Vincent Lafforgue and Sergey Lysenko}
\title{Geometrizing the minimal representations of even orthogonal groups}
\date{}
\maketitle
\begin{abstract}
\noindent{\scshape Abstract}\hskip 0.8 em Let $X$ be a smooth projective curve. Write $\Bun_{\SO_{2n}}$ for the moduli stack of $\SO_{2n}$-torsors on $X$. We give a geometric interpretation of the automorphic function $f$ on $\Bun_{\SO_{2n}}$ corresponding to the minimal representation. Namely, we construct a perverse sheaf $\cK_H$ on $\Bun_{\SO_{2n}}$ such that $f$ should be equal to the trace of Frobenius of $\cK_H$ plus some constant function. The construction is based on some explicit geometric formulas for the Fourier coefficients of $f$ on one hand, and on the geometric theta-lifting on the other hand. Our construction makes sense for more general simple algebraic groups, we formulate the corresponding conjectures. They could provide a geometric interpretation of some unipotent automorphic representations in the framework of the geometric Langlands program.
\end{abstract}

\bigskip

\centerline{\scshape 1. Introduction}

\bigskip\noindent
1.1 The theory of minimal representations has been developped (at least since 1989) in the works of D. Kazhdan, G. Savin, W.T. Gan, D. Ginzburg, S. Rallis, D. Soudry and others (cf. \cite{GS} for a recent survey) in several settings, over finite, local and global fields. In the theory of automorphic forms they are of special interest as they allow to prove some particular cases of Langlands functoriality via `generalized theta correspondences'. 

 The first example of a minimal representation is the Weil representation of the metaplectic group. In \cite{L1} a geometric version of the corresponding automorphic theta-function was constructed. In the present paper we develop a similar geometric theory for the minimal representations of even orthogonal groups. One of our motivations is a hope that the automorphic sheaves corresponding to the minimal representations will yield new cases of the geometric Langlands functoriality, as in the classical theory. For example, this should be the case for the dual pair $(\SO_3, \SO_{2n-3})$ in $\SO_{2n}$.   
 
 The place of minimal representations becomes clearer from the perspective of Arthur conjectures \cite{A}, they are particular examples of unipotent automorphic representations. 
 
 Let $k$ be a finite field. Let $X$ be a smooth projective geometrically connected curve over $k$. Let $H$ be a simple split group. Let $T\subset B\subset H$ be a maximal torus and a Borel subgroup. Write $\Lambda$ for the coweight lattice of $T$. Write $\check{H}$ for the Langlands dual group of $H$ over $\Qlb$. Set $F=k(X)$. Let $\AA$ be the adeles ring of $F$, $\cO\subset\AA$ be the entire adeles. For $x\in X$ write $F_x$ for the completion of $F$ at $x$.
 
 The unipotent representations we are interested in have been studied, in particular, in \cite{Mo}. Moeglin considers irreducible representations $\pi$ of $H(\AA)$ appearing as a direct summand in $L^2(H(F)\backslash H(\AA))$, which are everywhere nonramified and satisfy the following assumption. There is a character $\chi: T(\AA)/T(\cO)\to \Qlb^*$ that decomposes as 
$$
T(\AA)/T(\cO)\,\iso\, \Div(X)\otimes\Lambda\toup{\deg\otimes\id}\Lambda\to \Qlb^*
$$
such that $\pi$ appears in the induced representation $ind_{B(\AA)}^{H(\AA)}\chi$. They are expected to correspond to homomoprhisms $\phi: \SL_2\to \check{H}$ such that the corresponding unipotent $\check{H}$-orbit does not intersect any proper Levi subgroup. 

 Namely, let $\phi: \SL_2\to \check{H}$ satisfy this property. Let $\pi_{\phi}=\otimes'_{x\in X} \pi_x$, where $\pi_x$ is the spherical representation of $H(F_x)$ with Langlands parameter
$$
\phi\left(
\begin{array}{cc}
\mid t_x\mid^{\frac{1}{2}} & 0\\
0 & \mid t_x\mid^{-\frac{1}{2}}
\end{array}
\right),
$$
where $t_x\in F_x$ is a uniformizer. If $H$ is one of the split groups $\SO_{2n}, \SO_{2n+1}, \Sp_{2n}$ then, as Moeglin proved, $\pi_{\phi}$ appears in $L^2(H(F)\backslash H(\AA))$ as a direct summand with multiplicity one. We also expect this to hold in type $E_n$. 

 Let $\Bun_H$ denote the stack of $H$-torsors on $X$. The problem we are interested in is to find an object $K_{\phi}\in \D(\Bun_H)$ of the derived category of $\Qlb$-sheaves on $\Bun_H$, which is a geometric analog of $\pi_{\phi}$. Let $\sigma:\Gm\to \check{H}$ denote the restriction of $\phi$ under the map $\Gm\to\SL_2$, $x\mapsto \diag(x, x^{-1})$. Then $K_{\phi}$ should be a $\sigma$-Hecke eigensheaf as defined in (\cite{L2}, Definition~1). 

 If $H$ is of type $D_n$ or $E_n$ then the subregular unipotent orbit in $\check{H}$ does not intersect any proper Levi subgroup, and the corresponding representation $\pi_{\phi}$ of $H(\AA)$ is the minimal one. Its Arthur parameter is the homomorphism $\id\times \phi: \pi_1(X)\times\SL_2\to \check{H}$, where $\phi$ corresponds to the subregular unipotent orbit.
 
  In Appendix~A we introduce a notion of an almost constant local system on $\Bun_H$. We think they are nothing but the automorphic sheaves on $\Bun_H$ corresponding to the Arthur parameters of the form $\alpha\times\phi_p: \pi_1(X)\times\SL_2\to \check{H}$, where $\phi_p: \SL_2\to\check{H}$ is principal, and $\alpha$ factors through the center $Z(\check{H})$ of $\check{H}$. Conjecturally, any local system on $\Bun_H$ is almost constant (after passing from $k$ to its algebraic closure). Denote by $\D(\Bun_H)_{ls}\subset \D(\Bun_H)$ the full triangulated subcategory generated by the almost constant local systems. It is preserved by Hecke functors, so they also act on the quotient category $\D(\Bun_H)/\D(\Bun_H)_{ls}$. 
  
  Assume $H=\SO_{2n}$ split with $n\ge 4$. Let $\phi: \SL_2\to \check{H}$ correspond to the subregular unipotent orbit. Our main results are Theorems~\ref{Th1} and \ref{Th_2}, they provide a perverse sheaf $\cK_H\in \D(\Bun_H)$ irreducible on each connected component of $\Bun_H$ and such that its image in $\D(\Bun_H)/\D(\Bun_H)_{ls}$ satisfies the Hecke property for the Arthur parameter $\id\times\phi: \pi_1(X)\times\SL_2\to \check{H}$. So, the object $K_{\phi}$, which we are not able to find yet, will have the same image as $\cK_H$ in $\D(\Bun_H)/\D(\Bun_H)_{ls}$. 
  
  In the classical setting one could take the orthogonal complement to the vector space generated by the `almost constant functions' on $\Bun_H$, and find a function in this orthogonal complement which coincides with the trace of Frobenius of $\cK_H$ modulo the `almost constant functions' on $\Bun_H$. In the geometric setting the problem of lifting of $\cK_H\in \D(\Bun_H)/\D(\Bun_H)_{ls}$ to $K_{\phi}\in \D(\Bun_H)$ looks more difficult.
   
\medskip\noindent
1.2 One can construct the minimal representation $\pi$ of $H(\AA)$ as the restricted tensor product of minimal representations $\pi_x$ of $H(F_x)$, $x\in X$. Each $\pi_x$ is the local theta-lift of the trivial representation of $\SL_2(F_x)$ for the dual pair $(\SL_2, H)$ in the metaplectic group $\wt\Sp_{4n}$. The corresponding global theta-lift is divergent, and one has to take a residue of the corresponding series to obtain $\pi$. Another way is to realize $\pi$ as a residue of Eisenstein series for some parabolic subgroups of $H$ (see \cite{GRS}). Neither of these constructions admits an evident geometrization. It is not clear what a residue of a geometric Eisenstein series (or a divergent theta-series) should be in general, this was one of the key technical difficulties in this paper. Our construction suggests, a posteriori, a possible geometric approach to residues of Eisenstein series at least in the simplest cases (see Section~1.3).

 Our construction of $\cK_H$ is based on the theta-lifting. Let $\Omega$ be the canonical line bundle on $X$. Let $G_1$ be the group scheme on $X$ of automorphisms of $\cO_X\oplus\Omega$ acting trivially on $\det(\cO_X\oplus\Omega)$. Let $P\subset H$ be a parabolic subgroup preserving some isotropic $n$-dimensional subspace in the standard representation of $H$. 
 
 In (\cite{L2}, Definition~2) the theta-lifting functor $F_H: \D^-(\Bun_{G_1})_!\to \D^{\prec}(\Bun_H)$ has been introduced, it is given by the kernel $\Aut_{G_1, H}$ on $\Bun_{G_1}\times\Bun_H$, which is the restriction of the theta-sheaf for $\wt\Sp_{4n}$. Let $q: \Bun_{G_1}\times\Bun_H\to\Bun_H$ be the projection.  In \select{loc.cit.} we considered the morphism $\kappa: \check{G}_1\times\SL_2\to \check{H}$ given by 
$$
\SO_3\times\SL_2\toup{\id\times \phi_p} \SO_3\times \SO_{2n-3}\to \SO_{2n},
$$
where the latter map is given by an orthogonal direct sum, and $\phi_p: \SL_2\to \SO_{2n-3}$ is principal. By (\cite{L2}, Theorem~3), $F_H$ commutes with the Hecke functors with respect to $\kappa$. So, if we had $\Qlb\in \D^-(\Bun_{G_1})_!$ then $F_H(\Qlb)$ would be the automorphic sheaf $K_{\phi}$ we are looking for. However, $\Qlb$ is not in $\D^-(\Bun_{G_1})_!$, and the complex $F_H(\Qlb)$, which is $q_!\Aut_{G_1,H}$ up to a shift, does not make sense in the existing formalism. It is not bounded from above neither from below. 
  
  One can however, formally look at its Fourier coefficients with respect to $P$. The stack $\Bun_P$ of $P$-torsors on $X$ is the stack classifying $U\in\Bun_n$ and an exact sequence $0\to \wedge^2 U\to ?\to \cO_X\to 0$. Let $\cY_P$ be the stack classifying $U\in\Bun_n$ and a section $v: \wedge^2 U\to\Omega$. Then $\cY_P$ and $\Bun_P$ are dual generalized vector bundles over $\Bun_n$, the stack of rank $n$ vector bundles on $X$. From the explicit formulas for $\Aut_{G_1,H}$ in the Shr\"odinger model we noticed that \select{all the infiniteness} of the Fourier transform $\Four(q_!\Aut_{G_1,H})$ is concentrated on the zero section of $\cY_P$. This, together with the results about minimal representations (\cite{GRS}) has led to our construction of $\cK_H$ via the $P$-model (cf. Theorem~\ref{Th1}). 

  Let $Q$ (resp., $R$) be the parabolic subgroup of $H$ preserving a 1-dimensional (resp., $2$-dimensional) isotropic subspace in the standard representation. We also propose conjectural constructions of the same perverse sheaf $\cK_H$ via $Q$ and $R$-models. While the $Q$-model is a part of more general Conjecture~\ref{Con_general_descent} for simple groups admitting a parabolic subgroup with an abelian unipotent radical, the $R$-model plays a separate role. This is a parabolic subgroup of $H$ referred to as \select{Heisenberg parabolic} in \cite{GS}, such parabolic exists for any simple algebraic group. We hope that our explicit construction via $R$-model will generalize to cover the cases when there is no parabolic subgroup with an abelian unipotent radical (like $E_8$, for example). 

 In all three cases we construct some complexes $K_{P,\psi}$, $K_{Q,\psi}$, $K_{R,\psi}$ on $\Bun_P, \Bun_Q$ and $\Bun_R$ respectively given by some explicit formulas. We expect that each of this complexes is the restriction of $\cK_H$ from $\Bun_H$ (over suitable open substacks). Here $\psi: \Fp\to \Qlb^*$ is a nontrivial additive character. For the parabolic $P$ this is true and is a part of our construction of $\cK_H$ (cf. Theorem~\ref{Th1}).

   Our $R$-model uses as an input a new ingredient, the extended theta sheaf. This is a perverse sheaf interesting on its own ground, as it is a geometric analog of the matrix coefficient of the Weil representation of the semi-direct product $\wt\Sp_{2n}\rtimes H_n$, where $H_n$ is the Heisenberg group. Our definition of the perverse sheaf $K_{R,\psi}$ on $\Bun_R$ (cf. Section~2.3.3) is motivated by our compatibility result between $P$ and $R$-models (cf. Corollary~\ref{Cor_PR_compatibility}).
   
   In Section~8 we propose one more conjectural construction of $\cK_H$ via residues of geometric Eisenstein series (this approach is formalized in Section~1.3). We then apply this to calculate $\cK_H$ explicitly in the cases of genus zero and one (Propositions~\ref{Pp_answer_for_g=0} and \ref{Pp_case_g=1}). This calculation shows, in particular, that the generic rank of $\cK_H$ (for a given connected component of $\Bun_H$) does depend on the genus of $X$. In this sense the sheaf $\cK_H$ is not of local nature (as opposed to the case of the theta-sheaf for the metaplectic group).
   
   In the main body of the paper we work with \'etale $\Qlb$-sheaves in positive characteristic, however our Theorem~\ref{Th1} holds also for $\cD$-modules in characteristic zero (in Appendix~C we briefly explain how to change the proof of Theorem~\ref{Th1} in this case).
   
\medskip\noindent
1.3 {\scshape A geometric approach to residues} Let $G$ be a simple, simply-connected group, $P\subset G$ a maximal parabolic subgroup with Levi quotient $M$. The connected components of $\Bun_P$ are indexed by $\pi_1(M)$ naturally, write $\Bun_P^m$ for the connected component given by $m\in\pi_1(M)$. Note that $\pi_1(M)\,\iso\, \ZZ$. Set $\cS^m=(\nu_P^m)_!\Qlb$, where $\nu_P^m: \Bun_P^m\to\Bun_G$ is the natural map.
  
  Consider the induced representation $ind_{P(\AA)}^{G(\AA)}(\chi^s)$, $s\in \CC$, where 
$$
\chi: \PP(\AA)\to (M/[M,M])(\AA)\,\iso\, \AA^*\to \QQ^*,
$$
the last map sends $a\in \AA^*$ to $\mid\! a\!\mid$. The simplest residual representations appear inside these induced representations as non ramified subquotients at points $s\in\CC$ of reducibility. These points, for an appropriate normalization, correspond also to the (simple) poles of the Eisenstein series $E_P^G(s)$. 
  
  We suggest that for such $s$ there is an affine function $\alpha: \ZZ\to\ZZ$, $\alpha(m)=am+b$ such that for $m$ small enough, the perverse sheaf $^p\cH^{\alpha(m)}(\cS^m)$ stabilizes (or at least, contains the same irreducible perverse sheaf $\cK$ as a subquotient). Then say that the sequence $\cS^r$ \select{has a residue in the direction $\alpha$}.   
              
\bigskip

\centerline{\scshape 2. Main results}

\bigskip\noindent
2.1 {\scshape Notation} Let $k$ be an algebraically closed field of characteristic $p>2$ (except in Section~3.3 and 7.1-7.4, where we assume $k=\Fq$ finite with $q$ odd). All the schemes or stacks we consider are defined over $k$.

 Let $X$ be a smooth projective geometrically connected curve of genus $g$. Write $\Omega$ for the canonical line bundle on $X$. Fix a prime $\ell\ne p$. For a stack $S$ locally of finite type write $\D(S)$ for the category introduced in (\cite{LO}, Remark~3.21) and denoted $\D_c(S,\Qlb)$ in \select{loc.cit.} It should be thought of as the unbounded derived category of constructible $\Qlb$-sheaves on $S$. For $*=-,b$ one has the full subcategory $\D^*(S)\subset \D(S)$ denoted $\D^*_c(S,\Qlb)$ in \select{loc.cit.} Write $\D^{\prec}(S)\subset \D(S)$ for the full subcategory of complexes $K\in\D(S)$ such that for any open substack of finite type $U\subset S$ we have $K\mid_U\in \D^-(U)$. Write $\P(S)\subset \D(S)$ for the full subcategory of perverse sheaves.
 
  Fix a nontrivial character $\psi: \Fp\to\Qlb^*$ and denote by $\cL_{\psi}$ the corresponding Artin-Shreier sheaf on $\A^1$. Fix a square root $\Qlb(\frac{1}{2})$ of the sheaf $\Qlb(1)$ on $\Spec k$. If $k=\Fq$ the isomorphism class of such correspond to square roots of $q$ in $\Qlb$. For a morphism of stacks $f: Y\to Z$ write $\dimrel(f)$ for the function of a connected component $C$ of $Y$ given by $\dim C-\dim C'$, where $C'$ is the connected component of $Z$ containing $f(C)$.
  
  If $V\to S$ and $V^*\to S$ are dual rank $r$ vector bundles over a stack $S$, write $\Four_{\psi}: D^{\prec}(V)\to \D^{\prec}(V^*)$ for the Fourier transform given by $\Four_{\psi}(K)=(p_{V^*})_!(\xi^*\cL_{\psi}\otimes p_V^*K)[r](\frac{r}{2})$, where $p_V,p_{V^*}$ are the projections and $\xi: V\times_S V^*\to\A^1$ is the pairing.
  
  Write $\Bun_r$ for the stack of rank $r$ vector bundles on $X$. Our conventions about $\ZZ/2\ZZ$-gradings are those of (\cite{L1}, 3.1). For a group scheme $G$ over $X$ denote by $\Bun_G$ the stack of $G$-torsors on $X$. 
  
  For a connected reductive group $\cG$ over $\Qlb$ write $\Rep(\cG)$ for the category of finite-dimensional $\Qlb$-representations of $\cG$. Write $\DP(k)=\oplus_{d\in\ZZ} \P(\Spec k)[d]\subset \D(\Spec k)$ for the full subcategory in $\D(\Spec k)$. By definition, we have an equivalence of tensor categories $\Loc: \Rep(\Gm)\,\iso\, \DP(k)$ sending $\St^{\otimes m}$ to $\Qlb[m]$. Here $\St$ is the standard representation of $\Gm$.   
  
  If $k=\Fq$ then we denote $F=k(X)$, $\AA$ the adeles of $X$ and $\cO\subset\AA$ the entire adeles.
   
\medskip\noindent
2.2  {\scshape Extended theta sheaf}  For $n>0$ let $M_0=\cO^n\oplus\Omega^n$, write $G_n$ for the group scheme on $X$ of automorphisms of $M_0$ preserving the natural symplectic form $\wedge^2 M_0\to\Omega$. The stack $\Bun_{G_n}$ classifies $M\in\Bun_{2n}$ with symplectic form $\wedge^2 M\to\Omega$. Let $H_n=M_0\oplus\Omega$ be the corresponding Heisenberg group scheme on $X$, write $\GG_n=G_n\rtimes H_n$ for the corresponding semi-direct product (cf. Section~3.2 for details).
 
 Write $\cA_{G_n}$ for the line bundle on $\Bun_{G_n}$ with fibre $\det\RG(X,M)$ at $M$. We view it as a $\ZZ/2\ZZ$-graded purely of degree zero. Denote by $\Bunt_{G_n}\to\Bun_{G_n}$ the $\mu_2$-gerb of square roots of $\cA_{G_n}$. Write $\Aut$ for the perverse theta-sheaf on $\Bunt_{G_n}$ introduced in (\cite{L1}, Definition~1).
 
 The stack $\Bun_{\GG_n}$ classifies $M_1\in \Bun_{2n+2}$ with symplectic form $\wedge^2 M_1\to \Omega$ and a section $v: \Omega\hook{} M_1$ whose image is a subbundle. For a point of $\Bun_{\GG_n}$ write $L_{-1}$ for the orthogonal complement to $\Omega$, so $M=L_{-1}/\Omega\in \Bun_{G_n}$. Write $\rho_{\GG}: \Bun_{\GG_n}\to\Bun_{G_n}$ for the map sending $(M_1,v)$ to $M$. Set $\Bunt_{\GG_n}=\Bunt_{G_n}\times_{\Bun_{G_n}}\Bun_{\GG_n}$.

 Let $_0\Bun_{\GG_n}\subset \Bun_{\GG_n}$ be the open substack given by $\H^0(X,M)=0$, define $_0\Bun_{G_n}$, $_0\Bunt_{G_n}$, $_0\Bunt_{\GG_n}$ similarly. Write $\Bun_{\Omega}$ for the stack classifying exact sequences 
\begin{equation}
\label{seq_cO_by_Omega} 
0\to\Omega\to ?\to\cO\to 0
\end{equation}
on $X$. Let $\ev_{\Omega}: \Bun_{\Omega}\to\A^1$ be the map sending (\ref{seq_cO_by_Omega}) to the corresponding element of $\H^1(X,\Omega)$. We have canonically 
\begin{equation}
\label{iso_0_Bunt_GGn_is}
_0\Bunt_{\GG_n}\,\iso\, {_0\Bunt_{G_n}}\times\Bun_{\Omega}
\end{equation}
\begin{Def} Write $\Aut^e_{\psi}$ for the intermediate extension of 
$$
(\Aut\boxtimes \ev_{\Omega}^*\cL_{\psi})\otimes(\Qlb[1](\frac{1}{2}))^{1-g}
$$ 
under the open immersion $_0\Bunt_{\GG_n}\hook{}\Bunt_{\GG_n}$. Here $e$ stands for `extended', we call $\Aut^e_{\psi}$ \select{the extended theta-sheaf}.
\end{Def}
 
  Let $P_n\subset G_n$ be the parabolic group subscheme preserving $\cO^n\subset M_0$. Set $\PP_n=P_n\rtimes H_n$. 
The stack $\Bun_{\PP_n}$ classifies $\cL\in\Bun_n$ included into an exact sequence on $X$
\begin{equation}
\label{seq_cL_by_Omega}  
0\to\Omega\to\bar\cL\to\cL\to 0
\end{equation}
and an exact sequence on $X$
\begin{equation}
\label{seq_Omega_by_Sym2barcL}  
0\to\Sym^2\bar\cL\to ?\to\Omega\to 0
\end{equation}

 Let $\cT_n$ be the stack classifying $\cL\in\Bun_n$ and an exact sequence (\ref{seq_cL_by_Omega}) on $X$. Let $\cZ_{\cT_n}$ be the stack classifying a point of $\cT_n$ and a splitting $s: \bar\cL\to\Omega$ of (\ref{seq_cL_by_Omega}). 
 
 Write $\cZ_{2,\cT_n}$ for the stack classifying a point of $\cT_n$ as above and a section $\bar s: \Sym^2\bar\cL\to \Omega^2$. The map $h_{\cT}: \cZ_{\cT_n}\to \cZ_{2,\cT_n}$ over $\cT_n$ given by $\bar s=s\otimes s$ is a closed immersion. One has a diagram of dual generalized vector bundles $\cZ_{2,\cT_n}\to\cT_n\gets \Bun_{\PP_n}$ over $\cT_n$. Let $\Four_{\cZ_{\cT}, \psi}: \D^{\prec}(\cZ_{2,\cT_n})\to \D^{\prec}(\Bun_{\PP_n})$ denote the corresponding Fourier transform functor. Set
$$
K_{\PP_n,\psi}=\Four_{\cZ_{\cT}, \psi}h_{\cT !}\Qlb[a],
$$
where $a$ is a function of a connected component of $\cT_n$ given by $\dim\cT_n-\chi(\cL)$, $\cL\in\Bun_n$.

 Let $^0\Bun_n\subset\Bun_n$ be the open substack classifying $\cL\in\Bun_n$ with $\H^0(X,\Sym^2 \cL)=0$. Write $^0\Bun_{\PP_n}$ for the preimage of $^0\Bun_n$ under the map $\Bun_{\PP_n}\to \Bun_n$ sending the above point to $\cL$. We will define a morphism $\tilde\nu_{\PP}: \Bun_{\PP_n}\to \Bunt_{\GG_n}$ whose restriction to $^0\Bun_{\PP_n}$ is smooth (cf. Section~3.4.3).
 
\begin{Pp} 
\label{Pp_explicit_formula_Aute}
There is an isomorphism over $^0\Bun_{\PP_n}$
\begin{equation}
\label{iso_explicit_Aute}
\tilde\nu_{\PP}^*\Aut^e_{\psi}\otimes(\Qlb[1](\frac{1}{2}))^{\dimrel(\tilde\nu_{\PP})}\,\iso\, K_{\PP_n,\psi}
\end{equation}
\end{Pp}
\begin{Rem} We conjecture that the isomorphism (\ref{iso_explicit_Aute}) holds over the whole of $\Bun_{\PP_n}$.
\end{Rem}

 We also introduce a finite-dimensional analog of the extended theta-sheaf and calculate all the $*$-fibres of $\Aut^e_{\psi}$ (cf. Section~3).

In the case $k=\Fq$ we show that $\Aut^e_{\psi}$ is a geometric analog of the following matrix coefficient of the Weil representation. Let $\chi: \Omega(\AA)/\Omega(F)\to \Qlb^*$ be the character 
\begin{equation}
\label{character_chi_for_Fq}
\chi(\omega)=\psi(\sum_{x\in X} \tr_{k(x)/k} \Res \omega_x)
\end{equation}
Denote by $(\rho, \cS_{\psi})$ a (unique up to isomorphism) irreducible representation of $H_n(\AA)$ over $\Qlb$ with central character $\chi$. Let $\hat G_n(\AA)$ be the metaplectic group defined by this representation
\begin{multline*}
\hat G_n(\AA)=\{(g,\sigma)\mid g\in G_n(\AA), \sigma\in\Aut \cS_{\psi}, \; 
\rho(gm, \omega)\comp \sigma=\sigma \comp \rho(m,\omega)\;\mbox{for}\; (m,\omega)\in H_n(\AA)\}
\end{multline*}
The sequence is exact 
\begin{equation}
\label{seq_metaplectic_classical_def}
1\to \Qlb^*\to \hat G_n(\AA)\to G_n(\AA)\to 1
\end{equation}
Then $\cS_{\psi}$ is naturally a representation of $\hat \GG_n(\AA)=\hat G_n(\AA)\rtimes H_n(\AA)$. For a subgroup $\cK\subset G_n(\AA)$ write $\hat\cK$ for its preimage in $\hat G_n(\AA)$. One checks that $\cS_{\psi}$ admits a unique up to a multiple non zero $H_n(F)$-invariant functional $\Theta:\cS_{\psi}\to\Qlb$. The group $\hat G_n(F)$ acts naturally on the space of such functionals, this gives a splitting of (\ref{seq_metaplectic_classical_def}) over $G_n(F)$. View $G_n(F)\subset \hat G_n(\AA)$ as a subgroup. The representation $\cS_{\psi}$ also admits a unique up to a multiple $H_n(\cO)$-invariant vector $v_0$, it similarly yields a splitting of (\ref{seq_metaplectic_classical_def}) over $G_n(\cO)$. This yields the subgroups $\GG_n(F)$ and $\GG_n(\cO)$ of $\hat\GG_n(\AA)$. Let 
\begin{equation}
\label{autom_form_phi}
\phi: \GG_n(F)\backslash \hat\GG_n(\AA)/\GG_n(\cO)\to\Qlb
\end{equation}
be given by $\phi(g)=\Theta(gv_0)$, $g\in \hat\GG_n(\AA)$. Then $\Bunt_{\GG_n}$ can be thought of as a geometric analog of 
$$
\GG_n(F)\backslash \hat\GG_n(\AA)/\GG_n(\cO),
$$
and $\Aut^e_{\psi}$ is a geometric counterpart of $\phi$.

\medskip\noindent
2.3 {\scshape Models of the minimal sheaf} 

\medskip\noindent
Fix $n\ge 2$, let $H=\SO_{2n}$ be the split orthogonal group over $k$. We write $H_n$ when we need to express the dependence on $n$. Write $V_0$ for the standard representation of $H$. Write $\check{H}$ for the Langlands dual group over $\Qlb$, so $\check{H}\,\iso\,\SO_{2n}$.

The stack $\Bun_H$ classifies $V\in\Bun_{2n}$ with a non degenerate symmetric form $\Sym^2(V)\to\cO_X$ and a compatible trivialization $\det V\,\iso\, \cO_X$. One has $\pi_1(H)\,\iso\,\ZZ/2\ZZ$, and the connected components $\Bun_H^{\theta}$ of $\Bun_H$ are naturally indexed by $\theta\in \pi_1(H)$. 
 
\medskip\noindent
2.3.1 {\scshape $P$-model\ } Fix an $n$-dimensional isotropic subspace $U_0\subset V_0$. Let $P\subset H$ be the parabolic subgroup preserving $U_0$. The stack $\Bun_P$ classifies $U\in \Bun_n$ and an exact sequence  of $\cO_X$-modules
\begin{equation}
\label{seq_cO_by_wedge2U}
0\to \wedge^2 U\to ?\to \cO_X\to 0
\end{equation} 

 Let $\cY_P$ be the stack classifying $U\in \Bun_n$ and a section $v: \wedge^2 U\to\Omega$. So, $\cY_P$ and $\Bun_P$ are dual generalized vector bundles over $\Bun_n$. 
 
 Let $\cS_P$ be the stack classifying $U\in\Bun_n$, $M\in \Bun_{G_1}$ and a morphism $s: U\to M$. Let $\pi_P: \cS_P\to \cY_P$ be the morphism sending $(U, M,s)$ to $(U, v)$, where $v$ is the composition
$$
\wedge^2 U\toup{\wedge^2 s} \wedge^2 M\,\iso\,\Omega
$$
Let $\cZ_P\subset\cY_P$ be the closed substack classifying $(U,v)$ such that the generic rank of $v: U\to U^*\otimes\Omega$ is at most 2. This is equivalent to requiring that $\wedge^3 v: \wedge^3 U\to \wedge^3(U^*\otimes\Omega)$ vanishes. Clearly, $\pi_P$ factors as 
$$
\cS_P\toup{\pi_P} \cZ_P\hook{} \cY_P
$$

 Let $\ocY_P\subset\cY_P$ be the open substack given by the condition that $v\ne 0$. Let $\ocS_P$ and $\ocZ_P$ be the preimages of $\ocY_P$ in $\cS_P$ and $\cZ_P$ respectively. 
 
 For $d\ge 0$ write $X^{(d)}$ for the $d$-th symmetric power of $X$. Stratify $\ocZ_P$ by locally closed substacks $\cZ_{P, m}$ indexed by $m\ge 0$. Here $\cZ_{P,m}$ is given by the condition that there exists an effective divisor $D\in X^{(m)}$ such that $v: \wedge^2 U\to \Omega(-D)$ is surjective. Note that $\cZ_{P,0}\subset \ocZ_P$ is an open substack.  The stack $\cZ_{P,m}$ can be seen as the stack classifying $D\in X^{(m)}$, $U\in \Bun_n$, $M'\in\Bun_2$ together with a surjective morphism of $\cO_X$-modules $U\to M'$, and an isomorphism $\det M'\,\iso\, \Omega(-D)$.
 
  Write $\Bun_n^d$ for the connected component of $\Bun_n$ classifying $U\in\Bun_n$ with $\deg U=d$. Let $\Bun_P^d$, $\cZ^d_P$ and so on denote the preimage of $\Bun_n^d$ in the corresponding stack.
 
 The stack $\ocS_P$ is smooth. The restriction 
\begin{equation}
\label{map_pi_nice}
\pi_P: \ocS_P\to\ocZ_P
\end{equation}
of $\pi_P$ is representable, proper and surjective, this is an isomorphism over $\cZ_{P,0}$. For each $d\in\ZZ$ the stack $\ocS{}_P^d$ is irreducible, so $\ocZ{}_P^d$ is also irreducible.
 
\begin{Pp} 
\label{Pp_small_map}
1)  If $n\ge 4$ then the map (\ref{map_pi_nice}) is small, and one has canonically
\begin{equation}
\label{complex_pi_P_!} 
(\pi_P)_!\IC(\ocS_P)\,\iso\, \IC(\ocZ_P)
\end{equation}
2) If $n=3$ then (\ref{map_pi_nice}) is semi-small, and $\mathop{\oplus}\limits_{m\ge 0} \IC(\cZ_{P,m})$ is a direct summand of $(\pi_P)_!\IC(\ocS_P)$. 
\end{Pp}
 
 Let $^e\Bun_n\subset\Bun_n$ be the open substack given by two conditions $\H^0(X, \wedge^2 U)=0$ and $\H^0(X,\Omega\otimes\wedge^2 U)=0$ for $U\in\Bun_n$. Let $^e\cS_P$, $^e\Bun_P$, $^e\cY_P$ and so on be the preimages of $^e\Bun_n$ in the corresponding stacks. Write $\nu_P: \Bun_P\to\Bun_H$ for the map induced by $P\hook{}H$. Its restriction $^e\Bun_P\to\Bun_H$ is smooth.
 
  Let $Z(e,P)$ be the set of
$d\in\ZZ$ such that $^e\ocZ{}^d_P$ is not empty. There is $N\in\ZZ$ such that if $d\le N$ then $d\in Z(e,P)$. Since $\ocZ{}^d_P$ is irreducible, if $d\in Z(e,P)$ then $^e\ocZ{}^d_P\cap \cZ_{P,0}$ is also nonempty.
 
 Write $\Four_{\cY_P, \psi}: \D^{\prec}(\cY_P)\to \D^{\prec}(\Bun_P)$ for the Fourier transform functor over $\Bun_n$. Set 
$$
K_{P,\psi}=\Four_{\cY_P,\psi}\IC(\cZ_P)
$$

 Assume $n\ge 4$. We will construct a perverse sheaf $\cK_H$ on $\Bun_H$ irreducible on each connected component and defined up to a unique isomorphism (cf. Section~2.3.4 and Definition~\ref{Def_true_cK_H} in Section~7.3.1).
Here is our main result. 
 
\begin{Th}
\label{Th1}
For each $d\in Z(e,P)$ there exists an isomorphism over $^e\Bun_P^d$ 
\begin{equation}
\label{iso_Th1_on_eBunP}
\nu_P^*\cK_H\otimes(\Qlb[1](\frac{1}{2})^{\dimrel(\nu_P)}\,\iso\, K_{P,\psi}
\end{equation} 
\end{Th}

 Let us formulate a conjectural version of the Hecke property
of $\cK_H$. Write 
$$
\H^{\la}_H: \Rep(\check{H})\times \D(\Bun_H)\to \D(X\times\Bun_H)
$$ 
for the Hecke functors on $\Bun_H$ (cf. \cite{L2}, Section~2.2.1 for a precise definition).
Let $\sigma: \Gm\to \check{H}$ be the composition $\Gm\to \SL_2\to\check{H}$, where the second map corresponds to the subregular unipotent orbit, and the first one is the standard maximal torus. Let $E_0: \Rep(\check{H})\to \DP(\Spec k)$ be the functor $W\mapsto \Loc(\Res^{\sigma}(W))$, here $\Res^{\sigma}: \Rep(\check{H})\to \Rep(\Gm)$ is the restiction via $\sigma$. 

\begin{Con}
\label{Con_1}
There is a functor $E_1: \Rep(\check{H})\to \DP(\Spec k)$ and an isomorphism in $\D(X\times\Bun_H)$
$$
\H^{\la}_H(W, \cK_H)\,\iso\, (E_0(W)\otimes\cK_H)[1](\frac{1}{2})\oplus E_1(W)
$$
functorial in $W\in\Rep(\check{H})$.
\end{Con} 

 In Appendix~A we introduce a notion of an almost constant local system on $\Bun_H$. In Section~7.6 we prove the following weaker form of Conjecture~\ref{Con_1}.

\begin{Th}
\label{Th_2}
 Let $x\in X$. There is a functor $E_1: \Rep(\check{H})\to \D(\Bun_H)$ 
with the following properties. If $W\in \Rep(\check{H})$ then $E_1(W)$ is a direct sum of shifted almost constant local systems on $\Bun_H$. There is an isomorphism functorial in $W\in \Rep(\check{H})$
$$
_x\H^{\la}_H(W, \cK_H)\,\iso\, E_0(W)\otimes \cK_H \oplus E_1(W)
$$
\end{Th}

\begin{Rem} Consider $\Four_{\cY_P,\psi}^{-1}(\nu_P^*\cK_H)$ over the whole of $\cY_P$,
we expect that  it is the extension by zero from $\cZ_P$.
\end{Rem}
 
\medskip\noindent
2.3.2 {\scshape $Q$-model \ } Let $W_0\subset U_0$ be a 1-dimensional subspace. Let $Q\subset H$ be the parabolic subgroup preserving $W_0$. The stack $\Bun_Q$ classifies $V'\in \Bun_{H_{n-1}}$, $W\in \Bun_1$ and an exact sequence of $\cO_X$-modules
\begin{equation}
\label{seq_V'_by_W}
0\to W\to ?\to V'\to 0
\end{equation} 

 Let $\cY_Q$ be the stack classifying $W\in\Bun_1, V'\in \Bun_{H_{n-1}}$ and $t: W\to V'\otimes\Omega$. So, $\cY_Q$ and $\Bun_Q$ are generalized vector bundles over $\Bun_1\times\Bun_{H_{n-1}}$.
 
 Let $\cZ_Q\subset \cY_Q$ be the closed substack given by the condition that the composition
$$
W^{\otimes 2} \toup{t\otimes t} \Sym^2 (V'\otimes\Omega)\to \Omega^2,
$$  
vanishes, that is, the image of $t$ is isotropic. 
 
 Let $^u(\Bun_1\times \Bun_{H_{n-1}})\subset \Bun_1\times \Bun_{H_{n-1}}$ be the open substack given by $\H^0(X, V'\otimes W)=0$ and $\H^0(X, V'\otimes W\otimes\Omega)=0$ for $W\in\Bun_1, V'\in \Bun_{H_{n-1}}$. 
 
 Write $^u\cY_Q$, $^u\Bun_Q$, $^u\cZ_Q$ and so on for the preimages of $^u(\Bun_1\times \Bun_{H_{n-1}})$ in the corresponding stack.  Let $\nu_Q: \Bun_Q\to\Bun_H$ be the map induced by $Q\hook{}H$. Its restriction $^u\Bun_Q\to\Bun_H$ is smooth.
 
 Write $\Four_{\cY_Q, \psi}: \D^{\prec}(\cY_Q)\to \D^{\prec}(\Bun_Q)$ for the Fourier transform functor. Set
$$
K_{Q,\psi}=\Four_{\cY_Q,\psi}\IC(\cZ_Q)
$$ 
 
\begin{Con}
\label{Con2}
There exists an isomorphism over $^u\Bun_Q$
\begin{equation}
\label{iso_Th2_on_eBunQ}
\nu_Q^*\cK_H\otimes(\Qlb[1](\frac{1}{2}))^{\dimrel(\nu_Q)}\,\iso\, K_{Q,\psi}
\end{equation}
\end{Con} 

  We prove that $P$ and $Q$-models are compatible.  Namely, in Section~5 we define an open substack $^{\diamond \!}\Bun_{P\cap Q}\subset \Bun_{P\cap Q}$ and show that 
the restrictions of $K_{Q,\psi}$ and of $K_{P,\psi}$ to $^{\diamond \!}\Bun_{P\cap Q}$ are isomorphic up to a shift
(cf. Proposition~\ref{Pp_comparison_PQ}). Note that $P\cap Q$ contains a Borel subgroup of $H$.
This implies that the pointwise Euler-Poincare characteristics of $K_{P,\psi}$ (resp., of $K_{Q,\psi}$) are constant along the fibres of the projection $\nu_P: {^e\Bun_P}\to \Bun_H$ (resp., $\nu_Q: {^u\Bun_Q}\to\Bun_H$), cf. Proposition~\ref{Pp_Euler_char_pointwise}. 

\begin{Rem} Consider $\Four_{\cY_Q,\psi}^{-1}(\nu_Q^*\cK_H)$ over the whole of $\cY_Q$, we expect it to be the extension by zero from $\cZ_Q$.
\end{Rem}
 
\medskip\noindent
2.3.3 {\scshape $R$-modelÊ\ } Fix a 2-dimensional subspace $U_{0,2}\subset U_0$. Let $R\subset H$ be the parabolic subgroup preserving $U_{0,2}$, we call it the Heisenberg parabolic. 
The stack $\Bun_{R}$ classifies $V\in\Bun_H$ with an isotropic subbundle $U_2\subset V$, where $U_2\in\Bun_2$.
For $(U_2\subset V)\in\Bun_{R}$ write $V'=V_{-2}/U_2\in\Bun_{H_{n-2}}$, where $V_{-2}$ is the orthogonal complement of $U_2$ in $V$. We also need the stack $\Bun_{P_{n-2}}$ classifying $(U'\subset V')$, where $V'\in \Bun_{H_{n-2}}$ and $U'$ is an isotropic subbundle of rank $n-2$.

 Let $\cY_{R}$ be the stack classifying $(U_2\subset V)\in\Bun_R$ and a section $v_2: \wedge^2 U_2\to \Omega$. 
Let $f_{R}: \cY_{R}\to \Bun_R$ be the projection forgetting $v_2$. Write $j_{R}: \ocY_{R}\hook{} \cY_{R}$ for the open substack given by $v_2\ne 0$. 
 
 Let $\cX_{R}$ be the stack classifying $(U_2\subset V)\in\Bun_{R}$ and an upper modification $s_2: U_2\subset M$ equipped with $\det M\,\iso\,\Omega$. Here $M\in \Bun_2$ and $s_2$ is an inclusion of coherent $\cO_X$-modules.
Let 
$$
\pi_{R}: \cX_{R}\to \ocY_{R}
$$ 
be the map over $\Bun_{R}$ given by $v_2=\wedge^2 s_2$. The map $\pi_{R}$ is representable and proper. In Section~6.2.2 we define a natural map
\begin{equation}
\label{map_cXP^H_to_BuntG}
\tilde\rho_{R}: \cX_{R}\to \Bunt_{\GG_{2n-4}}
\end{equation}

 Let $^{\flat}(\Bun_2\times\Bun_{P_{n-2}})\subset \Bun_2\times\Bun_{P_{n-2}}$ be the open substack given by 
\begin{equation}
\label{conditions_flat}
\H^0(X, \Sym^2(U_2\otimes U'))=\H^0(X, \Omega\otimes \wedge^2 U')=\H^0(X, U_2\otimes U'^*\otimes\Omega)=0
\end{equation}
for $(U_2, U'\subset V')\in \Bun_2\times\Bun_{P_{n-2}}$. The projection is smooth
\begin{equation}
\label{projection_ww}
\id\times\nu_P: {^{\flat}(\Bun_2\times\Bun_{P_{n-2}})}\to \Bun_2\times\Bun_{H_{n-2}}
\end{equation}
 
 Let $^w(\Bun_2\times\Bun_{H_{n-2}})\subset \Bun_2\times\Bun_{H_{n-2}}$ be the intersection of the image of (\ref{projection_ww}) with the open substack given by the conditions
\begin{equation}
\label{cond_w_for_P^H} 
\H^0(X, \wedge^2 U_2)=\H^0(X, \Omega\otimes \wedge^2 U_2)=\H^0(X, \Omega\otimes U_2\otimes V')=0
\end{equation}
Informally, $(U_2, V')\in {^w(\Bun_2\times\Bun_{H_{n-2}})}$ if $U_2$ is sufficiently `negative' and $V'$ is `sufficiently stable' compared to $U_2$.

Let $^w\ocY_{R}$, $^w\cY_{R}$, $^w\Bun_{R}$ and so on denote the preimage of $^w(\Bun_2\times\Bun_{H_{n-2}})$ in the corresponding stack.
  
  Let $\nu_{R}: \Bun_{R}\to\Bun_H$ be the map induced by $R\hook{} H$. Its restriction to $^w\Bun_{R}$ is smooth.   
The following is proved in Section~6.3.   
\begin{Pp} 
\label{Pp_perversity_for_PH-model}
The complex 
\begin{equation}
\label{perv_sheaf_over_c_ocY_PH}
(\pi_{R})_!\tilde\rho_{R}^*\Aut^e_{\psi}\otimes(\Qlb[1](\frac{1}{2}))^{\dimrel(\tilde\rho_R)}
\end{equation}
is perverse over the open substack $^w\ocY_{R}\subset \ocY_{R}$. 
\end{Pp}
  
  Let $\cF_{R,\psi}$ be the intermediate extension of (\ref{perv_sheaf_over_c_ocY_PH}) under $j_{R}: {^w\ocY_{R}}\hook{} {^w\cY_{R}}$. For the projection $f_{R}: {^w\cY_{R}}\to {^w\Bun_{R}}$
we set
\begin{equation}
\label{complex_K_PH_psi}
K_{R,\psi}=(f_{R})_!\cF_{R,\psi} \in\D^{\prec}(^w\Bun_R)
\end{equation}
  
\begin{Con}
\label{Con3} The complex (\ref{complex_K_PH_psi}) is a perverse sheaf on $^w\Bun_{R}$, and there exists an isomorphism over $^w\Bun_{R}$
\begin{equation}
\label{iso_Th3_on_eBunP^H}
\nu_R^*\cK_H\otimes(\Qlb[1](\frac{1}{2}))^{\dimrel(\nu_R)}\,\iso\, K_{R,\psi}
\end{equation}
\end{Con}
  
A partial evidence for Conjecture~\ref{Con3} is provided in Section~6. Namely, in Section~6.3 we define an open substack $^{w\flat}\Bun_{P\cap R}\subset \Bun_{P\cap R}$ and show that the restrictions of $K_{P,\psi}$ and of $K_{R,\psi}$ to $^{w\flat}\Bun_{P\cap R}$ are isomorphic up to a shift (cf. Corollary~\ref{Cor_PR_compatibility}). Note that here $P\cap R$ contains a Borel subgroup of $H$. 

\medskip
\noindent
2.3.4 {\scshape Actual construction of $\cK_H$\  } We don't know in general if $\cK_H$ is nonzero at the generic point of each connected component of $\Bun_H$. 
For this reason it is not clear if the isomorphisms of Theorem~\ref{Th1} or Conjectures~\ref{Con2},\ref{Con3} characterize $\cK_H$ up to a unique isomorphism. 
  
  Our actual construction of $\cK_H$ is via the theta-lifting for the dual pair $(G_1, H)$. Write $\Aut_{G_1,H}$ for the complex on $\Bun_{G_1}\times\Bun_H$ introduced in (\cite{L2}, Definition~2). This is the kernel of the theta-lifting functor from $\Bun_{G_1}$ to $\Bun_H$. 
  
  Let $q_H: \Bun_{G_1}\times\Bun_H\to\Bun_H$ be the projection. Since $\Bun_{G_1}$ is not of finite type, the complex $q_{H !}\Aut_{G_1,H}$ does not make sense literally (at the level of functions for $k=\Fq$ the corresponding integral is also divergent). 
  
  For $a\in\ZZ$ write $_a\Bun_{G_1}\subset\Bun_{G_1}$ for the open substack classifying $M\in\Bun_{G_1}$ such that for any line bundle $L$ on $X$ with $\deg(L)\le a$ one has $\Hom(M, L)=0$. The stack $_a\Bun_{G_1}$ is of finite type.
Let $_aq: {_a\Bun_{G_1}}\times\Bun_H\to\Bun_H$ be the projection. Set
\begin{equation}
\label{complex_a_tildeK_on_BunH} 
_a\tilde K=(_aq_!)\Aut_{G_1,H}\in \D^{\prec}(\Bun_H)
\end{equation}

 For $a\in\ZZ$ write $_a\Bun_n\subset\Bun_n$ for the open substack classifying $U\in \Bun_n$ such that for any line bundle $L$ on $X$ of degree $\le a$ one has $\Hom(U, L)=0$. Write $^e_a\Bun_n$, $^e_a\Bun_P$ and so on for the preimage of $_a\Bun_n$ in the corresponding stack.
 
  First, using the compatibilities of $P$ and $Q$-models, we find an open substack of finite type $\cU_H\subset\Bun_H$ large enough such that the sheaf $\cK_H$ we are looking for should not vanish over $\cU_H^b$ for each $b\in \ZZ/2\ZZ$, where $\cU_H^b=\cU_H\cap \Bun_H^b$. We also find $a\in\ZZ$ small enough with the following property. Let 
$$
\tilde\cK_{\cU}={^p\H^0(_a\tilde K)\mid_{\cU_H}}
$$  
We prove that for each $b\in \ZZ/2\ZZ$ there is a unique irreducible subquotient $\cK_{\cU, b}$ in $\tilde\cK_{\cU}\mid_{\cU_H^b}$ with the following properties. Let $\cK_H$ be the intermediate extension of $\cK_{\cU,0}\oplus\cK_{\cU,1}$ under $\cU_H\hook{}\Bun_H$. Then Theorem~\ref{Th1} holds for this $\cK_H$. Moreover, all the other irreducible subquotients of $\tilde\cK_{\cU}$ are shifted almost constant local systems on $\Bun_H$ (cf. Proposition~\ref{Pp_K^d_coincide}). The notion of an almost constant local system is introduced in Appendix~A.
  
  Moreover, for any $i\ne 0$ and any irreducible subquotient $F$ of $^p\H^i(_a\tilde K)\mid_{\cU_H}$, let $\bar F$ be the intermediate extension of $F$ under $\cU_H\hook{}\Bun_H$. Then $\bar F$ is a shifted almost constant local system on $\Bun_H$ (cf. Remark~\ref{Rem_about_H^i_for_i_not_zero}). 
  
\begin{Rem} Most of our results hold also for $k=\Fq$, the corresponding changes are left to the reader. Part 2) of Proposition~\ref{Pp_small_map} holds over $k$ algebraically closed only as it uses the decomposition theorem (\cite{BBD}, Corollary~5.4.6).
\end{Rem} 

\medskip\noindent
2.4 In Section~8 we propose one more conjectural construction of the perverse sheaf $\cK_H$ via residues of geometric Eisenstein series. We do not completely check that it indeed produces $\cK_H$ except in cases $g=0$ and $g=1$. In the present paper, the construction via Eisenstein series for us is rather a way to calculate $\cK_H$ explicitly in these particular cases.

\medskip\noindent
2.4.1 {\scshape Case $g=0$} \  Recall that $\Bun_H$ admits the Shatz stratification indexed by dominant coweights $\Lambda^+_H$ of $H$ (cf. Section~8.7). Write $Shatz^{\lambda}$ for the Shatz stratum corresponding to $\lambda\in\Lambda^+_H$. For $b\in\ZZ/2\ZZ$ write $OSh^b$ for the open Shatz stratum in $\Bun_H^b$. For $b=0$ (resp., $b=1$) let $\lambda=(1,1,0,\ldots, 0)$ (resp., $\lambda=(1,1,1,0,\ldots,0)$). We show that for each $b\in\ZZ/2\ZZ$ the stack $\Bun_H^b-OSh^b$ is irreducible, and its open Shatz stratum is $Shatz^{\lambda}$. Call $Shatz^{\lambda}$ the \select{subregular Shatz stratum} in $\Bun_H^b$ by analogy with subregular unipotent orbits.

\begin{Pp} 
\label{Pp_answer_for_g=0}
Assume $g=0$.
For each $b\in\ZZ/2\ZZ$ one has $\cK_H\,\iso\, \IC(Shatz^{\lambda})$ over $\Bun_H^b$, where $Shatz^{\lambda}$ is the subregular Shatz stratum in $\Bun_H^b$.
\end{Pp}

\smallskip\noindent
2.4.2 {\scshape Case $g=1$} \  In Section~8.8 we introduce an open substack $\cW^0_H\subset \Bun^0_H$ with the following property. Let $T\subset H$ be the standard maximal torus, $W$ be the Weyl group of $(T,H)$. Let $\nu_T^0: \Bun_T^0\to \Bun_H^0$ be the natural map. Then over $\cW^0_H$ the map $\nu_T^0$ is a Galois covering with Galois group $W$. For an irreducible representation $\sigma$ of $W$ write $\cL_{\sigma}$ for the perverse sheaf on $\Bun_H$, the intermediate extension under $\cW^0_H\hook{}\Bun_H^0$ of the isotypic component of $(\nu_T^0)_!\Qlb\mid_{\cW^0_H}$ corresponding to $\sigma$. Then $\cL_{\sigma}$ is an irreducible perverse sheaf.

 Recall that $\H^1(X, \mu_2)\,\iso\, (\ZZ/2\ZZ)^2$. Let $\tau^1: \Spec k\to \Bun_{H_2}^1$ be the map given by 
$$
V=\sum_{\alpha\in \H^1(X, \mu_2)} \cA_{\alpha},
$$ 
where the symmetric form on $V$ is the orthogonal sum of the canonical forms $\cA_{\alpha}^2\,\iso\,\cO_X$. Here $\cA_{\alpha}$ denotes the line bundle obtained from $\alpha$ via extension of scalars $\mu_2\subset \Gm$. Note that $V$ is semistable, it admits no isotropic line (or rank 2) subbundles of degree zero. The map $\tau^1$ is \'etale, write $\cW^1_{H_2}$ for its image. Actually, $\cW^1_{H_2}$ is the classifying stack of $(\ZZ/2\ZZ)^3$.

  Let $f^1: \Bun_{H_2}^1\times\Bun_{H_{n-2}}^0\to \Bun_H^1$ be the map sending $(V,V')$ to $V\oplus V'$, the symmetric form being the orthogonal direct sum of forms for $V, V'$. The restriction of $f^1$ to $\cW^1_{H_2}\times \cW^0_{H_{n-2}}$ is \'etale, write $\cW^1_H$ for the image of $\cW^1_{H_2}\times \cW^0_{H_{n-2}}$ under $f^1$. 
  
  Recall that the symmetric group $S_n$ is naturally a quotient of $W$. The standard representation of $S_n$ in $\Qlb^n$ by permutations decomposes as $\Qlb\oplus \sigma_n$, where $\sigma_n$ is an irreducible $(n-1)$-dimensional representation. View $\sigma_n$ as a representation of $W$. Note that $\sigma_n$ is also a subrepresentation of the standard representation of $\check{H}\,\iso\,\SO_{2n}$ restricted to its Weyl group.
  
\begin{Pp} 
\label{Pp_case_g=1}
Assume $g=1$. Over $\Bun_H^0$ the perverse sheaf $\cK_H$ is isomorphic to $\cL_{\sigma_n}$. The restriction of $\cK_H$ under the etale map 
$$
f^1: \cW^1_{H_2}\times \cW^0_{H_{n-2}}\to \Bun_H^1
$$ 
is isomorphic to the irreducible perverse sheaf $\Qlb\boxtimes \cL_{\sigma_{n-2}}$. In particular, over $\cW^0_H$ (resp., $\cW^1_H$) the perverse sheaf $\cK_H$ is a (shifted) local system of rank $n-1$ (resp., $n-3$). For $n=4$ this local system over $\cW^1_H$ is of order two. 
\end{Pp}

  Since $\cK_H$ does not vanish at each generic point of $\Bun_H$, the isomorphism (\ref{iso_Th1_on_eBunP}) of Theorem~\ref{Th1} in the case $g=1$ determines $\cK_H$ up to a unique isomorphism over each connected component of $\Bun_H$ (cf. Proposition~\ref{Pp_fibres_of_nu_P_connected}). 
  
\begin{Rem} The isomorphism (\ref{iso_Th1_on_eBunP}) does not hold over the whole of $^e\Bun_P$. Indeed, Proposition~\ref{Pp_case_g=1} shows that for $g=1$ it does not hold over $^e\Bun_P^0$ and $0\notin Z(e,P)$.   
\end{Rem}

\medskip\noindent
2.5 In Section~9 we propose Conjecture~\ref{Con_general_descent} generalizing our construction for other simple algebraic groups, which admit a maximal parabolic subgroup with an abelian unipotent radical. Conjecture~\ref{Con_general_descent} should lead, in particular, to the geometric analogs of the minimal representations for $E_6$ and $E_7$. 

\bigskip

\centerline{\scshape 3. Extended theta sheaf} 
 
\bigskip\noindent 
3.1.1  Keep notations of Section~2.2. Let $V$ be a $k$-vector space of dimension $d$. Recall that $\Sym^2(V)$ is the quotient of $V\otimes V$ by the subspace generated by the vectors $v_1\otimes v_2-v_2\otimes v_1$, $v_i\in V$. Its dual identifies canonically with the space $\ST^2(V^*)$ of symmetric tensors in $V^*\otimes V^*$. View $\ST^2(V^*)$ as the space of symmetric bilinear forms on $V$. 

 Let $p_V: V\to\Sym^2 V$ be the map sending $v$ to $v\otimes v$. It is finite, a Galois $S_2$-covering over its image $\Im p_V$ with zero removed. Consider the diagram
$$
\A^1\getsup{\ev_V} V^*\times V\toup{\id\times p_V} V^*\times\Sym^2(V),
$$
where $\ev_V$ is the natural pairing. Write $\Four_{2,\psi}: \D^b(V^*\times\Sym^2(V))\,\iso\, \D^b(V^*\times \ST^2(V^*))$ for the Fourier transform along $\Sym^2 V$. Set
\begin{equation}
\label{def_S_psi_e}
S^e_{\psi}=\Four_{2,\psi}((\id\times p_V)_!\ev_V^*\cL_{\psi})[2d](d),
\end{equation}
where $e$ stands for `extended'. This is a perverse sheaf on $V^*\times \ST^2(V^*)$, and $\DD(S^e_{\psi})\,\iso\, S^e_{\psi^{-1}}$ canonically.

 Let $\sigma$ be the automorphism of $V$ of multiplication by $-1$. Then $\sigma^*\ev_V^*\cL_{\psi}\,\iso\, \ev_V^*\cL_{\psi}^*$. Since the characteristic is not 2, $(\ev_V^*\cL_{\psi})^{\otimes 2}$ is nontrivial for $d>0$, and  
$$
(p_V\times\id)_!\ev_V^*\cL_{\psi}[2d]
$$
is an irreducible perverse sheaf. Thus, $S^e_{\psi}$ is also irreducible. 
 
 Let $i: \Spec k\hook{} V^*$ be the zero section. For the map $i\times\id: \ST^2(V^*)\hook{} V^*\times\ST^2(V^*)$ set 
$$
S_{\psi}=(i\times\id)^*S^e_{\psi}[-d](-\frac{d}{2})
$$  
Then $S_{\psi}$ identifies canonically with the perverse sheaf introduced in (\cite{L1}, Section~4.1). In this sense $S^e_{\psi}$ extends $S_{\psi}$.

 Note that $S^e_{\psi}$ is $\GL(V)$-equivariant. 
Let $\rho_V: V^*\times\ST^2(V^*)\to (V^*\times\ST^2(V^*))/\GL(V)$ denote the stack quotient. Denote by $\SSS^e_{\psi}$ an irreducible perverse sheaf on $(V^*\times\ST^2(V^*))/\GL(V)$ equipped with an isomorphism $\rho_V^*\SSS^e_{\psi}[d^2](\frac{d^2}{2})\,\iso\, S^e_{\psi}$. The perverse sheaf $\SSS^e_{\psi}$ is defined up to a unique isomorphism.

\begin{Rem} For $b: V\to V^*$ with $b^*=b$ let $\beta_b: V\to\A^1$ be the map $v\mapsto \<v, bv\>$. One has a usual ambiguity in identifying $\ST^2(V^*)$ with $\Sym^2 (V^*)$, namely, $b$ goes to $\beta_b$ or $\frac{1}{2}\beta_b$. 
\end{Rem}
 
\medskip\noindent
3.1.2 Let $Q_i(V)\subset \ST^2(V^*)$ be the locally closed subscheme of $b: V\to V^*$ symmetric such that $\dim\Ker b=i$. Let $Q'_i(V)\subset V^*\times Q_i(V)$ be the closed subscheme of pairs $(v^*, b)$ such that $v^*\in (V/V_0)^*$, where $V_0=\Ker b$. 

\begin{Pp} The $*$-restriction of $S^e_{\psi}$ to $V^*\times Q_i(V)$ is the extension by zero from $Q'_i(V)$ of a $\GL(V)$-equivariant rank one local system placed in usual cohomological degree $-d+i-d(d+1)/2$.
\end{Pp}
\begin{Prf} 
For $b: V\to V^*$ with $b^*=b$ and $v^*\in V^*$ 
let $\beta_{b,v^*}: V\to\A^1$ be the map $v\mapsto \<v, bv\>+\<v,v^*\>$. Let $V_0=\Ker b$ and $p_0: V\to V/V_0$ be the projection. Then $(p_0)_!\beta_{b,v^*}^*\cL_{\psi}$ will vanish unless $v^*\in (V/V_0)^*$. 

 We are reduced to the case $V_0=0$. In this case, over the algebraic closure $\bar k$, in a suitable affine coordinates of $V$ the quadratic form $v\mapsto \beta_{b,v^*}(v)$ writes as $(x_1,\ldots, x_d)\mapsto x_1^2+\ldots+x_d^2+a$ for some $a\in \bar k$. Our assertion follows now from (\cite{L1}, Lemma~3). 
\end{Prf} 

\medskip

 Let $\phi_V: V^*\times Q_0(V)\to\A^1$ be the map $(v^*,b)\mapsto \frac{1}{4}\<v^*, b^{-1}v^*\>$. Here we have viewed $b\in Q_0(V)$ as a symmetric isomorphism $b: V\,\iso\, V^*$. Let $q_V$ denote the composition $V^*\times Q_0(V)\to Q_0(V)\hook{} \ST^2(V^*)$, where the first map is the projection.
 
\begin{Pp} There is a canonical isomorphism over $V^*\times Q_0(V)$
$$
\phi_V^*\cL_{\psi}\otimes S^e_{\psi}\,\iso\, q_V^* S_{\psi}[d](\frac{d}{2})
$$
\end{Pp}
\begin{Prf}
Let $(v^*,b)\in V^*\times Q_0(V)$ and $v\in V$.
The change of variables $w=v+b^{-1}v^*/2$ gives
$\<v, bv\>+\<v, v^*\>=\<w, bw\>-\frac{1}{4}\<b^{-1}v^*, v^*\>$. Our assertion follows.
\end{Prf}
 
\medskip\noindent 
3.1.3 We will need a relative version of the above construction. Let $S$ be a smooth stack locally of finite type. Let $\cV\to S$ be a vector bundle of rank $d$. For the diagram as above
$$ 
\A^1\getsup{\ev_{\cV}} \cV^*\times_S \cV\toup{\id\times p_{\cV}} \cV^*\times_S\Sym^2(\cV)
$$
write $\Four_{2,\psi}: \D^b(\cV^*\times_S \Sym^2 \cV)\,\iso\, \D^b(\cV^*\times_S \ST^2(\cV^*))$ for the Fourier transform along $\Sym^2 \cV$. By abuse of notation, we also denote
\begin{equation}
\label{def_S_psi_e_relative}
S^e_{\psi}=\Four_{2,\psi}((\id\times p_{\cV})_!\ev_{\cV}^*\cL_{\psi})\otimes(\Qlb[1](\frac{1}{2}))^{2d+\dim S}
\end{equation}
Note that $\cV$ yields a morphism of stacks 
$$
\alpha_{\cV}: \cV^*\times_S \ST^2(\cV^*)\to (V^*\times\ST^2(V^*))/\GL(V)
$$ 
and (\ref{def_S_psi_e_relative}) is canonically isomorphic to $\alpha_{\cV}^*\SSS^e_{\psi}\otimes(\Qlb[1](\frac{1}{2}))^{\dimrel(\alpha_{\cV})}$.
 
\medskip\noindent
3.2 Write $G_n$ for the group scheme on $X$ of automorphisms of $M_0=\cO^n\oplus\Omega^n$ preserving the natural symplectic form $\wedge^2 M_0\to\Omega$. Write $H_n=M_0\oplus\Omega$ for the corresponding Heisenberg group with operation 
$$
(m_1,\omega_1)(m_2,\omega_2)=(m_1+m_2, \omega_1+\omega_1+\frac{1}{2}\<m_1, m_2\>)
$$
The group scheme $G_n$ acts on $H_n$ by group automorphisms so that $g\in G_n$ sends $(m,\omega)\in H_n$ to $(gm,\omega)$. Write an element of $\GG_n=\GG_n\rtimes H_n$ as $(g, (m,\omega))$ with $g\in G_n$, $(m,\omega)\in H_n$ then the product in $\GG_n$ is given by
$$
(g_1, (m_1,\omega_1))(g_2, (m_2,\omega_2))=(g_1g_2, (g_2^{-1}m_1, \omega_1)(m_2,\omega_2))
$$
Recall the subgroup $\PP_n\subset \GG_n$ (cf. Section~2.2). The stack $\Bun_{P_n}$ classifies $\cL\in\Bun_n$ and an exact sequence $0\to\Sym^2\cL\to ?\to\Omega\to 0$ on $X$, it gives rise to an exact sequence $0\to\cL\to M\to \cL^*\otimes\Omega\to 0$ with $M\in\Bun_{G_n}$. 

\medskip
\noindent
3.3 For Section~3.3 we assume $k=\Fq$. 
Write $L=\cO^n$, this is a lagrangian subbundle in $M_0=L\oplus L^*\otimes\Omega$. 
 
 Write $\cS(L^*\otimes\Omega(\AA))$ for the Schwarz space of locally constant $\Qlb$-valued functions with compact support on $L^*\otimes\Omega(\AA)$. This is a model of the Weil representation of $\H_n(\AA)$, in which the metaplectic extension (\ref{seq_metaplectic_classical_def}) naturally splits over $\PP_n(\AA)$. The purpose of this section is to give an explicit formula for the restriction 
$$
\phi_{\PP}: \PP_n(F)\backslash \PP_n(\AA)/\PP_n(\cO)\to\Qlb
$$
of (\ref{autom_form_phi}). Recall the character $\chi: \Omega(\AA)/\Omega(\cO)\to\Qlb^*$ given by (\ref{character_chi_for_Fq}). The action of $\PP_n(\AA)$ in $\cS(L^*\otimes\Omega(\AA))$ is given by the following formulas. 
 
 For $l_1\in L(\AA), l_1^*\in L^*\otimes\Omega(\AA),
\omega\in \Omega(\AA)$ and $f\in \cS(L^*\otimes\Omega(\AA))$ the action of $(l_1+l_1^*, \omega)\in H_n(\AA)$ on $f$ is the function
$$
l^*\in L^*\otimes\Omega(\AA)\mapsto 
\chi(\omega+\<l^*, l_1\>+\frac{1}{2}\<l_1^*, l_1\>)f(l^*+l_1^*)
$$ 
 
 Write $\AA^*$ for the ideles of $X$. For $a\in\AA^*$ write $\mid a\mid\in\Qlb^*$ for the absolute value of $a$. For a vector bundle $\cW$ on $X$ and $g\in\GL(\cW)(\AA)$ write $\mid g\mid=\mid \det g\mid$.
 
 Let $a\in\GL(L)(\AA)$, $b\in \Hom(L^*\otimes\Omega, L)(\AA)$. The action of 
\begin{equation}
\label{g_formula}
g=\left(
\begin{array}{cc}
a & b\\
0 & {a^*}^{-1}
\end{array}
\right)\in P_n(\AA)
\end{equation}
on $f$ is the function 
$$
l^*\in L^*\otimes\Omega(\AA)\mapsto \;
\mid a\mid^{\frac{1}{2}}
\chi(\frac{1}{2}\<a^*l^*, b^*l^*\>)f(a^*l^*)
$$
Let $g$ be given by (\ref{g_formula}) and $m=l_1+l_1^*\in M_0(\AA)$, $\omega\in\Omega(\AA)$. It follows that the action of $(g, (m,\omega))\in \PP_n(\AA)$ on $f$ is the function
$$
l^*\in L^*\otimes\Omega(\AA)\mapsto \;
\mid a\mid^{\frac{1}{2}}
\chi(\frac{1}{2}\<a^*l^*, b^*l^*\>)\chi(\omega+\<a^*l^*, l_1\>+\frac{1}{2}\<l_1^*, l_1\>)f(a^*l^*+l_1^*)
$$  

The theta-functional $\Theta: \cS(L^*\otimes\Omega(\AA))\to\Qlb$ (cf. Section~2.2) sends $f$ to 
$$
\sum_{l^*\in L^*\otimes\Omega(F)} f(l^*)
$$ 
Let $\phi_0$ be the characteristic function of $L^*\otimes\Omega(\cO)$, this is a unique up to a multiple $H_n(\cO)$-invariant vector in $\cS(L^*\otimes\Omega(\AA))$. So, the value of $\phi_{\PP}$ on $(g, (m,\omega))\in\PP_n(\AA)$ is
\begin{equation}
\label{value_first}
\sum_{l^*\in L^*\otimes\Omega(F)} \mid a\mid^{\frac{1}{2}}
\chi(\frac{1}{2}\<a^*l^*, b^*l^*\>)\chi(\omega+\<a^*l^*, l_1\>+\frac{1}{2}\<l_1^*, l_1\>)\phi_0(a^*l^*+l_1^*)
\end{equation}

One has a canonical bijection 
\begin{equation}
\label{bij_for_Bunn(k)}
\{\cL\in\Bun_n(k), \alpha: \cL(F)\,\iso\, L(F)\}\;\iso\; \GL(L)(\AA)/\GL(L)(\cO)
\end{equation}
One also has a canonical bijection 
$$
\Bun_{\PP_n}(k)\,\iso\, \PP_n(F)\backslash \PP_n(\AA)/\PP_n(\cO),
$$
where $\Bun_{\PP_n}(k)$ is the set of isomorphism classes of $\PP_n$-torsors on $X$. Recall that $\Bun_{\PP_n}$ is the stack classifying pairs of exact sequences (\ref{seq_cL_by_Omega}) and (\ref{seq_Omega_by_Sym2barcL}) on $X$ (cf. Section~2.2).

 Consider a point $\cF_{\PP_n}\in\Bun_{\PP_n}$ given by this pair of exact sequences and corresponding to the double class of $(g, (m,\omega))\in \PP_n(\AA)$. We assume that $g$ is given by (\ref{g_formula}). Let $\cL\in\Bun_n$ together with a trivialization $\alpha: \cL(F)\,\iso\, L(F)$ correspond to $a\GL(L)(\cO)\in \GL(L)(\AA)/\GL(L)(\cO)$ via (\ref{bij_for_Bunn(k)}). 
 
 For each closed point $x\in X$ write $F_x$ for the completion of $F$ at $x$, write $\cO_x\subset F_x$ for the completion of $\cO_X$ at $x$. For $x\in X$ we have a diagram
$$
\begin{array}{ccc}
L^*(F_x) & \toup{\alpha^*} & \cL^*(F_x)\\
\cup && \cup\\
a^{*-1}L^*(\cO_x) & \toup{\alpha^*}& \cL^*(\cO_x),
\end{array}
$$
where the horizonal arrows are isomorphisms. 

  Recall that $\H^1(X,\cL)\,\iso\,\cL(\AA)/(\cL(F)+\cL(\cO))$ canonically.
In particular, 
\begin{equation}
\label{H1_description_cLstar}
\H^1(X, \cL^*\otimes\Omega)\,\iso\, L^*\otimes\Omega(\AA)/(a^{*-1}(L^*\otimes\Omega)(\cO)+(L^*\otimes\Omega)(F))
\end{equation}
The extension (\ref{seq_cL_by_Omega}) is given by the image of $a^{*-1}l_1^*\in L^*\otimes\Omega(\AA)$ in (\ref{H1_description_cLstar}). Clearly, (\ref{value_first}) vanishes unless there is $l^*\in L^*\otimes\Omega(F)$ with $a^*l^*+l_1^*\in L^*\otimes\Omega(\cO)$. That is, the image of $a^{*-1}l_1^*$ in $\H^1(X, \cL^*\otimes\Omega)$ vanishes and (\ref{seq_cL_by_Omega}) splits. 
So, $\phi_{\PP}(\cF_{\PP_n})=0$ unless (\ref{seq_cL_by_Omega}) splits. 

 Now it is convenient to assume that $l_1^*=0$. Fix a splitting $\bar \cL\,\iso\, \Omega\oplus\cL$ of (\ref{seq_cL_by_Omega}). Since $\Omega^{-1}\otimes\Sym^2(\bar\cL)\,\iso\, \Omega\oplus (\Omega^{-1}\otimes\Sym^2\cL) \oplus\cL$, the datum of $\cF_{\PP_n}$ becomes a datum of three exact sequences 
(\ref{seq_cO_by_Omega}) given by $\xi\in \H^1(X,\Omega)$, 
\begin{equation}
\label{seq_Omega_by_Sym2cL}
0\to\Sym^2 \cL\to ?\to\Omega\to 0
\end{equation}
given by $\gamma\in \H^1(X, \Omega^{-1}\otimes\Sym^2\cL)$ and 
\begin{equation}
\label{seq_cO_by_cL}
0\to \cL\to ?\to \cO\to 0
\end{equation}
given by $\delta\in\H^1(X,\cL)$. Note that $\delta$ corresponds to the image of $l_1$ in $\H^1(X,\cL)$, and $\xi$ is the image of $\omega$ in $\Omega(\AA)/(\Omega(F)+\Omega(\cO))$.
 
Note that 
$$
\{l^*\in L^*\otimes\Omega(F)\mid a^*l^*\in L^*\otimes\Omega(\cO)\} \;\toup{\alpha^*}\; \H^0(X,\cL^*\otimes\Omega)
$$
is a bijection, so
\begin{equation}
\label{value_almost_geom}
\phi_{\PP}(\cF_{\PP_n})=\sum_{s\in \H^0(X, \cL^*\otimes\Omega)} \mid a\mid^{\frac{1}{2}}\psi(\<s\otimes s, \gamma\>+\<s, \delta\>+\omega)
\end{equation}

\medskip\noindent
3.4.1 {\scshape Geometrization}  
Let $f_{\PP}: \Bun_{\PP_n}\to\Bun_n$ be the map sending a pair (\ref{seq_cL_by_Omega}) and (\ref{seq_Omega_by_Sym2barcL}) to $\cL$. Write $_c\Bun_n\subset\Bun_n$ for the open substack of $\cL\in\Bun_n$ given by $\H^0(X,\cL)=0$. 
Write $_c\Bun_{\PP_n}\subset\Bun_{\PP_n}$ for the preimage of $_c\Bun_n$ under $f_{\PP}$.

 Let $\nu_{\cY}: \cY\to {_c\Bun_{\PP_n}}$ be the stack classifying a point of $_c\Bun_{\PP_n}$ as above together with a splitting of (\ref{seq_cL_by_Omega}). Note that $\nu_{\cY}$ is a torsor under a vector bundle on $_c\Bun_{\PP_n}$ with fibre $\Hom(\cL,\Omega)$.

 The stack $\cY$ can be seen as the stack classifying $\cL\in {_c\Bun_n}$ and exact sequences 
(\ref{seq_cO_by_Omega}) given by $\xi\in\H^1(X,\Omega)$, (\ref{seq_Omega_by_Sym2cL}) given by $\gamma\in \H^1(X, \Omega^{-1}\otimes\Sym^2\cL)$
and (\ref{seq_cO_by_cL}) given by $\delta\in\H^1(X,\cL)$.

 Let $p_{\cX}: \cX\to\cY$ be the stack over $\cY$ classifying the same objects as $\cY$ together with $s\in\Hom(\cL,\Omega)$. Let $\ev_{\cX}: \cX\to\A^1$ be the map sending the above point to $\xi+\<s\otimes s, \gamma\>+\<s, \delta\>$. It is understood that $s\otimes s\in \Hom(\Sym^2 L, \Omega^2)$. Set
$$ 
K_{\cY,\psi}=p_{\cX !}\ev_{\cX}^*\cL_{\psi}\otimes(\Qlb[1](\frac{1}{2})^{\dim\cX},
$$ 
where $\dim\cX$ is the dimension of the corresponding connected component of $\cX$. Then $K_{\cY,\psi}$ is a geometric analog of (\ref{value_almost_geom}).

\medskip\noindent
3.4.2  Let $\cV\to {_c\Bun_n}$ be the vector bundle 
with fibre $\Hom(\cL,\Omega)$ at $\cL\in{_c\Bun_n}$.
The dual vector bundle $\cV^*\to{_c\Bun_n}$ classifies $\cL\in{_c\Bun_n}$ and an extension (\ref{seq_cO_by_cL}) on $X$. Set $\cW=\cV\times_{_c\Bun_n}\cV^*$. Let $\cW_2\to \cV^*$ be the stack classifying a point of $\cV^*$ together with an element of $\Hom(\Sym^2\cL, \Omega^2)$. Let $\cW_2^*\to \cV^*$ be the stack classifying a point of $\cV^*$ as above together with an exact sequence (\ref{seq_Omega_by_Sym2cL}) on $X$. Write $\Four_{\cW, \psi}:\D(\cW_2)\to\D(\cW_2^*)$ for the Fourier transform over $\cV^*$.

 Let $p_{2,\cW}: \cW\to\cW_2$ be the map over $\cV^*$ sending $s\in\Hom(\cL,\Omega)$ to $s\otimes s\in \Hom(\Sym^2\cL, \Omega^2)$. The map $p_{2,\cW}$ is finite, an $S_2$-covering over the image of $\Im p_{2,\cW}$ with removed zero section. Let $\ev_{\cW}: \cW\to \A^1$ be the natural pairing between $\cV$ and $\cV^*$. Then 
$\cY\,\iso\, \cW_2^*\times\Bun_{\Omega}$ naturally. 
By definition, 
$$
K_{\cY,\psi}\,\iso\, \Four_{\cW,\psi}((p_{2,\cW})_!\ev_{\cW}^*\cL_{\psi})\boxtimes\ev_{\Omega}^*\cL_{\psi} 
\otimes(\Qlb[1](\frac{1}{2}))^{\dim\cW+\dim\Bun_{\Omega}}
$$
As in Section~3.1.1, one shows now that $K_{\cY,\psi}$ is a perverse sheaf irreducible over each connected component of $\cY$, and $\DD(K_{\cY,\psi})\,\iso\, K_{\cY,{\psi}^{-1}}$ canonically.

  There is a natural map
$$
f_{\cW}: \cW_2^*\to \cV^*\times_{_c\Bun_n} \ST^2(\cV^*)
$$
defined as follows. The transpose to the linear map $\Sym^2 \H^0(X,\cL^*\otimes\Omega)\to \Hom(\Sym^2\cL,\Omega^2)$ is a map $\H^1(X, \Omega^{-1}\otimes\Sym^2\cL)\to \ST^2(\H^1(X,\cL))$ denoted $\gamma\mapsto\bar\gamma$. Then $f_{\cW}$ sends $(\cL,\gamma,\delta)$ to $(\cL, \bar\gamma, \delta)$. For the perverse sheaf $S^e_{\psi}$ on $\cV^*\times_{_c\Bun_n} \ST^2(\cV^*)$ defined in Section~3.1.3 one gets an isomorphism
$$
K_{\cY,\psi}\,\iso\, (f_{\cW}^*S^e_{\psi} \boxtimes\ev_{\Omega}^*\cL_{\psi})
\otimes(\Qlb[1](\frac{1}{2}))^{\dimrel(f_{\cW})+\dim\Bun_{\Omega}}
$$

\begin{Pp} There is a canonical isomorphism over $\cY$
$$
\nu^*_{\cY}K_{\PP_n, \psi}\otimes(\Qlb[1](\frac{1}{2}))^{\dimrel(\nu_{\cY})}\,\iso\, K_{\cY,\psi}
$$
\end{Pp}
\begin{Prf}
Let $_c\cT_n$ be the preimage of $_c\Bun_n$ in $\cT_n$, and similarly for $_c\cZ_{\cT_n}$, $_c\cZ_{2, \cT_n}$. Let $\cV_2\to {_c\Bun_n}$ be the stack classifying $\cL\in{_c\Bun_n}$ and a section 
$$
\bar s: \Sym^2(\cL\oplus\Omega)\to\Omega^2
$$ 
Let $\bar h_{\cT}: \cV\to \cV_2$ be the morphism over $_c\Bun_n$ sending $s_1\in \Hom(\cL,\Omega)$ to $\bar s=s\otimes s$ with $s=(s_1,\id): \cL\oplus\Omega\to\Omega$. Let $\nu_{\cT}: {_c\Bun_n}\to {_c\cT_n}$ be the map sending $\cL$ to $\bar\cL=\cL\oplus\Omega$. After the base change by $\nu_{\cT}$, the diagram
$$
\begin{array}{ccccc}
_c\cZ_{\cT_n} & \to & {_c\cT_n} & \gets & _c\Bun_{\PP_n}\\
& \searrow\lefteqn{\scriptstyle h_{\cT}} & \uparrow\\
&& _c\cZ_{2, \cT_n}
\end{array}
$$
identifies with the diagram
$$
\begin{array}{ccccc}
\cV & \to & _c\Bun_n & \gets & \cY\\
& \searrow\lefteqn{\scriptstyle \bar h_{\cT}} & \uparrow\\
&& \cV_2
\end{array}
$$
We have $\cV\times_{_c\Bun_n}\cY\,\iso\,\cX$ naturally. 
The stacks $\cV_2$ and $\cY$ are dual (generalized) vector bundles over $_c\Bun_n$, write $\ev_{\cV\cY}: \cV_2\times_{_c\Bun_n}\cY\to\A^1$ for the natural pairing. The diagram commutes
$$
\begin{array}{ccc}
\cV\times_{_c\Bun_n}\cY &\iso & \cX\\
\downarrow\lefteqn{\scriptstyle \bar h_{\cT}\times\id} && \downarrow\lefteqn{\scriptstyle \ev_{\cX}}\\
\cV_2\times_{_c\Bun_n}\cY &\toup{\ev_{\cV\cY}}&\A^1
\end{array}
$$
Our assertion follows.
\end{Prf}

\medskip
\noindent
3.4.3 Let $\nu_{\PP}: \Bun_{\PP_n}\to\Bun_{\GG_n}$ be the morphism induced by the inclusion $\PP_n\to \GG_n$. We lift it to a morphism $\tilde\nu_{\PP}: \Bun_{\PP_n}\to\Bunt_{\GG_n}$ sending a point (\ref{seq_cL_by_Omega}) and (\ref{seq_Omega_by_Sym2barcL}) of $\Bun_{\PP_n}$ to $(\Omega\hook{v} M_1, \cB)$. Here 
$\cB=\det\RG(X, \cL^*\otimes\Omega)$ is equipped with the $\ZZ/2\ZZ$-graded isomorphism
$$
\cB^2\,\iso\, \det\RG(X,M)
$$
given by the exact sequence $0\to\cL\to M\to \cL^*\otimes\Omega\to 0$. We have denoted here $M=L_{-1}/\Omega$, where $L_{-1}$ is the orthogonal complement of $\Omega$ in $M_1$.
  
  Recall the open substack $^0\Bun_n\subset {_c\Bun_n}$ introduced in Section~2.2. The restriction $\nu_{\PP}: {^0\Bun_{\PP_n}}\to\Bun_{\GG_n}$ of $\nu_{\PP}$ is smooth, hence $\tilde\nu_{\PP}: {^0\Bun_{\PP_n}}\to\Bunt_{\GG_n}$ is also smooth. 
      
  Recall that $_0\Bun_{\GG_n}$ is the preimage of $_0\Bun_{G_n}$ under $\rho_{\GG}: \Bun_{\GG_n}\to\Bun_{G_n}$. For a point $(\Omega\subset L_{-1}\subset M_1)$ of $_0\Bun_{\GG_n}$ the exact sequence $0\to \Omega\to L_{-1}\to M\to 0$ splits canonically, this yields an isomorphism $_0\Bun_{\GG}\,\iso\, {_0\Bun_G}\times\Bun_{\Omega}$ sending the above point to $M=L_{-1}/\Omega$ and (\ref{seq_cO_by_Omega}), which is the push-forward of $0\to L_{-1}\to M_1\to\cO\to 0$ by $L_{-1}\to\Omega$. This in turn gives the isomorphism (\ref{iso_0_Bunt_GGn_is}). 
  
\medskip\noindent
\begin{Prf}\select{of Proposition~\ref{Pp_explicit_formula_Aute}}

\medskip\noindent
Write $^0\cY$ for the preimage of $^0\Bun_{\PP_n}$ under $\nu_{\cY}: \cY\to {_c\Bun_{\PP_n}}$. Then $\nu_{\cY}: {^0\cY}\to {^0\Bun_{\PP_n}}$ is a torsor under the vector bundle $\cV\times_{_c\Bun_n}{^0\Bun_{\PP_n}}$.
Since both sides of (\ref{iso_explicit_Aute}) are perverse, it suffices to establish an isomorphism over $^0\cY$
\begin{equation}
\label{iso_main_over_0_cY}
\nu_{\cY}^*\tilde\nu_{\PP}^*\Aut^e_{\psi}\otimes(\Qlb[1](\frac{1}{2}))^{\dimrel(\nu_{\cY})+\dimrel(\nu_{\PP})}\,\iso\, K_{\cY,\psi}
\end{equation}
Write $^0_0\cY$ for the preimage of $_0\Bunt_{\GG_n}$ under $\tilde\nu_{\PP}\comp\nu_{\cY}: {^0\cY}\to\Bunt_{\GG_n}$. Since $K_{\cY,\psi}$ is irreducible over each connected component of $\cY$, and $\Aut^e_{\psi}$ is the intermediate extension from $_0\Bunt_{\GG_n}$, it suffices to establish the isomorphism (\ref{iso_main_over_0_cY}) over $^0_0\cY$. 

 Let $_0\Bun_{P_n}$ be the preimage of $_0\Bun_{G_n}$ in $\Bun_{P_n}$, and similarly for $_0^0\Bun_{\PP_n}$. 
 
 The stack $^0_0\Bun_{\PP_n}$ classifies $\cL\in{^0\Bun_n}$, an exact sequence (\ref{seq_cO_by_Omega}) given by $\xi_1\in\H^1(X,\Omega)$ and a point (\ref{seq_Omega_by_Sym2barcL}) of $_0\Bun_{P_n}$ given by $\gamma_1\in \H^1(X, \Omega^{-1}\otimes\Sym^2\cL)$. Now $^0_0\cY$ can be seen as a stack classifying
a point $(\cL, \gamma_1,\xi_1)\in{^0_0\Bun_{\PP_n}}$ and a section $s_1: \cL\to\Omega$ (here $s_1$ gives a new splitting of the exact sequence $0\to\Omega\to \Omega\oplus\cL\to\cL\to 0$).

 Another description of $^0_0\cY$ was given in Section~3.4.1, namely, it is a stack classifying $\cL\in{^0\Bun_n}$ and the exact sequences (\ref{seq_cO_by_Omega}) given by $\xi\in\H^1(X,\Omega)$, (\ref{seq_Omega_by_Sym2cL}) given by $\gamma\in \H^1(X, \Omega^{-1}\otimes\Sym^2\cL)$
and (\ref{seq_cO_by_cL}) given by $\delta\in\H^1(X,\cL)$.

 Given a point in the first description of $^0_0\cY$, the corresponding point $(\cL, \xi,\gamma,\delta)\in {^0_0\cY}$ in the second description is as follows. We have to take the push-forward of 
$$
(\gamma_1,\xi_1)\in \Ext^1(\Omega, \Sym^2(\cL\oplus\Omega))
$$ 
by $\epsilon\otimes\epsilon: \Sym^2(\cL\oplus\Omega)\to \Sym^2(\cL\oplus\Omega)$. Here $\epsilon$ is the automorphism of $\cL\oplus\Omega$ acting trivially on $\Omega$ and whose restriction to $\cL$ is $(\id, s_1): \cL\to \cL\oplus\Omega$. Thus, $\gamma=\gamma_1$, $\delta$ is the push-forward of $\gamma_1$ by $s_1\otimes \id+\id\otimes s_1: \Sym^2 \cL\to \cL\otimes\Omega$, and $\xi=\xi_1+\<s_1\otimes s_1, \gamma_1\>$. 

To simplify notations, we give the rest of the proof at the level of functions, the geometrization is straightforward. 
By (\cite{L1}, Proposition~7), the LHS of (\ref{iso_main_over_0_cY}) at $(\cL, \gamma_1,\xi_1, s_1)\in{^0_0\cY}$ equals
$$
\psi(\xi_1)\sum_{u\in\Hom(\cL,\Omega)} \psi(\<u\otimes u,\gamma_1\>)
$$
and the RHS of (\ref{iso_main_over_0_cY}) equals
$$
\sum_{u\in\Hom(\cL,\Omega)} \psi(\xi+\<\gamma, u\otimes u\>+\<\delta, u\>)=\sum_{u\in\Hom(\cL,\Omega)} \psi(\xi_1+\<(s_1+u)\otimes(s_1+u), \gamma_1\>)
$$
We are done. 
\end{Prf}
  
\medskip
  
\begin{Rem} The isomorphism (\ref{iso_explicit_Aute}) is not canonical, it depends on a choice of an isomorphism in (\cite{L1}, Proposition~7).
\end{Rem}
  
\bigskip

\centerline{\scshape 4. $P$-model and theta-lifting}  
  
\bigskip\noindent
4.1 Keep notations of Section~2.3. Let $^{sm}\Bun_n\subset\Bun_n$ be the open substack classifying $U\in\Bun_n$ such that $\H^0(X, \Omega\otimes\wedge^2 U)=0$. Write $^{sm}\Bun_P$ for the preimage of $^{sm}\Bun_n$ in $\Bun_P$. The restriction $^{sm}\Bun_P\to\Bun_H$ of $\nu_P$ is smooth, hence $^e\Bun_P\to\Bun_H$ is also smooth.

 To see that $\ocS_P$ is smooth, first consider the stack classifying $M^*\in \Bun_2$ equipped with an isomorphism $\det M^*\,\iso\,\Omega^{-1}$, a coherent sheaf $F$ of generic rank $n-2$ on $X$ and an exact sequence $0\to M^*\to L\to F\to 0$ on $X$. This stack is smooth, and its open substack given by the condition that $L$ is locally free identifies with $\ocS_P$.
 
\medskip\noindent 
\begin{Prf}\select{of Proposition~\ref{Pp_small_map}} 
The connected components of $\ocZ_P$ are $\ocZ{}_P^d$ for $d\in\ZZ$, and $\ocZ{}_P^d$ is irreducible. 

The stack $\cZ_{P,m}$ is smooth for any $m\ge 0$. Consider its connected component $\cU$ containing a point $\eta=(s: U\to M', D\in X^{(m)})$, where $U\in \Bun_n$, $M'\in\Bun_2$ is equipped with $\det M'\,\iso\,\Omega(-D)$, and $s$ is surjective. One checks that the dimension of this  connected component is
$$
m(1-n)-2\deg U+(n^2+3)(g-1)
$$
So, the  codimension of $\cU$ in the corresponding connected component of $\ocZ_P$ equals $m(n-1)$. The fibre of (\ref{map_pi_nice}) over $\eta$ is the scheme of upper modifications $M'\subset M$ such that $\div(M/M')=D$. This fibre is connected and its dimension is $m$. 
Our assertion follows.
\end{Prf}  
  
\medskip

 Our construction of $\cK_H$ is based on the following explicit formula for $\Aut_{G_1,H}$. Let $f_{\cS}: \cS_P\to \Bun_{G_1}\times\cY_P$ be the map sending $(U\toup{s}M)\in\cS_P$ to 
the collection $M\in \Bun_{G_1}$, $(U,v)\in\cY_P$ with $v: \wedge^2 U\to\Omega$, where $(U,v)$ is the image of $(U,M,s)$ under $\pi_P: \cS_P\to\cY_P$. As in Section~2.3.1, by some abuse of notation, write 
$$
\Four_{\cY_P,\psi}: \D^{\prec}(\Bun_{G_1}\times \cY_P)\to \D^{\prec}(\Bun_{G_1}\times\Bun_P)
$$ 
for the Fourier transform over $\Bun_{G_1}\times\Bun_n$.
The following is an immediate consequence of (\cite{L3}, Proposition~1).
  
\begin{Pp} 
\label{Pp_explicit_from_Wald_periods}
For the map $\id\times\nu_P: \Bun_{G_1}\times\Bun_P\to \Bun_{G_1}\times\Bun_H$ there is an isomorphism
\begin{equation}
\label{iso_explicit_for_Aut_GH}
(\id\times\nu_P)^*\Aut_{G_1,H}\otimes(\Qlb[1](\frac{1}{2}))^{\dimrel(\nu_P)}\,\iso\, \Four_{\cY_P,\psi}(f_{\cS !}(\Qlb[1](\frac{1}{2}))^b),
\end{equation}
where $b$ is a function of a connected component of $\cS_P$ whose value at $(U,M,s)\in\cS_P$ equals
${\dim\Bun_n+}\dim\Bun_{G_1}+\chi(U^*\otimes M)$. \QED
\end{Pp} 

 Note that for the function $b$ from Proposition~\ref{Pp_explicit_from_Wald_periods} its restriction to $\ocS_P$ equals $\dim\ocS_P$.
  
  Recall that for $a\in\ZZ$ one has the open substack $_a\Bun_n\subset\Bun_n$ introduced in Section~2.3.4. If $a'\le a$ then $_a\Bun_n\subset {_{a'}\Bun_n}$ is open. One checks that $\mathop{\cup}\limits_{a\in\ZZ} {_a\Bun_n}=\Bun_n$.   
Similarly, if $a'\le a$ then $_a\Bun_{G_1}\subset {_{a'}\Bun_{G_1}}$ is open and $\mathop{\cup}\limits_{a\in\ZZ} {_a\Bun_{G_1}}=\Bun_{G_1}$. For the complex $_a\tilde K$ given by (\ref{complex_a_tildeK_on_BunH}) this implies the following.
  
\begin{Cor} 
\label{Cor_passage_to_pH^0_first}
1) For all $a\in\ZZ$ there is an isomorphism over $_a\ocY_P$.
\begin{equation}
\label{iso_Cor1_over_a^eBunP}
\Four_{\cY_P,\psi}^{-1}
\nu_P^*(_a\tilde K)\otimes(\Qlb[1](\frac{1}{2}))^{\dimrel(\nu_P)}\,\iso\, \IC(\cZ_P)
\end{equation}
2) The isomorphism (\ref{iso_Cor1_over_a^eBunP}) still holds 
over $^e_a\ocY_P$ with $_a\tilde K$ replaced by $^p\cH^0(_a\tilde K)$.
\end{Cor}
\begin{Prf}  
1) Consider the restriction $_a\pi_P: \cS_P\times_{\Bun_{G_1}}{_a\Bun_{G_1}}\to \cZ_P$ of the map $\pi_P$ from Section~2.3.1. For the function $b$ as in Proposition~\ref{Pp_explicit_from_Wald_periods}, this proposition yields an isomorphism over the whole of $\Bun_P$
$$
\nu_P^*(_a\tilde K)\otimes(\Qlb[1](\frac{1}{2}))^{\dimrel(\nu_P)}\,\iso\, \Four_{\cY_P,\psi}((_a\pi_P)_!\Qlb)[b](\frac{b}{2})
$$ 

 Let us establish a canonical isomorphism
\begin{equation}
\label{iso_a_and_without_a} 
(_a\pi_P)_!\Qlb[b](\frac{b}{2})\,\iso\, \IC(\cZ_P)
\end{equation}
over the open substack $_a\ocY_P$ of $\cY_P$. Consider a $k$-point $(s: U\to M)$ of $\ocS_P$. Assume $U\in {_a\Bun_n}$ then for any line bundle $L$ on $X$ with $\deg L\le a$ and any morphism $y: M\to L$ the composition $U\toup{s}M\toup{y} L$ vanishes. Since $s$ is surjective at the generic point of $X$, $y$ also vanishes and $M\in {_a\Bun_{G_1}}$. Thus, $(_a\pi_P)_!\Qlb\,\iso\, (\pi_P)_!\Qlb$ over $_a\ocY_P$, and (\ref{iso_a_and_without_a}) follows from Proposition~\ref{Pp_small_map}. Part 1) follows.

\medskip\noindent   
2) Since $\nu_P: {^e_a\Bun_P}\to\Bun_H$ is smooth, the functor $\Four_{\cY_P,\psi}^{-1}\nu_P^*[\dimrel(\nu_P)]$ followed by restriction to $^e_a\ocY_P$ is exact for the perverse t-structures.
\end{Prf}

\medskip  
  
  The stack $\Bun_n$ is smooth, its connected components are indexed by $d\in\ZZ$. Namely, the connected component $\Bun_n^d$ of $\Bun_n$ classifies $U\in\Bun_n$ with $\deg U=d$. Write $\Bun_P^d$, $^e\Bun_P^d$ and so on for the preimage of $\Bun_n^d$ in the corresponding stack.
  
 Write $C(e,P)$ for the set of $d\in\ZZ$ such that the stack $^e\Bun_n^d$ from Section~2.3.1 is nonempty. For $a\in\ZZ$ write $_a^e\Bun_n^d$ for the preimage of $_a\Bun_n$ in $^e\Bun_n^d$. Given $d\in C(e,P)$, the stack $^e_a\Bun_n^d$ is nonempty for $a$ small enough. It is easy to see that for $d\in C(e,P)$ and $g=0$ (resp., $g\ge 1$) one has $d\le -n/2$ (resp., $d\le -(g-1)n/2$). 
 
  Write $Z(e,P)$ for the set of $d\in\ZZ$ such that $^e\ocZ{}^d_P$ is not empty. Clearly, $Z(e,P)\subset C(e,P)$. If $d\in Z(e,P)$ then $^e_a\ocZ{}^d_P$ is not empty for $a$ small enough. There is $N=N(g)$ such that if $d\le N$ then $d\in Z(e,P)$. 
 
\begin{Def} 
\label{Def_cK^d_H}
Let $a,d\in\ZZ$ be such that $^e_a\ocZ{}^d_P$ is nonempty. Then by Corollary~\ref{Cor_passage_to_pH^0_first} and Lemma~\ref{Lm_general_ab_categories} below, there is a unique irreducible subquotient $_a\cK^d_H$ of the perverse sheaf $^p\cH^0(_a\tilde K)$ equipped with an isomorphism
\begin{equation}
\label{iso_def_aK^d_H}
\Four_{\cY_P,\psi}^{-1}\nu_P^*(_a\cK^d_H)\otimes(\Qlb[1](\frac{1}{2}))^{\dimrel(\nu_P)}\,\iso\, \IC(\cZ_P)
\end{equation}
over $^e_a{}\ocY{}_P^d$. The perverse sheaf $_a\cK^d_H$ is defined up to a unique isomorphism. We can not conclude for the moment that (\ref{iso_def_aK^d_H}) holds over $^e_a\cY_P^d$, as the LHS could apriori be a non irreducible perverse sheaf.
  
  If $a'\le a$ and $^e_a\ocZ{}^d_P$ is nonempty then $^e_{a'}\ocZ{}^d_P$ is also nonempty. The open immersion $_a\Bun_{G_1}\hook{} {_{a'}\Bun_{G_1}}$ yields a morphism $_a\tilde K\to {_{a'}\tilde K}$, hence also a morphism of perverse sheaves 
\begin{equation}
\label{morphism_alpha_from_Def3}
\alpha: {^p\cH^0(_a\tilde K)}\to {^p\cH^0(_{a'}\tilde K)}
\end{equation}
After applying the functor 
$$
\Four_{\cY_P,\psi}^{-1}
\nu_P^*[\dimrel(\nu_P)]
$$ 
followed by restriction to $^e_a\ocY{}^d_P$, the map $\alpha$ becomes an isomorphism. By Lemma~\ref{Lm_general_ab_categories} below, $\alpha$ induces a natural isomorphism $_a\cK^d_H\,\iso\, {_{a'}\cK^d_H}$. For $d\in Z(e,P)$ define a perverse sheaf 
$$
\cK^d_H\in \P(\Bun_H)
$$ 
as $_a\cK^d_H$ for any $a$ small enough (such that $^e_a\ocZ{}^d_P$ is nonempty). We see that $\cK^d_H$ is defined up to a unique isomorphism. The perverse sheaf $\cK^d_H$ is equipped with an isomorphism over $^e\ocY{}^d_P$
$$
\Four_{\cY_P,\psi}^{-1}\nu_P^*(\cK^d_H)\otimes(\Qlb[1](\frac{1}{2}))^{\dimrel(\nu_P)}\,\iso\, \IC(\cZ_P) 
$$
\end{Def} 

\begin{Lm} 
\label{Lm_general_ab_categories}
Let $f: \cA\to\cB$ be an exact functor between  abelian categories. Let $F,F'$ be two objects of $\cA$ which are of finite length and $\alpha: F\to F'$ a morphism in $\cA$. Assume that $R=f(F)$ is an irreducible object of $\cB$, and $f(\alpha): f(F)\to f(F')$ is an isomorphism. Then $F$ admits a biggest subobject $F_0$ such that $f(F_0)=0$, let $F'_0\subset F'$ be the corresponding biggest subobject of $F'$. 
Then $F/F_0$ admits a unique irreducible subobject $F_1$, and $f(F_1)\,\iso\, R$. Define $F'_1\subset F'/F'_0$ similarly. Then $\alpha: F_0\to F'_0$ and the induced map $\alpha: F_1\to F'_1$ is an isomorphism. We refer to $F_1$ as \rm{the subquotient canonically associated to $(f, F)$}. 
\end{Lm}
\begin{Prf}
Let $G\subset F$ be a subobject such that $f(G)=0$ and maximal with this property. Let $G_1\subset F$ be another subobject such that $f(G_1)=0$. Write $\bar G_1$ for the image of $G_1$ in $F/G$, let $\bar G$ be the preimage of $\bar G_1$ under the projection $F\to F/G$. Then $f(\bar G)=0$, so $\bar G=G$. Thus, $G$ is the biggest subobject of $F$ such that $f(G)=0$. 

 If $F_1\subset F/F_0$ is an irreducible subobject then $f(F_1)\,\iso\, R$. Since $f$ is exact this $F_1$ is unique. 
Since $\alpha: F_0\to \alpha(F_0)$ is surjective, $f(\alpha): f(F_0)\to f(\alpha(F_0))$ is also surjective, hence $\alpha(F_0)\subset F'_0$. Our assertion follows.
\end{Prf}
 
\medskip

We will see in Section~7 that for all $d\in Z(e,P)$ of the same parity the perverse sheaves $\cK^d_H$ are canonically isomorphic to each other (cf. Proposition~\ref{Pp_K^d_coincide}). 

\bigskip

\centerline{\scshape 5. Comparison of $P$ and $Q$-models}

\bigskip\noindent
5.1 Keep notations of Section~2.3. The stack $\Bun_{P\cap Q}$ classifies a point (\ref{seq_cO_by_wedge2U}) of $\Bun_P$ together with an exact sequence on $X$
\begin{equation}
\label{seq_U'_by_W}
0\to W\to U\to U'\to 0
\end{equation}
with $W\in\Bun_1$, $U'\in\Bun_{n-1}$. Write $\nu_{P,Q}: \Bun_{P\cap Q}\to \Bun_Q$ and $\nu_{Q,P}: \Bun_{P\cap Q}\to\Bun_P$ for the natural maps. Write $^{\diamond \!}(\Bun_1\times\Bun_{n-1})\subset \Bun_1\times\Bun_{n-1}$ for the open substack given by 
\begin{equation}
\label{cond_diamond_P_and_Q}
\begin{array}{c}
\H^0(X, U'\otimes W)=\H^0(X,\wedge^2 U')=\Hom(U',W)=0\\ \\
\Hom(U', W\otimes\Omega)=\H^0(X, \Omega\otimes\wedge^2 U')=0
\end{array}
\end{equation} 
for $W\in\Bun_1$, $U'\in\Bun_{n-1}$. Write $^{\diamond \!}\Bun_{P\cap Q}$ for the preimage of $^{\diamond \!}(\Bun_1\times\Bun_{n-1})$ in $\Bun_{P\cap Q}$. 

  Our purpose is to prove the following.
 
\begin{Pp} 
\label{Pp_comparison_PQ}
There exists an isomorphism on $^{\diamond \!}\Bun_{P\cap Q}$
\begin{equation}
\label{iso_comparison_PQ_main}
\nu_{P,Q}^*K_{Q,\psi}\otimes(\Qlb[1](\frac{1}{2}))^
{\dimrel(\nu_{P,Q})}\,\iso\, \nu_{Q,P}^*K_{P,\psi}\otimes(\Qlb[1](\frac{1}{2}))^{\dimrel(\nu_{Q,P})}
\end{equation}
\end{Pp}  
\begin{Rem} Recall that $K_{P,\psi}$ is perverse over the open substack of $\Bun_P$ given by $\H^0(X,\wedge^2 U)=0$ for a point (\ref{seq_cO_by_wedge2U}) of $\Bun_P$, and $K_{Q,\psi}$ is perverse over the open substack of $\Bun_Q$ given by $\H^0(X, W\otimes V')=0$ for a point (\ref{seq_V'_by_W}) of $\Bun_Q$. The restrictions of $\nu_{P,Q}$ and of $\nu_{Q,P}$ to $^{\diamond \!}\Bun_{P\cap Q}$ are smooth, so both sides of (\ref{iso_comparison_PQ_main}) are perverse.  We will see in the course of the proof of Proposition~\ref{Pp_first_reduction_for_PQ} below that both sides of (\ref{iso_comparison_PQ_main}) are irreducible over each connected component of $^{\diamond \!}\Bun_{P\cap Q}$.
\end{Rem}
 
\medskip\noindent
5.2  The stack $\Bun_{P\cap Q}$ can also be seen as the stack classifying exact sequences on $X$
\begin{equation}
\label{seq_cO_by_wedge2U'}
0\to \wedge^2 U'\to ?\to \cO\to 0
\end{equation}
and (\ref{seq_V'_by_W}), where we denoted by $V'$ the corresponding point of $\Bun_{H_{n-1}}$.

 Write $\cS$ for the stack classifying a point $(W,U')\in {^{\diamond \!}(\Bun_1\times\Bun_{n-1})}$ together with the exact sequences (\ref{seq_U'_by_W}) and (\ref{seq_cO_by_wedge2U'}) on $X$. 
 
 Let $\cT$ be the stack over $\cS$ with fibre $\Hom(W\otimes U',\Omega)$. The conditions (\ref{cond_diamond_P_and_Q}) imply that $\Ext^1(W\otimes U', \Omega)=0$, so $\cT$ is a vector bundle. The natural projection $^{\diamond \!}\Bun_{P\cap Q}\to \cS$ is a torsor under the vector bundle $\cT^*$. Denote by
\begin{equation}
\label{seq_cO_by_cT^*}
0\to \cT^*\to \cE^*\to \cO\to 0
\end{equation}
the corresponding exact sequence of $\cO_{\cS}$-modules, so $^{\diamond \!}\Bun_{P\cap Q}$ is the preimage of $1$ in $\cE^*$.
 
  Let $\cT_{Q}$ be stack over $\cS$ with fibre $\Hom(W,  V'\otimes \Omega)$. The conditions (\ref{cond_diamond_P_and_Q}) imply that $\cT_Q$ is a vector bundle over $\cS$, and for a point of $\cS$ the sequence is exact
\begin{equation}
\label{seq_first_on_cS}
0\to \Hom(W, U'\otimes \Omega)\to \Hom(W, V'\otimes \Omega)\to \Hom(W\otimes U',  \Omega)\to 0
\end{equation}
Indeed, (\ref{cond_diamond_P_and_Q}) imply that $\Ext^1(W, U'\otimes\Omega)=0$.
  
  Define a full subcategory $\P^W(\cT_Q)\subset\P(\cT_Q)$ singled out by the following equivariance condition. Let $\cV\cT_Q$ be the vector bundle over $\cS$ classifying a point of $\cS$ and $t_1\in \Hom(W, U'\otimes\Omega)$. So, the sequence of vector bundles on $\cS$ is exact 
\begin{equation}
\label{seq_cT_by_cVcT_Q}
0\to \cV\cT_Q\to \cT_Q\to\cT\to 0
\end{equation}
The vector bundle $\cV\cT_Q$ acts on $\cT_Q$ by translations over $\cS$. Write $\ev_{\cV\cT_Q}: \cV\cT_Q\to\A^1$ for the map sending $(W,U', (\ref{seq_cO_by_wedge2U'}), (\ref{seq_U'_by_W}), t_1)$ to the pairing of $t_1$ with (\ref{seq_U'_by_W}). Define $\P^W(\cT_Q)\subset\P(\cT_Q)$ as the full subcategory of $(\cV\cT_Q, \ev_{\cV\cT_Q}^*\cL_{\psi^{-1}}$)-equivariant perverse sheaves on $\cT_Q$. For $F\in\P(\cT_Q)$ this means that
for the action and the projection maps $\act,\pr: \cV\cT_Q\times_{\cS}\cT_Q\to \cT_Q$ there is an isomorphism
$$
\act^*F\,\iso\, \pr^*F\otimes \ev_{\cV\cT_Q}^*\cL_{\psi^{-1}}
$$
whose restriction to the unit section is the identity, and it satisfies the corresponding associativity condition. If such an isomorphism exists then it is unique. Write $\D^W(\cT_Q)\subset\D^{\prec}(\cT_Q)$ for the full subcategory of complexes whose all perverse cohomology sheaves lie in $\P^W(\cT_Q)$.

  Let $\cT_P$ be the stack over $\cS$ with fibre $\Hom(\wedge^2U, \Omega)$. The conditions (\ref{cond_diamond_P_and_Q}) imply that $\cT_P$ is a vector bundle over $\cS$. For a point of $\cS$ the exact sequence $0\to  W\otimes U'\to \Lambda^{2}U\to \Lambda^{2}U'\to 0$ yields a sequence
\begin{equation}
\label{seq_second_on_cS}
0\to \Hom(\wedge^2 U',  \Omega)\to \Hom(\wedge^2 U,  \Omega)\to\Hom(W\otimes U',  \Omega)\to 0,
\end{equation}
which is exact because of (\ref{cond_diamond_P_and_Q}). Let $\cV\cT_P$ be the vector bundle over $\cS$ classifyig a point of $\cS$ and $v_1\in\Hom(\wedge^2 U',  \Omega)$, so 
\begin{equation}
\label{seq_cT_by_cVcT_P}
0\to \cV\cT_P\to \cT_P\to \cT\to 0
\end{equation}
is an exact sequence of vector bundles on $\cS$.

 Write $\ev_{\cV\cT_P}: \cV\cT_P\to\A^1$ for the map sending $(W, U',  (\ref{seq_cO_by_wedge2U'}), (\ref{seq_U'_by_W}), v_1)$ to the pairing of $v_1$ with (\ref{seq_cO_by_wedge2U'}). As above, one defines the category $\P^W(\cT_P)$ of $(\cV\cT_P, \ev^*_{\cV\cT_P}\cL_{\psi^{-1}})$-equivariant perverse sheaves on $\cT_P$, similarly for $\D^W(\cT_P)$.
  
\begin{Lm} 
\label{Lm-cup}
1) The push-forward of the exact sequence (\ref{seq_cT_by_cVcT_P}) by the morphism 
$\cV\cT_P\to\cO_{\cS}$ given by pairing with the extension (\ref{seq_cO_by_wedge2U'}) is canonically isomorphic  to the exact sequence $0\to \cO\to \cE\to \cT\to 0$ on $\cS$ dual to (\ref{seq_cO_by_cT^*}).\\
2) The push-forward of the exact sequence (\ref{seq_cT_by_cVcT_Q}) by the morphism 
$\cV\cT_Q\to \cO_{\cS}$  given by pairing with the extension (\ref{seq_U'_by_W}) is canonically isomorphic to the exact sequence $0\to \cO\to \cE\to \cT\to 0$ on $\cS$. 
\end{Lm}
\begin{Prf}
1) Dualizing (\ref{seq_first_on_cS}) one gets the exact sequence $0\to \H^1(X, W\otimes U')\to \H^1(X, {W\otimes V'})\to \H^1(X, W\otimes U'^*)\to 0$ on $\cS$. Part 1) follows from the fact that $\Bun_{P\cap Q}$ is the stack classifying $U',W$ as above and exact sequences $0\to \wedge^2 U'\to ?\to \cO\to 0$, $0\to W\to ?\to V'\to 0$ on $X$. 

\medskip\noindent
2) Dualizing (\ref{seq_second_on_cS}) one gets the exact sequence $0\to \H^1(X, W\otimes U')\to \H^1(X, \wedge^2 U)\to \H^1(X, \wedge^2 U')\to 0$ on $\cS$. Part 2) follows from the fact that $\Bun_{P\cap Q}$ is the stack classifying $U',W$ as above and exact sequences $0\to W\to U\to U'\to 0$, $0\to \wedge^2 U\to ?\to \cO\to 0$ on $X$. 
\end{Prf}
  
\medskip

 As above, define $\P^W(\cE)$ as the category of perverse sheaves on $\cE$ which are $(\cO, \cL_{\psi^{-1}})$-equivariant, similarly for the derived category $\D^W(\cE)$. Lemma~\ref{Lm-cup} yields canonical equivalences 
$$
\D^W(\cT_P)\xrightarrow[\widetilde{}]{\epsilon_P} \D^W(\cE)\xleftarrow[\widetilde{}]{\epsilon_Q}\D^W(\cT_Q)
$$ 
exact for the perverse t-structures.
 
 The Fourier transform $\Four_{\cE, \psi}: \D^{\prec}(\cE)\,\iso\,\D^{\prec}(\cE^*)$ yields an equivalence between the full subcategories on both sides 
$$
\Four_{\cE,\psi}: \D^W(\cE)\,\iso\, \D^{\prec}(^{\diamond}\Bun_{P\cap Q})
$$
\noindent
5.3.1 Let $\bar\cT_Q$ be the stack classifying $(W, U')\in {^{\diamond \!}(\Bun_1\times\Bun_{n-1})}$, an exact sequence (\ref{seq_cO_by_wedge2U'}) on $X$, and $t\in\Hom(W, V'\otimes\Omega)$. Here $V'\in\Bun_{H_{n-1}}$ is given by (\ref{seq_cO_by_wedge2U'}). 
 The projection $\bar\cT_Q\to \cY_Q$ is smooth, we set 
$$
\bar\cT\cZ_Q=\bar\cT_Q\times_{\cY_Q}\cZ_Q
$$ 
 
 Define the partial Fourier transform along $\Hom(W, U'\otimes \Omega)$ as the following equivalence
$$
\Four_{Q,\psi}:\D^{\prec}(\bar\cT_Q)\,\iso\, \D^W(\cT_Q)
$$ 
{\scshape Notation.} Write $(\alpha)$ for an exact sequence (\ref{seq_cO_by_wedge2U'}), $(\gamma)$ for an exact sequence (\ref{seq_U'_by_W}). 

Consider the diagram
$$
\begin{array}{ccc}
\cT_Q\getsup{p_Q} & \cV\cT_Q\times_{\cS}\cT_Q & \toup{a_Q} \bar\cT_Q\\
& \downarrow\lefteqn{\scriptstyle \ev_{\cV\cT_Q}}\\
&\AA^1,
\end{array}
$$
where $a_{Q}$ sends $t_1\in \Hom(W,U'\otimes\Omega)$, $(W,U', t,\alpha,\gamma)\in \cT_Q$ to $(W, U', \alpha, t+t_1)\in\bar\cT_Q$. The map $p_Q$ sends the same collection to $(W, U', t,\alpha,\gamma)\in \cT_Q$. The map $\ev_{\cV\cT_Q}$ sends the same collection to $\<t_1,\gamma\>$. Then
$$
\Four_{Q,\psi}(K)=(p_Q)_!(a_Q^*K\otimes\ev_{\cV\cT_Q}^*\cL_{\psi})\otimes(\Qlb[1](\frac{1}{2}))^{\dimrel(a_Q)}
$$
Let $\pr_Q: \cT_Q\to\bar\cT_Q$ be the projection forgetting $(\gamma)$. Note that $(\pr_Q)_!: \D^W(\cT_Q)\to \D(\bar\cT_Q)$ is quasi-inverse to $\Four_{Q,\psi}$. 

\medskip\noindent
5.3.2 Let $\bar\cT_P$ be the stack classifying $(W, U')\in {^{\diamond \!}(\Bun_1\times\Bun_{n-1})}$, an exact sequence (\ref{seq_U'_by_W}) on $X$, and $v\in\Hom(\wedge^2 U,\Omega)$. The projection $\bar\cT_P\to \cY_P$ is smooth. Set $\bar\cT\cZ_P=\bar\cT_P\times_{\cY_P}\cZ_P$.

 Define the partial Fourier transform along $\Hom(\wedge^2 U', \Omega)$ as the following equivalence
$$
\Four_{P,\psi}: \D^{\prec}(\bar\cT_P)\,\iso\, \D^W(\cT_P)
$$
Consider the diagram
$$
\begin{array}{ccc}
\cT_P\getsup{p_P} & \cV\cT_P\times_{\cS}\cT_P & \toup{a_P} \bar\cT_P\\
& \downarrow\lefteqn{\scriptstyle \ev_{\cV\cT_P}}\\
&\AA^1,
\end{array}
$$
where $a_P$ sends $v_1\in \Hom(\wedge^2 U',\Omega)$, $(W, U', v, \alpha,\gamma)\in \cT_P$ to $(W, U', v+v_1,\gamma)\in\bar\cT_P$. The map $p_P$ sends the same collection to $(W, U', v, \alpha,\gamma)\in \cT_P$. The map $\ev_{\cV\cT_P}$ sends the same collection to $\<v_1, \alpha\>$. Then
$$
\Four_{P,\psi}(K)=(p_P)_!(a_P^*K\otimes \ev_{\cV\cT_P}^*\cL_{\psi})\otimes(\Qlb[1](\frac{1}{2}))^{\dimrel(a_P)}
$$

Let us reduce Proposition~\ref{Pp_comparison_PQ} to the following result, whose proof is found in Section~5.5.
\begin{Pp} 
\label{Pp_first_reduction_for_PQ}
There is a canonical isomorphism in $\D^W(\cE)$
\begin{equation}
\label{iso_PQ_comparison_second_Pp}
\epsilon_Q\Four_{Q,\psi}(\IC(\bar\cT\cZ_Q))\,\iso\,\epsilon_P \Four_{P,\psi}(\IC(\bar\cT\cZ_P))
\end{equation}
\end{Pp}

\begin{Prf}\select{of Proposition~\ref{Pp_comparison_PQ}}

\medskip\noindent
It is formal to check that one has canonical isomorphisms
$$
\Four_{\cE,\psi}\epsilon_P\Four_{P,\psi}(\IC(\bar\cT\cZ_P))\,\iso\, \nu_{Q,P}^*K_{P,\psi}\otimes(\Qlb[1](\frac{1}{2}))^{\dimrel(\nu_{Q,P})}
$$
and
$$
\Four_{\cE,\psi}\epsilon_Q\Four_{Q,\psi}(\IC(\bar\cT\cZ_Q))\,\iso\, \nu_{P,Q}^*K_{Q,\psi}\otimes(\Qlb[1](\frac{1}{2}))^
{\dimrel(\nu_{P,Q})}
$$
Our assertion follows now from Propsition~\ref{Pp_first_reduction_for_PQ}.
\end{Prf}

\begin{Rem} 
\label{Rem_for_PQ_compatibility_density}
i) Let $^{sm}\cZ_Q\subset\cZ_Q$ be the open substack given by the condition that $W\hook{} V'\otimes\Omega$ is a subbundle. We set 
$$
^{sm}\bar\cT\cZ_Q=\bar\cT_Q\times_{\cY_Q}{^{sm}\cZ_Q}
$$ 
Since $^{sm}\cZ_Q$ is smooth,  $^{sm}\bar\cT\cZ_Q$ is also smooth. The conditions (\ref{cond_diamond_P_and_Q}) imply that $^{sm}\bar\cT\cZ_Q$ is dense in $\bar\cT\cZ_Q$. So, $\IC(\bar\cT\cZ_Q)$ is the intermediate extension from $^{sm}\bar\cT\cZ_Q$.

 Recall that $\cZ_{P,0}\subset\cZ_P$ denotes the open substack given by the condition that $v: \wedge^2 U\to\Omega$ is surjective. We set $^{sm}\bar\cT\cZ_P=\bar\cT_P\times_{\cY_P} {\cZ_{P,0}}$. The conditions (\ref{cond_diamond_P_and_Q}) imply that $^{sm}\bar\cT\cZ_P$ is dense in $\bar\cT\cZ_P$. So, $\IC(\bar\cT\cZ_P)$
is the intermediate extension from $^{sm}\bar\cT\cZ_P$.

 The connected components of $^{sm}\bar\cT\cZ_P$ are given by fixing the degrees of $W, U'$. The connected components of $^{sm}\bar\cT\cZ_Q$ are also given by fixing the degrees of $W,U'$. 
  
\noindent
ii) The open substack of $^{sm}\bar\cT\cZ_Q$ given by the condition that the composition $W\toup{t} V'\otimes\Omega\to U'^*\otimes\Omega$ is a subbundle is dense in $^{sm}\bar\cT\cZ_Q$.

 Similarly, the open substack of $^{sm}\bar\cT\cZ_P$ given by the condition that the composition $W\otimes U'\hook{} \wedge^2 U\toup{v}\Omega$ is surjective is dense in $^{sm}\bar\cT\cZ_P$.
\end{Rem}

\medskip\noindent
5.4 Let $^0\cT\subset\cT$ be the open substack classifying $(W, U', \alpha,\gamma)\in\cS$ and $s: W\to U'^*\otimes\Omega$ whose image is a rank one subbundle in $U'^*\otimes\Omega$. Let $^0\cT_Q$ (resp., $^0\cT_P$) be the preimage of $^0\cT$ under the projection $\cT_Q\to \cT$ (resp., under $\cT_P\to\cT$).

 Define a closed substack $^0\cX\subset {^0\cT}$ by the following conditions. A point $(W, U', \alpha,\gamma, s)$ of $^0\cT$ as above yields an exact sequence 
\begin{equation}
\label{seq_for_U'_n-2} 
0\to U'_{n-2}\to U'\toup{s} W^*\otimes\Omega\to 0
\end{equation} 
It induces the surjections  $U'^*\otimes W\to U'^*_{n-2}\otimes W$ and $\wedge^2 U'\to 
U'_{n-2}\otimes W^*\otimes\Omega$ of $\cO_X$-modules. Then $^0\cX$ is given by the conditions
\begin{itemize}
\item the image of $(\gamma)$ under $\H^1(X, U'^*\otimes W)\to \H^1(X, U'^*_{n-2}\otimes W)$ vanishes,
\item the image of $(\alpha)$ under $\H^1(X, \wedge^2 U')\to \H^1(X, U'_{n-2}\otimes W^*\otimes\Omega)$ vanishes.
\end{itemize}
Write $^0\cX_Q$ (resp., $^0\cX_P$) for the preimage of $^0\cX$ under $^0\cT_Q\to {^0\cT}$ (resp., under $^0\cT_P\to {^0\cT}$).

 Stratify $^0\cX$ by locally closed substacks $^0\cX_i$ indexed by $i\ge 0$ and given by the condition 
$$
\dim\Hom(U'_{n-2}, W)=i
$$ 
Write $^0\cX_{Q,i}$ (resp., $^0\cX_{P,i}$) for the preimage of $^0\cX_i$ in $^0\cX_Q$ (resp., in $^0\cX_P$).

\begin{Lm}
\label{Lm_supportQ}
The restriction of $\Four_{Q,\psi}(\IC_{\bar\cT\cZ_{Q}})$ to $^0\cT_Q$ is the extension by zero under $^0\cX_Q\hook{} {^0\cT_Q}$ of a perverse sheaf. 
This perverse sheaf is smooth along the stratification of $^0\cX_Q$ by $^0\cX_{Q,i}$, and its $*$-restriction to the stratum $^0\cX_{Q,i}$ is a shifted rank one local system. 
\end{Lm}
\begin{Prf}
A point of $\cS$ gives rise to the diagram, where the top line is an exact sequence
$$
\begin{array}{rcl}
0\to \Hom(W, U'\otimes \Omega)\to &\Hom(W, V'\otimes \Omega)& \xrightarrow{\beta}\Hom(W\otimes U',  \Omega)\to 0\\
&  \downarrow q   \\ 
&\Hom(W^2, \Omega^2)  
\end{array}
$$
For $s\in \Hom(W\otimes U',  \Omega)$ the restriction of $q$ to the affine subspace $\beta^{-1}(s)$ is affine, and the underlying linear map $\Hom(W, U'\otimes \Omega)\to \Hom(W, W^*\otimes\Omega^{2})$ is given by
the composition with $2s\in \Hom(U'\otimes  \Omega
, W^*\otimes  \Omega^2)$.

 Let $(W, U', \alpha, \gamma, s)\in {^0\cT}$ be such that the corresponding fibre of the composition 
$$
a_Q^{-1}(\bar\cT\cZ_Q)\toup{p_Q} \cT_Q\to \cT
$$ 
is non empty. Then there is $t: W\to V'\otimes\Omega$ extending $s: W\to U'^*\otimes\Omega$ such that the image of $t$ is isotropic. The map $s$ gives rise to the exact sequence (\ref{seq_for_U'_n-2}). Write $U'_n$ for the orthogonal complement of $U'_{n-2}$ in $V'$, so $U'_{n-2}\subset U'_n$ and $U'_n\in\Bun_n$. Moreover, $U'_n/U'_{n-2}\in\Bun_{H_1}$, so one has a canonical decomposition 
$$
U'_n/U'_{n-2}\,\iso\, (W^*\otimes\Omega) \oplus (W\otimes \Omega^{-1})
$$
as a sum of isotropic subbundles. We get the diagram
\begin{equation}
\label{diag_for_supp_Q}
\begin{array}{ccccc}
0\to U'_{n-2} \to & U'_n & \to & U'_n/U'_{n-2} & \to 0\\
&&Ê\nwarrow\lefteqn{\scriptstyle t} & \uparrow\\
&&& W\otimes\Omega^{-1},
\end{array}
\end{equation}
where the vertical arrow is the inclusion as an isotropic subbundle. This shows that the image of $(\alpha)$ under $\H^1(X, \wedge^2 U')\to \H^1(X, U'_{n-2}\otimes W^*\otimes\Omega)$ vanishes.  

  Now the fibre of $a_Q^{-1}(\bar\cT\cZ_Q)\toup{p_Q} \cT_Q$ over $(W, U', \alpha,\gamma, t)$ is the scheme of $t_1\in\Hom(W, U'\otimes\Omega)$ such that the image of $t+t_1: W\to V'\otimes\Omega$ is isotropic. Using the exact sequence 
$$
0\to U'_{n-2}\otimes W^*\otimes\Omega\to W^*\otimes U'\otimes\Omega\toup{s} W^{-2}\otimes\Omega^2\to 0,
$$
one identifies this scheme with $\H^0(X, U'_{n-2}\otimes W^*\otimes\Omega)$. For any such $t_1$, the image of $t+t_1$ is an isotropic subbundle in $V'\otimes\Omega$. So, one has to integrate over $\Hom(W, U'_{n-2}\otimes\Omega)$
the restriction of $\cL_{\psi}$ under the composition
$$
\Hom(W, U'_{n-2}\otimes\Omega)\hook{}\Hom(W, U'\otimes \Omega)\toup{\gamma} \A^1
$$
This local system is trivial iff the image of $\gamma$ under $\H^1(X, W\otimes U'^*)\to \H^1(X, U'^*_{n-2}\otimes W)$ vanishes.

 Note that $\Hom(W, U'_{n-2}\otimes\Omega)^*\,\iso\, \H^1(X, U'^*_{n-2}\otimes W)$, and 
$$
\chi(U'^*_{n-2}\otimes W)=\chi(W\otimes U'^*)-\chi(W^2\otimes\Omega^{-1})
$$ 
is fixed on each connected component of $^0\cT$. So, the stratification of $^0\cX$ by $^0\cX_i$, $i\ge 0$ coincides with the one given by fixing $\dim\H^0(X, U'_{n-2}\otimes W^*\otimes\Omega)$.
\end{Prf}
  
\medskip

 The stack $^0\cT$ is smooth, its connected components are given by fixing the degrees of $W, U'$.
 
\begin{Lm}
\label{Lm_non_empty}
Consider a connected component $\cC$ of $^0\cT$ given by $\deg U'=a_U$, $\deg W=a_W$. Assume $a_U<0$ and $a_W$ sufficiently small compared to $a_U$ (it suffices to require $(n-3)a_W\le a_U+(n-4)(g-1)$ and $a_W\le g-2$). Then the open substack of $\cC$ given by $\Hom(U'_{n-2}, W)=0$ is non empty.   
\end{Lm}
\begin{Prf}
Write $\cB$ for the connected component of $\Bun_1\times\Bun_{n-1}$ given by $\deg W=a_W, \deg U'=a_U$. 
Write $\cP$ for the stack classifying $U'_{n-2}\in\Bun_{n-2}$, $W\in\Bun_1$ with $\deg U'_{n-2}=a_U+a_W-(2g-2)$, $\deg W=a_W$, and an exact sequence (\ref{seq_for_U'_n-2}) on $X$. The stack $\cP$ is smooth an irreducible.

 For a point of $\cP$ one has $\chi(W\otimes U'^*_{n-2})\le 0$. So, the open substack $^0\cP\subset \cP$ given by $\Hom(U'_{n-2}, W)=0$ is non empty. 
 
 Let $_0\cB\subset\cB$ be the open substack given by $\H^0(X, W\otimes U'(x))=0$ for any $x\in X$. Under our assumptions, for $(W,U')\in\cB$ one has $\chi(W\otimes U'(x))\le 0$, so $_0\cB$ is nonempty. The stack $_0\cB$ is contained in the image of the map $\xi: \cP\to \cB$ sending the above point to $(W, U')$. So, $\xi$ is dominant. Write $^{\diamond\!}\cB\subset\cB$ for the preimage of $^{\diamond\!}(\Bun_1\times\Bun_{n-1})$ in $\cB$. Assume $^{\diamond\!}\cB$ non empty. Let $^{\diamond\!}\cP=\xi^{-1}(^{\diamond\!}\cB)$. Since $\cP$ is irreducible, $^0\cP\cap {^{\diamond\!}\cP}$ is non empty. Our assertion follows.
\end{Prf}

\begin{Lm}  
\label{Lm_supportP}
The restriction of $\Four_{P,\psi}(\IC_{\bar\cT\cZ_{P}})$ to $^0\cT_P$ is the extension by zero under $^0\cX_P\to {^0\cT_P}$ of a perverse sheaf. This perverse sheaf is smooth along the stratification of $^0\cX_P$ by $^0\cX_{P,i}$, and its $*$-restriction to the stratum $^0\cX_{P,i}$ is a shifted rank one local system. 
\end{Lm}
\begin{Prf}
Consider a point $(W,U',\alpha, \gamma, s)\in{^0\cT}$ such that the fibre over this point of the composition $a_P^{-1}(\bar\cT\cZ_P)\toup{p_P}\cT_P\to\cT$ is non empty. Then there is $v:\wedge^2 U\to\Omega$ extending $s: U'\otimes W\to\Omega$ such that $(U, v)\in\cZ_P$. The map $s$ gives rise to the exact sequence (\ref{seq_for_U'_n-2}).

 Note that $(U,v)\in \cZ_{P,0}$, that is, $v: \wedge^2 U\to\Omega$ is surjective, because its restriction to $U'\otimes W$ is already surjective. This point of $\cZ_{P,0}$ gives rise to $M\in\Bun_{G_1}$ with an exact sequence $0\to U_{n-2}\to U\toup{\tilde s} M\to 0$ such that $v=\wedge^2 \tilde s$. By our assumption, the composition $W\to U\to M\to M^*\otimes\Omega\to U^*\otimes\Omega$ is a rank one subbundle, so $W\subset M$ is also a subbundle, and we get a diagram
\begin{equation}
\label{diag_for_P-model}
\begin{array}{ccccccc}
0\to  & W &\to &M & \to &W^*\otimes\Omega& \to 0\\  
& \uparrow\lefteqn{\scriptstyle\id} && \uparrow\lefteqn{\scriptstyle\tilde s} && \uparrow\lefteqn{\scriptstyle s}\\
0\to & W & \to & U & \to & U' & \to 0\\
&&& \uparrow && \uparrow\\
&&& U_{n-2} && U'_{n-2}
\end{array}
\end{equation}
This diagram induces an isomorphism $U_{n-2}\,\iso\, U'_{n-2}$, so the image of $(\gamma)$ in $\H^1(X, U'^*_{n-2}\otimes W)$ vanishes. So, the fibre of $a_P^{-1}(\bar\cT\cZ_P)\toup{p_P}\cT_P$ over $(W,U',\alpha,\gamma, v)\in {^0\cT_P}$ identifies with the scheme of sections $U'_{n-2}\to U$ making the following diagram commutative
$$
\begin{array}{ccccc}
0\to W \to  & U & \to &U' &\to 0\\
& & \nwarrow & \uparrow\\
&&& U'_{n-2}
\end{array}
$$
The group $\Hom(U'_{n-2}, W)$ acts freely and transitively on this fibre. The local system $\ev_{\cV\cT_P}^*\cL_{\psi}$ changes under this action by the character $\Hom(U'_{n-2}, W)\subset \Hom(\wedge^2 U',\Omega)\toup{\alpha} \A^1$.
This character is trivial iff the image of $\alpha$ under $\H^1(X, \wedge^2 U')\to \H^1(X, W^*\otimes U'_{n-2}\otimes\Omega)$ vanishes. Clearly, over the locus of $^0\cT_{P,i}$ one gets a shifted rank one local system. 
\end{Prf}
  
\medskip


  
\medskip\noindent
5.5 \select{Proof of Proposition~\ref{Pp_first_reduction_for_PQ}}

\medskip\noindent
5.5.1 By Remark~\ref{Rem_for_PQ_compatibility_density} ii), the perverse sheaf $\Four_{Q,\psi}(\IC(\bar\cT\cZ_Q))$ is the intermediate extension under $^0\cT_Q\hook{}\cT_Q$, and $\Four_{P,\psi}(\IC(\bar\cT\cZ_P))$ is the intermediate extension under $^0\cT_P\hook{}\cT_P$. So, it suffices to establish (\ref{iso_PQ_comparison_second_Pp}) over $^0\cE=\cE\times_{\cT}{^0\cT}$. 
  
  First, let us define a full subcategory 
$
\P^W(\cT\times_{\cS} {^{\diamond\!}\Bun_{P\cap Q}})\subset \P(\cT\times_{\cS} {^{\diamond\!}\Bun_{P\cap Q}})
$. 
Write $\ev_{\cT}: \cT\times_{\cS}\cT^*\times_{\cS}\to\A^1$ for the natural pairing between $\cT$ and $\cT^*$. Recall that $^{\diamond\!}\Bun_{P\cap Q}\to \cS$ is a torsor under $\cT^*$. As in Section~5.2, one defines the category $\P^W(\cT\times_{\cS} {^{\diamond\!}\Bun_{P\cap Q}})$ of $(\cT^*, \ev_{\cT})$-equivariant perverse sheaves on $\cT\times_{\cS} {^{\diamond\!}\Bun_{P\cap Q}}$. Similarly for the derived category $\D^W(\cT\times_{\cS} {^{\diamond\!}\Bun_{P\cap Q}})$.
  
  One has a canonical equivalence 
\begin{equation}
\label{eq_cE_one_more}  
\varepsilon: \D^W(\cE)\,\iso\, \D^W(\cT\times_{\cS}{^{\diamond\!}\Bun_{P\cap Q}})
\end{equation}
exact for the perverse t-structures. It is characterised by the following. Write $\ev_{\cE}$ for the composition 
$$
\cE\times_{\cS}{^{\diamond\!}\Bun_{P\cap Q}}\hook{} \cE\times_{\cS}\cE^*\to\A^1,
$$
where the second map is the natural pairing. Then (\ref{eq_cE_one_more}) sends $K$ to the complex $K'$ equipped with an isomorphism $\pr_1^*K\otimes \ev_{\cE}^*\cL_{\psi}\,\iso\, (q_{\cE}\times\id)^*K'[1](\frac{1}{2})$, where $q_{\cE}: \cE\to\cT$ is the natural surjection, and $\pr_1: \cE\times_{\cS}{^{\diamond\!}\Bun_{P\cap Q}}\to \cE$ is the projection.
Such $K'$ is defined up to a unique isomorphism.
  
\medskip\noindent
5.5.2 Let us define a morphism 
$$
\ev_{Q,i}: {^0\cX_i}\times_{\cS}{^{\diamond\!}\Bun_{P\cap Q}}\to \A^1
$$ 
Consider a point of $^0\cX_i\times_{\cS}{^{\diamond\!}\Bun_{P\cap Q}}$ given by $(U',W, \alpha,\gamma,s)\in {^0\cX_i}$ and an exact sequence (\ref{seq_V'_by_W}) giving rise to $V\in {^{\diamond\!}\Bun_{P\cap Q}}$. Since the image of $\alpha$ in $\H^1(X, U'_{n-2}\otimes W^*\otimes\Omega)$ vanishes, we amy pick a lifting of $s$ to $t: W\to V'\otimes\Omega$ such that the image of $t$ is isotropic. Such $t$ is defined uniquely up to adding an element $t_1\in\Hom(W, U'_{n-2}\otimes\Omega)$. The map $\ev_{Q,i}$ sends this point to the pairing of $t$ with (\ref{seq_V'_by_W}). This is well-defined, because the image of $\gamma$ in $\H^1(X, U'^*_{n-2}\otimes W)$ vanishes.
  
 Let us define a morphism
$$
\ev_{P,i}: {^0\cX_i}\times_{\cS}{^{\diamond\!}\Bun_{P\cap Q}}\to \A^1
$$ 
Consider a point of $^0\cX_i\times_{\cS}{^{\diamond\!}\Bun_{P\cap Q}}$ given by $(U',W, \alpha,\gamma,s)\in {^0\cX_i}$ and an exact sequence (\ref{seq_cO_by_wedge2U}) giving rise to $V\in {^{\diamond\!}\Bun_{P\cap Q}}$. Since the image of $\gamma$ in $\H^1(X, U'^*_{n-2}\otimes W)$ vanishes, we may pick a lifting $v: \wedge^2 U\to\Omega$ of $s$ such that $(U,v)\in \cZ_{P,0}$. Such $v$ is uniquely defined up to adding an element $v_1\in \Hom(U'_{n-2}, W)\subset\Hom(\wedge^2 U', \Omega)\subset \Hom(\wedge^2 U,\Omega)$. Let $\ev_{P,i}$ send this point to the pairing of $v$ with (\ref{seq_cO_by_wedge2U}). The result is well-defined, because $\<v_1, \alpha\>=0$. Note that $\ev_{P,i}=\ev_{Q,i}$.
  
\medskip\noindent  
5.5.3  The $*$-restriction of $\varepsilon\epsilon_Q\Four_{Q,\psi}(\IC(\bar\cT\cZ_Q))$ to ${^0\cX_i}\times_{\cS}{^{\diamond\!}\Bun_{P\cap Q}}$ identifies (up to a shift and a twist) with $\ev_{Q,i}^*\cL_{\psi}$. The $*$-restriction of $\varepsilon\epsilon_P \Four_{P,\psi}(\IC(\bar\cT\cZ_P))$ to ${^0\cX_i}\times_{\cS}{^{\diamond\!}\Bun_{P\cap Q}}$ identifies (up to a shift and a twist) with $\ev_{P,i}^*\cL_{\psi}$. 
    
  After applying $\varepsilon$, it suffices to establish (\ref{iso_PQ_comparison_second_Pp}) over $^0\cT\times_{\cS}{^{\diamond\!}\Bun_{P\cap Q}}$. For each connected component of $^0\cX$ there is $i$ such that $^0\cX_i$ is dense in this component. This concludes the proof of Proposition~\ref{Pp_first_reduction_for_PQ}. 

\medskip\noindent
5.6 {\scshape Pointwise Euler characteristics}
Note that the maps $\nu_P: {^e\Bun_P}\to\Bun_H$ and $\nu_Q: {^u\Bun_Q}\to\Bun_H$ are surjective.
 
\begin{Pp} 
\label{Pp_Euler_char_pointwise}
There is a function $E_{\cK}: \Bun_H(k)\to\ZZ$ with the following properties. \\
1) For any $k$-point $\eta\in {^e\Bun_P}$ over $V\in \Bun_H(k)$ one has
$$
\chi(K_{P,\psi}\mid_{\eta})=(-1)^{\dimrel(\nu_P)}E_{\cK}(V)
$$ 
2) For any  $k$-point $\eta\in{^u\Bun_Q}$ over $V\in\Bun_H(k)$ one has
$$
\chi(K_{Q,\psi}\mid_{\eta})=(-1)^{\dimrel(\nu_Q)}E_{\cK}(V)
$$
\end{Pp}
\begin{Prf}
For $r\ge 1$ consider the stack $\cD_r$ classifying collections: $(W_i\subset U_i\subset V)\in {^{\diamond}\Bun_{P\cap Q}}$ for $1\le i\le r$, here $W_i\in\Bun_1$, $U_i\in\Bun_n$, $V\in\Bun_H$, and inclusions $W_i\subset U_{i+1}$ whose image is a subbundle such that $(W_i\subset U_{i+1}\subset V)\in {^{\diamond}\Bun_{P\cap Q}}$ for $1\le i<r$.

 Let $f_r: \cD_r\to \Bun_P\times_{\Bun_H}\Bun_P$ be the map sending the above point to $(U_1\subset V, U_r\subset V)$.  
The union of the images of $f_r$ for all $r\ge1$ contains
$^e\Bun_P\times_{\Bun_H}{^e\Bun_P}$. If $(U\subset V, U'\subset V)$ is in the image of some $f_r$ then, by Proposition~\ref{Pp_comparison_PQ}, the pointwise Euler characteristics of $K_{P,\psi}$ at $(U\subset V)$ and $(U'\subset V)$ coincide.
Since $\nu_P: {^e\Bun_P}\to\Bun_H$ is surjective, part 1) follows.

 Let $g_r: \cD_r\to \Bun_Q\times_{\Bun_H}\Bun_Q$ be the map sending a point of $\cD_r$ to $(W_1\subset V, W_r\subset V)$. Using $g_r$ one similarly proves part 2).
\end{Prf}

\bigskip

\centerline{\scshape 6. Comparison of $P$ and $R$-models}

\bigskip\noindent
6.1.1 Keep notations of Section~2.3. Recall that $\Bun_R$ classifies $V\in\Bun_H$ and an isotropic subbundle $U_2\subset V$ with $U_2\in\Bun_2$. Write $V_{-2}$ for the orthogonal complement of $U_2$ in $V$, so $V'=V_{-2}/U_2\in \Bun_{H_{n-2}}$. Let $^{sm}(\Bun_2\times\Bun_{H_{n-2}})\subset \Bun_2\times\Bun_{H_{n-2}}$ be the open substack given by 
$$
\H^0(X, \Omega\otimes \wedge^2 U_2)=\H^0(X, \Omega\otimes U_2\otimes V')=0
$$ 
for $(U_2, V')\in \Bun_2\times\Bun_{H_{n-2}}$. Let $^{sm}\Bun_R$ be the preimage of $^{sm}(\Bun_2\times\Bun_{H_{n-2}})$ in $\Bun_R$. Recall the map $\nu_R: \Bun_R\to\Bun_H$ from Section~2.3.3. The restriction $^{sm}\Bun_R\to\Bun_H$ of $\nu_R$ is smooth. So, $^w\Bun_R\to\Bun_H$ is also smooth.  
  
  The map $f_R: \cY_R\to \Bun_R$ is a vector bundle over the open substack $^w\Bun_R$.
  
\medskip\noindent  
6.1.2  Write $\bar R$ for the quotient of $R$ by the center of the unipotent radical of $R$. The stack $\Bun_{\bar R}$ classifies $V'\in\Bun_{H_{n-2}}$, $U_2\in\Bun_2$ and an exact sequence 
\begin{equation}
\label{seq_V'_by_U2}
0\to U_2\to V_{-2}\to V'\to 0
\end{equation}  
Write $\cY_{\bar R}$ for the stack classifying a point of $\Bun_{\bar R}$ as above and an exact sequence on $X$
\begin{equation}
\label{seq_cO_by_wedge^2U_2}
0\to \wedge^2 U_2\to ?\to \cO\to 0
\end{equation}
Then $\cY_{\bar R}$ is a group stack over $\Bun_{\bar R}$, it acts on $\Bun_R$ over $\Bun_{\bar R}$ as follows.
If an $R$-torsor $\cF$ on $X$ is given by a collection $(U_2\subset V)$ as above, the sheaf $\cA_{\cF}$ of automorphisms of $\cF$ acting trivially on $\cF\times_{R}\bar  R$ identifies canonically with $\wedge^2 U_2$. The action map $\cY_{\bar R}\times_{\Bun_{\bar R}}\Bun_R\to\Bun_R$ sends $(\cF, (\ref{seq_cO_by_wedge^2U_2}))$ to $\cF\times_{\cA_{\cF}} \cF'$, where $\cF'$ is the $\cA_{\cF}$-torsor given by (\ref{seq_cO_by_wedge^2U_2}). In more elementary terms, $\cF\times_{\cA_{\cF}} \cF'$ is given by the exact sequence $0\to V_{-2}\to \tilde V\to U_2^*\to 0$, which is the sum of $0\to V_{-2}\to V\to U_2^*\to 0$ with the push-forward via $U_2\subset V_{-2}$ of the sequence $0\to U_2\to ?\to U_2^*\to 0$ given by (\ref{seq_cO_by_wedge^2U_2}).

 Write $a_R: \cY_{\bar R}\times_{\Bun_{\bar R}}\cY_R\to \cY_R$ for the action map defined similarly. This action on a point $(U_2\subset V,\, v_2)\in\cY_R$ does not change $v_2: \wedge^2 U_2\to\Omega$. 
 
 As in Section~2.3.3, we write $^w\Bun_{\bar R}$, $^w\cY_{\bar R}$ and so on for the preimage of $^w(\Bun_2\times\Bun_{H_{n-2}})$ in the corresponding stack.
The projection $^w\cY_{\bar R}\to {^w\Bun_{\bar R}}$ is a vector bundle. One checks that $^w\Bun_R\to {^w\Bun_{\bar R}}$ is a torsor under this vector bundle (for the above action).
 
 Write $\ev_R: \cY_{\bar R}\times_{\Bun_{\bar R}}\cY_R\to \A^1$ for the map sending $(U_2\subset V, v_2, (\ref{seq_cO_by_wedge^2U_2}))$ to the natural pairing of $v_2$ with (\ref{seq_cO_by_wedge^2U_2}). 
 
 As in Section~5.2, one defines the category $\P^W(^w\cY_R)$ of $(\cY_{\bar R}, \ev_R^*\cL_{\psi})$-equivariant perverse sheaves on $^w\cY_R$. This is the category of perverse sheaves $F$ on $^w\cY_R$ equipped with an isomorphism
$$
a_R^*F\,\iso\,\pr_2^*F\otimes\ev_R^*\cL_{\psi}
$$
over $^w\cY_{\bar R}\times_{\Bun_{\bar R}}{^w\cY_R}$
whose restriction to the unit section is the identity, and satisfying the corresponding associativity condition. Write $\D^W(^w\cY_R)\subset \D^{\prec}(^w\cY_R)$ for the full subcategory of complexes whose all perverse cohomology sheaves lie in $\P^W(^w\cY_R)$. 

 The Fourier transform 
\begin{equation}
\label{Four_R_psi-1} 
 \Four_{R,\psi^{-1}}: \D^{\prec}(^w\Bun_R)\,\iso\,\D^W(^w\cY_R)
\end{equation}
is the following equivalence. Consider the diagram
$$
\begin{array}{ccc}
^w\cY_R \getsup{p_R} & \cY_{\bar R}\times_{\Bun_{\bar R}}{^w\cY_R} & \toup{a_R} {^w\cY_R}\toup{f_R} {^w\Bun_R}\\
& \downarrow\lefteqn{\scriptstyle \ev_R}\\
& \A^1,
\end{array}
$$
where $p_R$ is the projection. We set 
$$
\Four_{R,\psi^{-1}}(K)=(p_R)_!(a_R^*f_R^*K\otimes \ev_R^*\cL_{\psi})\otimes(\Qlb[1](\frac{1}{2}))^{\dimrel(f_R\comp a_R)}
$$
It is exact for the perverse t-structures and the functor $(f_R)_ !: \D^W(^w\cY_R)\,\iso\, \D^{\prec}(^w\Bun_R)$ is quasi-inverse to $\Four_{R,\psi^{-1}}$.

One similarly defines the category $\D^W(^w\ocY_R)$. Note that for any $K\in \P^W(^w\ocY_R)$ one has $(j_R)_{!*}(K)\in \P^W(^w\cY_R)$ for the open immersion $j_R: \ocY_R\hook{}\cY_R$ from Section~2.3.3.

\medskip\noindent
6.2.1 {\scshape the map $e_R$}

\medskip\noindent
Given a vector bunbdle $\cM$ on $X$, a line bundle $\cA$ on $X$ and a symplectic form $\wedge^2\cM\to\cA$, we write $H(\cM)=\cM\oplus\cA$ for the Heisenberg group scheme on $X$ with operation
$$
(m_1, a_1)(m_2,a_2)=(m_1+m_2, a_1+a_2+\frac{1}{2}\<m_1, m_2\>)
$$
(The line bundle $\cA$ is usually clear from the context, and we omit it in our notation).

 Given $(U_2, V')\in \Bun_2\times\Bun_{H_{n-2}}$, the vector bundle $U_2\otimes V'$ is equipped with a natural symplectic form $\wedge^2(U_2\otimes V')\to \wedge^2 U_2$, so one gets the corresponding Heisenberg group $H(U_2\otimes V')$.

 Now $\Bun_R$ identifies canonically with the stack classifying $(U_2, V')\in \Bun_2\times\Bun_{H_{n-2}}$ and a torsor on $X$ under the group scheme $H(U_2\otimes V')$.
 
 Write $\Mod_2$ for the stack classifying $U_2\in\Bun_2$ with an upper modification $s_2: U_2\hook{} M$, here $M\in\Bun_2$ and $s_2$ is an inclusion of coherent $\cO_X$-modules.
 
  Consider the stack $\Mod_2\times_{\Bun_2}\Bun_R$ classifying $(U_2\subset V)\in\Bun_R$ and $(s_2: U_2\hook{}M)\in\Mod_2$. Let us define a morphism 
$$
e_R: \Mod_2\times_{\Bun_2}\Bun_R\to\Bun_R
$$
For a point of the source write $V'=V_{-2}/U_2$, where $V_{-2}$ is the orthogonal complement of $U_2$ in $V$. The map $s_2$ yields an inclusion of coherent $\cO_X$-modules $H(U_2\otimes V')\subset H(M\otimes V')$, which is a homomorphism of group schemes over $X$. View $(U_2\subset V)\in \Bun_R$ as a triple $(U_2, V', \cF)$, where $\cF$ is a torsor on $X$ under $H(U_2\otimes V)$. Let $\tilde\cF$ be  the torsor under $H(M\otimes V')$ on $X$ obtained from $\cF$ by the extension of the structure group $H(U_2\otimes V')\to H(M\otimes V')$. Then $(M, V', \tilde\cF)\in \Bun_R$ is given by some pair $(M\subset \tilde V)\in\Bun_R$. By definition, $e_R$ sends $(U_2\subset V, U_2\subset M)$ to $(M\subset\tilde V)$. 

\begin{Rem} Let $(U_2\subset V, U_2\hook{s_2}M)\in \Mod_2\times_{\Bun_2}\Bun_R$ and $(M\subset \tilde V)$ be its image by $e_R$. Let $U\subset V$ be an isotropic subbundle of rank $n$ such that $U_2\subset U$. Define $U'$ by the exact sequence $0\to U_2\to U\to U'\to 0$. Let 
\begin{equation}
\label{seq_U'_by_M}
0\to M\to \tilde U\to U'\to 0
\end{equation}
be the push-forward of the latter exact sequence by $s_2: U_2\to M$. The point $(U\subset V)\in\Bun_P$ is given by an exact sequence (\ref{seq_cO_by_wedge2U}). Let 
\begin{equation}
\label{seq_cO_by_wedge2tildeU} 
0\to\wedge^2 \tilde U\to ?\to \cO\to 0
\end{equation}
be the push-forward of this exact sequence by $\wedge^2 U\hook{} \wedge^2 \tilde U$. Then (\ref{seq_cO_by_wedge2tildeU}) together with $M\subset \tilde U$ is a point of $\Bun_{P\cap R}$ whose image in $\Bun_R$ identifies canonically with $(M\subset \tilde V)$.
\end{Rem}
   
\medskip\noindent  
6.2.2 Recall the stack $\cX_R$ from Section~2.3.3. Let us define a morphism 
$$
\rho_R: \cX_R\to\Bun_{\GG_{2n-4}}
$$ 
To do so, we introduce the following.
\begin{Def} Given $(U_2\subset V)\in\Bun_R$, the vector bundle $U_2\otimes V$ is equipped with a symplectic form $\wedge^2 (U_2\otimes V)\to\wedge^2 U_2$. Consider then $M_1=(\Sym^2 U_2)^{\perp}/\Sym^2 U_2$, where $(\Sym^2 U_2)^{\perp}$ is the orthogonal complement of $\Sym^2 U_2$ in $U_2\otimes V$. So, $M_1$ is equipped with a symplectic form $\wedge^2 M_1\to \wedge^2 U_2$ and a line subbundle $\wedge^2 U_2\subset M_1$. We will refer to $M_1$ with these structures as \select{the symplectic-Heisenberg bundle} associated to $(U_2\subset V)\in\Bun_R$.
\end{Def} 

Consider a point $(U_2\subset V, U_2\hook{s_2} M)\in \cX_R$, here $M$ is an upper modification of $U_2\in\Bun_2$ equipped with $\det M\,\iso\,\Omega$. Let $(M\subset \tilde V)
\in\Bun_R$ be the image of this point under $e_R$. By definition, $\rho_R$ sends the above point of $\cX_R$ to the symplectic-Heisenberg bundle $(\det M\subset M_1)$ associated to $(M\subset \tilde V)$. Since we are given an isomorphism $\det M\,\iso\,\Omega$, this symplectic-Heisenberg bundle is a point of $\Bun_{\GG_{2n-4}}$. Moreover, by (\cite{L2}, Lemma~1), for the above point of $\cX_R$ one has a canonical $\ZZ/2\ZZ$-graded isomorphism
\begin{equation}
\label{iso_detRG_MotimesV'}
\det\RG(X, M\otimes V')\,\iso\, \det\RG(X, V')^2\otimes\det\RG(X, M)^{2n-4}\otimes\det\RG(X, \cO)^{8-4n}
\end{equation}
We lift $\rho_R$ to a morphism (\ref{map_cXP^H_to_BuntG}) sending the above point of $\cX_R$ to the collection $(\Omega\subset M_1, \cB_1)$, where 
\begin{equation}
\label{def_cB_1}
\cB_1=\det\RG(X, V')\otimes\det\RG(X, M)^{n-2}\otimes\det\RG(X, \cO)^{4-2n}
\end{equation}
and $\cB_1^2$ is identified with $\det\RG(X, M_1)$ via (\ref{iso_detRG_MotimesV'}).
 
\medskip\noindent
6.3 The stack $\Bun_{P\cap R}$ classifies exact sequences 
\begin{equation}
\label{seq_U'_by_U2}
0\to U_2\to U\to U'\to 0
\end{equation}
and (\ref{seq_cO_by_wedge2U}) on $X$ with $U'\in\Bun_{n-2}, U_2\in\Bun_2$. Write $\nu_{P,R}: \Bun_{P\cap R}\to \Bun_R$ and $\nu_{R,P}: \Bun_{P\cap R}\to\Bun_P$ for the natural maps. 
  
  We have a diagram
\begin{equation}
\label{diag_forR_two_squares}
\begin{array}{ccccc}
^w\ocY_R & \getsup{\pi_R} & ^w\cX_R & \toup{\tilde\rho_R} & \Bunt_{\GG_{2n-4}}\\
\uparrow && \uparrow && \uparrow\lefteqn{\scriptstyle \tilde\nu_{\PP}}\\
^w\ocY_R\times_{\Bun_R}\Bun_{P\cap R} & \getsup{\pi_R\times\id} & ^w\cX_R\times_{\Bun_R}\Bun_{P\cap R} & \toup{\nu_{\PP,R}} & \Bun_{\PP_{2n-4}},
\end{array}
\end{equation}
where the map $\nu_{\PP, R}$ is defined as follows. Given 
a collection 
\begin{equation}
\label{point_of_cXR_timesoverBunR_BunPR}
(U_2\subset U\subset V, s_2: U_2\to M, \,\det M\,\iso\,\Omega)\in {^w\cX_R\times_{\Bun_R}\Bun_{P\cap R}}
\end{equation}
let (\ref{seq_U'_by_M}) be the push-forward of $0\to U_2\to U\to U'\to 0$ by $s_2: U_2\to M$. Let $(M\subset \tilde V)\in\Bun_R$ be defined as in Section~6.2.2 and $(\Omega\subset M_1)$ be the symplectic-Heisenberg bundle associated to $(M\subset\tilde V)$. Then $\bar\cL=(M\otimes\tilde U)/\Sym^2 M$ is a lagrangian subbundle in $M_1$, it fits in the exact sequence (\ref{seq_cL_by_Omega}) with $\cL=M\otimes U'$. One checks that the element of $\Ext^1(\cL,\Omega)\,\iso\,\Ext^1(U',M)$ corresponding to (\ref{seq_cL_by_Omega}) is given by (\ref{seq_U'_by_M}). By definition, $\nu_{\PP, R}$ sends (\ref{point_of_cXR_timesoverBunR_BunPR}) to $(\Omega\subset \tilde\cL\subset M_1)\in\Bun_{\PP_{2n-4}}$. 

\begin{Lm} The right square in (\ref{diag_forR_two_squares}) is canonically 2-commutative.
\end{Lm}
\begin{Prf}  
For a point (\ref{point_of_cXR_timesoverBunR_BunPR}) let $\cB=\det\RG(X, M\otimes U')$ and let $\cB_1$ be defined by (\ref{def_cB_1}). Recall that $V'=V_{-2}/U_2$, where $V_{-2}$ is the orthogonal complement of $U_2$ in $V$. We must upgrade the natural isomorphisms
$\cB^2\,\iso\, \det\RG(X, M\otimes V')\,\iso\, \cB_1^2$ to a compatible isomorphism $\cB\,\iso\,\cB_1$. 

 By (\cite{L3}, Lemma~1), there is a canonical $\ZZ/2\ZZ$-graded isomorphism 
$$
\cB\,\iso\, \frac{\det\RG(X, M)^{n-2}\otimes\det\RG(X, U')^2\otimes\det\RG(X, \Omega\otimes\det U')}
{\det\RG(X, \det U')\otimes\det\RG(X, \cO)^{2n-4}}
$$
Applying (\cite{L3}, Lemma~1) to the exact sequence $0\to U'\to V'\to U'^*\to 0$, we get the isomorphisms
\begin{multline*}
\det\RG(X, V')\,\iso\, \det\RG(X, U')\otimes\det\RG(X, U'^*)\,\iso\,\det\RG(X, U')\otimes\det\RG(X, U'\otimes\Omega)\,\iso\\
\det\RG(X, U')^2\otimes\det\RG(X, \Omega\otimes\det U')\otimes\det\RG(X, \det U')^{-1}
\end{multline*}  
They yield the desired isomorphism $\cB\,\iso\, \cB_1$.  
\end{Prf}  
  
\medskip

 Recall the stack $^{\flat}(\Bun_2\times\Bun_{P_{n-2}})$ from Section~2.3.3. Write $^{\flat}\Bun_{P\cap R}$ for the preimage of $^{\flat}(\Bun_2\times\Bun_{P_{n-2}})$ under the map $\Bun_{P\cap R}\to \Bun_2\times\Bun_{P_{n-2}}$ sending $(U_2\subset U\subset V)$ to $(U_2, U'\subset V')$.

 The restrictions of $\nu_{P,R}$ and of $\nu_{R,P}$ to $^{\flat}\Bun_{P\cap R}$ are smooth. For a point (\ref{point_of_cXR_timesoverBunR_BunPR}) the conditions (\ref{conditions_flat}) imply $\H^0(X, \Sym^2(M\otimes U'))=0$, so $\nu_{\PP,R}$ restricts to a morphism
$$
\nu_{\PP,R}: {^w\cX_R\times_{\Bun_R}{^{\flat}\Bun_{P\cap R}}}
\to{^0\Bun_{\PP_{2n-4}}}
$$

  To prove Proposition~\ref{Pp_perversity_for_PH-model}, we  establish an explicit formula for the restriction of (\ref{perv_sheaf_over_c_ocY_PH}) under the smooth projection $\pr_1: {^w\ocY_R}\times_{\Bun_R}{^{\flat}\Bun_{P\cap R}}\to{^w\ocY_R}$. By definition of $^w\Bun_R$, the latter map is surjective.
  
  Recall the stack $\cS_P$ from Section~2.3.1. The stack 
$\Bun_{P\cap R}\times_{\Bun_n}\cS_P$ classifies $(U_2\subset U\subset V)\in \Bun_{P\cap R}$ and a section $s: U\to M$ with $M\in\Bun_{G_1}$. Let $\cW_R\subset \Bun_{P\cap R}\times_{\Bun_n}\cS_P$ be the open substack given by the condition that the composition $U_2\hook{} U\toup{s} M$ is an inclusion of coherent $\cO_X$-modules (this composition is denoted $s_2$). 

 Write $\bar\cW_R\subset \Bun_{P\cap R}\times_{\Bun_n}\cY_P$ for the open substack classifying $(U_2\subset U\subset V)\in\Bun_{P\cap R}$ and $v: \wedge^2 U\to\Omega$ such that the composition $\wedge^2 U_2\hook{}\wedge^2 U\toup{v}\Omega$ is non zero (this composition is denoted $v_2$). Let $\pi_{\cW}: \cW_R\to \bar\cW_R$ be the morphism over $\Bun_{P\cap R}$ given by $v=\wedge^2 s$.
Write $\ev_{\bar\cW_R}: \bar\cW_R\to\A^1$ for the map sending the above point to the pairing of $v$ with the exact sequence (\ref{seq_cO_by_wedge2U}) defining $V$. We get a diagram
$$
\begin{array}{ccc}
\cW_R\;\toup{\pi_{\cW}} &\bar\cW_R & \toup{p_{\bar\cW}}\;\ocY_R\times_{\Bun_R}\Bun_{P\cap R}\\
& \downarrow\lefteqn{\scriptstyle \pr_{\cY}}\\
&\ocY_P
\end{array}
$$
where $p_{\bar W}$ sends a collection $(U_2\subset U\subset V, v)\in\bar\cW_R$ to $(U_2\subset U\subset V, v_2: \wedge^2 U_2\to\Omega)$, here $v_2$ is the restriction of $v$ to $\wedge^2 U_2$. We have denoted by $\pr_{\cY}: \bar\cW_R\to \ocY_P$ the projection sending the above point to $(U,v)$. 
  
\begin{Pp}
\label{Pp_restricting_17}
Over ${^w\ocY_R}\times_{\Bun_R}{^{\flat}\Bun_{P\cap R}}$ the complex 
\begin{equation}
\label{complex_restriction_of_17}
\pr_1^*(\pi_{R})_!\tilde\rho_{R}^*\Aut^e_{\psi}\otimes(\Qlb[1](\frac{1}{2}))^{\dimrel(\tilde\rho_R)+\dimrel(\pr_1)}
\end{equation}
identifies with 
\begin{equation}
\label{complex_RHS_for_62}
p_{\bar\cW !}(\ev_{\bar\cW_R}^*\cL_{\psi}\otimes \pr_{\cY}^*\IC(\cZ_P))\otimes(\Qlb[1](\frac{1}{2}))^{\dimrel(\pr_{\cY})}
\end{equation}
\end{Pp}
\begin{Prf}
By Proposition~\ref{Pp_explicit_formula_Aute}, diagram (\ref{diag_forR_two_squares}) yields an isomorphism between (\ref{complex_restriction_of_17}) and 
$$
(\pi_R\times\id)_! \nu_{\PP,R}^*K_{\PP_{2n-4},\psi}\otimes(\Qlb[1](\frac{1}{2}))^{\dimrel(\nu_{\PP,R})}
$$
over ${^w\ocY_R}\times_{\Bun_R}{^{\flat}\Bun_{P\cap R}}$.
By definition of $K_{\PP_{2n-4},\psi}$, the latter complex identifies canonically with 
$$
p_{\bar \cW !}((\pi_{\cW !}\Qlb)\otimes\ev_{\bar\cW_R}^*\cL_{\psi})\otimes(\Qlb[1](\frac{1}{2}))^{\dim \cW_R}
$$
By Proposition~\ref{Pp_small_map}, we have $\pi_{\cW !}(\Qlb[1](\frac{1}{2}))^{\dim \cW_R}\,\iso\, \pr_{\cY}^*\IC(\cZ_P)\otimes(\Qlb[1](\frac{1}{2}))^{\dimrel(\pr_{\cY})}$. We are done.
\end{Prf}

\medskip

 Write $^{w\flat}\Bun_{P\cap R}$ for the preimage of $^w\Bun_R$ under $\nu_{P,R}: {^{\flat}\Bun_{P\cap R}\to \Bun_R}$. As in Section~6.1.2, one defines the category $\P^W(\cY_R\times_{\Bun_R} {^{w\flat}\Bun_{P\cap R}})$ of $(\cY_{\bar R}, \ev_R^*\cL_{\psi})$-equivariant perverse sheaves on $\cY_R\times_{\Bun_R} {^{w\flat}\Bun_{P\cap R}}$, and the corresponding derived category 
$$
\D^W(\cY_R\times_{\Bun_R} {^{w\flat}\Bun_{P\cap R}})
$$ 
As for (\ref{Four_R_psi-1}), one defines an equivalence exact for the perverse t-structures (denoted by the same symbol by a slight abuse of notations)
$$
\Four_{R, \psi^{-1}}: \D^{\prec}(^{w\flat}\Bun_{P\cap R})\to
\D^W(\cY_R\times_{\Bun_R} {^{w\flat}\Bun_{P\cap R}})
$$
The functor $(f_R\times\id)_!: \D^W(\cY_R\times_{\Bun_R} {^{w\flat}\Bun_{P\cap R}})\to \D^{\prec}(^{w\flat}\Bun_{P\cap R})$ is quasi-inverse to the latter equivalence.

\begin{Pp} 
\label{Pp_twelve_twopoints}
1) The complex (\ref{complex_RHS_for_62}) is canonically isomorphic to the restriction of 
\begin{equation}
\label{one_more_perv_sheaf}
\Four_{R,\psi^{-1}}\nu_{R,P}^*K_{P,\psi}\otimes(\Qlb[1](\frac{1}{2}))^{\dimrel(\nu_{R,P})}
\end{equation}
to the open substack ${^w\ocY_R}\times_{\Bun_R}{^{\flat}\Bun_{P\cap R}}\hook{} \cY_R\times_{\Bun_R} {^{w\flat}\Bun_{P\cap R}}$. Here $\nu_{R,P}: {^{w\flat}\Bun_{R\cap P}\to \Bun_P}$ is the natural map.\\
2) Over each connected component of $\cY_R\times_{\Bun_R} {^{w\flat}\Bun_{P\cap R}}$ 
the complex (\ref{one_more_perv_sheaf}) is an irreducible perverse sheaf.
\end{Pp}
\begin{Prf} 1) This follows formally from the properties of Fourier transforms. Indeed, one calculates the Fourier transform over the vector space $\Hom(\wedge^2 U,\Omega)$ composed with the backwords Fourier transform over $\H^1(X, \wedge^2 U_2)$.\\
2) The map $\nu_{R,P}: {^{w\flat}\Bun_{R\cap P}\to \Bun_P}$ is smooth, and $\nu_{R,P}^*K_{P,\psi}[\dimrel(\nu_{R,P})]$ is an irreducible perverse sheaf on each connected component of $^{w\flat}\Bun_{R\cap P}$.
\end{Prf}
  
\bigskip\noindent
\begin{Prf}\select{of Proposition~\ref{Pp_perversity_for_PH-model}}

\medskip\noindent
Since $\pr_1: {^w\ocY_R}\times_{\Bun_R}{^{\flat}\Bun_{P\cap R}}\to{^w\ocY_R}$ is smooth and surjective, it suffices to show that (\ref{complex_restriction_of_17}) is perverse. This follows by combining Propositions~\ref{Pp_restricting_17} and \ref{Pp_twelve_twopoints}.
\end{Prf}  

\begin{Cor} 
\label{Cor_PR_compatibility}
Over $^{w\flat}\Bun_{P\cap R}$ there exists an isomorphism
$$
\nu_{R,P}^*K_{P,\psi}\otimes(\Qlb[1](\frac{1}{2}))^{\dimrel(\nu_{R,P})}\,\iso\,
\nu_{P,R}^*K_{R,\psi}\otimes(\Qlb[1](\frac{1}{2}))^{\dimrel(\nu_{P,R})}
$$
\end{Cor}
\begin{Prf} Combine Propositions~\ref{Pp_restricting_17}, \ref{Pp_twelve_twopoints} and use the irreducibility of (\ref{one_more_perv_sheaf}) on each connected component of $\cY_R\times_{\Bun_R} {^{w\flat}\Bun_{P\cap R}}$.
\end{Prf}

\bigskip

\centerline{\scshape 7. The perverse sheaf $\cK_H$}

\bigskip\noindent
7.1 Note that the results of Section~4 hold over a suitable finite subfield of $k$, in particular the perverse sheaves $_a\cK^d_H$ admit a Weyl structure for this finite subfield of $k$.
In Sections~7.1-7.4 we assume that the ground field is $k=\Fq$.  

Recall the stack $_a\Bun_{G_1}$ defined in Section~2.3.4. For $a<\min\{2g-2, 0\}$ let $_{un,a}\Bun_{G_1}$ be the stack classifying a line bundle $L$ with $\deg(L^*\otimes\Omega)=a$ and an exact sequence $0\to L\to M\to L^*\otimes\Omega\to 0$ on $X$. The map $_{un,a}\Bun_{G_1}\to 
\Bun_{G_1}$ sending the above point to $M$ is a locally closed immersion. Moreover, if $a<\min\{2g-2, 0\}$ then $\mathop{\cup}\limits_{b\le a}  {_{un, b}\Bun_{G_1}}$ is a stratification of $\Bun_{G_1}-\,{_a\Bun_{G_1}}$. 

 Write $_a\cU_H\subset\Bun_H$ for the open substack of $V\in\Bun_H$ such that for any $L\in\Bun_1$ with $\deg L\le a$ one has $\Hom(V, L)=0$. The stack $_a\cU_H$ is of finite type.
 
\begin{Lm} 
\label{Lm_restricting_Aut_GH}
For $a<\min\{2g-2, 0\}$ the $*$-restriction of $\Aut_{G_1,H}$ to $_{un, a}\Bun_{G_1}\times {_a\cU_H}$ identifies with 
$$
\Qlb[\dim\Bun_{G_1,H}-2n(g-1-a)]
$$
\end{Lm}
\begin{Prf}
Apply Proposition~\ref{Pp_explicit_from_Wald_periods} or (\cite{L1}, Theorem~1). For a point $(L\subset M, V)\in {_{un, a}\Bun_{G_1}\times {_a\cU_H}}$ we get $\H^0(X, V\otimes L^*\otimes\Omega)=0$ and $\H^0(X, M\otimes V)=\H^0(X, L\otimes V)$ is of dimension $2n(g-1-a)$.  
\end{Prf}

\medskip

 Lemma~\ref{Lm_restricting_Aut_GH} immediately yields the following. 
\begin{Cor}
\label{Cor_modulo_constant_complexes}
If $a<\min\{2g-2, 0\}$ then the cone of the natural map $_a\tilde K\to {_{a-1}\tilde K}$ over $_a\cU_H$
is a constant complex. \QED
\end{Cor} 

\medskip\noindent
7.2 For $b\in \ZZ/2\ZZ$ set $_a\cU_H^b={_a\cU_H}\cap \Bun_H^b$. Note that if $_a\Bun_n^d$ is not empty then $a<(d/n)+g$. So, one can not find $a\in\ZZ$ such that $^e_a\ocZ{}^d_P$ is not empty for all $d\in Z(e,P)$.

 Let $a$ be small enough so that the function $E_{\cK}$ defined in Proposition~\ref{Pp_Euler_char_pointwise} does not vanish over $_a\cU_H^b$ for each $b\in\ZZ/2\ZZ$. Set $\cU_H=
{_a\cU_H}$.

\begin{Lm} 
\label{Lm_something_finite}
The set of $d\in Z(e,P)$ such that $\cK^d_H$ vanishes over $\cU_H$ is at most finite. 
\end{Lm}
\begin{Prf}
Let $d\in Z(e,P)$, pick $a'$ such that $^e_{a'}\ocZ{}^d_P$ is not empty. The irreducible subquotient $\cK^d_H$ of $^p\H^0(_{a'}\tilde K)$ introduced in Definition~\ref{Def_cK^d_H} is characterised by the following property.
The perverse sheaf $\nu_P^*(\cK^d_H)[\dimrel(\nu_P)]$ 
over $^e\Bun_P^d$ contains the irreducible subquotient $K_{P,\psi}^d$.  

 If $\cK^d_H$ vanishes over $\cU_H$ then $K^d_{P,\psi}$ would vanish over $\nu_P^{-1}(\cU_H)\cap {^e\Bun_P^d}$. In particular, $E_{\cK}$ would vanish over $\cU_H\cap \nu_P(^e\Bun_P^d)$. 

 Let $\cI\subset Z(e,P)$ be the set of those $d$ for which $\cK^d_H$ vanishes over $\cU_H$. If $\cI$ is infinite then the union of $\cU_H\cap \nu_P(^e\Bun_P^d)$, $d\in\cI$ equals $\cU_H$, and $E_{\cK}$ would vanish over $\cU_H$. This contradiction shows that $\cI$ is finite.
\end{Prf}   

\medskip

 Using Lemma~\ref{Lm_something_finite}, we replace if necessary $a$ by a smaller integer and assume from now on that for all $d\in Z(e,P)$ the perverse sheaf $\cK^d_H$ does not vanish over $\cU_H$. We also assume $a<\min\{2g-2, 0\}$. Set $\tilde\cK_{\cU}={^p\H^0(_a\tilde K)}\mid_{\cU_H}$. 
 
\begin{Lm}
\label{Lm_already_appears}
For each $d\in Z(e,P)$ the perverse sheaf $\cK^d_H$ already appears as an irreducible subquotient in $\tilde\cK_{\cU}$. More precisely, let $a'\le a$ be such that 
$^e_{a'}\ocZ{}^d_P$ is not empty. Then there is a unique irreducible subquotient $\cK^d_{\cU}$ in $\tilde\cK_{\cU}$ such that 
$$
\alpha: {^p\H^0(_a\tilde K)}\to {^p\H^0(_{a'}\tilde K)}
$$ 
induces an isomorphism $\cK^d_{\cU}\,\iso\, {_{a'}\cK^d_H}$.
The subquotient $\cK^d_{\cU}$ of $\tilde\cK_{\cU}$ is characterised by the property that 
$$
\nu_P^*(\cK^d_{\cU})[\dimrel(\nu_P)]
$$ 
over $^e\Bun_P^d\cap \nu_P^{-1}(\cU_H)$ contains $K^d_{P,\psi}$ as an irreducible subquotient.
\end{Lm}
\begin{Prf}
By Corollary~\ref{Cor_modulo_constant_complexes}, the kernel and cokernel of $\alpha: {^p\H^0(_a\tilde K)}\to {^p\H^0(_{a'}\tilde K)}$ over $\cU_H$ are perverse sheaves, which are succesive extensions of constant perverse sheaves. Since  $\cK^d_H$ is not constant and does not vanish over $\cU_H$, our assertion follows.
\end{Prf}
    
\medskip\noindent
7.3.1 Let $F$ be an irreducible subquotient in $\tilde\cK_{\cU}$. Let $I_F\subset Z(e,P)$ be the set of $d\in Z(e,P)$ such that $F$ does not vanish over $\cU_H\cap \nu_P(^e\Bun_P^d)$. The set $I_F$ is infinite. Write $\bar F$ for the intermediate extension of $F$ under $\cU_H\hook{}\Bun_H$.

 Let $\nu_n: \Bun_P\to \Bun_n$ be the map sending (\ref{seq_cO_by_wedge2U}) to $U$. The morphism $\nu_n$ is smooth.
 
\begin{Lm} 
\label{Lm_two_cases}
For each $d\in Z(e,P)$ the perverse sheaf
$$
\Four_{\cY_P,\psi}^{-1}\nu_P^*(\bar F)\otimes(\Qlb[1](\frac{1}{2}))^{\dimrel(\nu_P)}
$$
over $^e\ocY{}^d_P$ either vanishes or identifies with $\IC(\cZ_P)$. In the first case there is a perverse sheaf $\cF^d\in \P(^e\Bun_n^d)$ and an isomorphism 
\begin{equation}
\label{iso_descent_on_Bun_P^d}
\nu_P^*(\bar F)\otimes(\Qlb[1](\frac{1}{2}))^{\dimrel(\nu_P)}
\,\iso\,
\nu_n^*\cF^d\otimes(\Qlb[1](\frac{1}{2}))^{\dimrel(\nu_n)}
\end{equation}
over $^e\Bun_P^d$. In the second case $F=\cK^d_{\cU}$.   
\end{Lm}
\begin{Prf}
We may assume $\bar F$ non constant. Let $a'\le a$ and $S={^p\H^0(_{a'}\tilde K)}$. By Corollary~\ref{Cor_modulo_constant_complexes}, 
the image of $\bar F$ under (\ref{morphism_alpha_from_Def3}) is a nonzero irreducible subquotient in $S$. By Corollary~\ref{Cor_passage_to_pH^0_first},
$$
\Four_{\cY_P,\psi}^{-1}\nu_P^*(S)\otimes(\Qlb[1](\frac{1}{2}))^{\dimrel(\nu_P)}\,\iso\, \IC(\cZ_P)
$$
over $^e_{a'}\ocY{}^d_P$. Since the union of $^e_{a'}\ocY{}^d_P$ for all $a'\le a$ equals $^e\ocY{}^d_P$, we are done.
\end{Prf}

\medskip

 In Appendix~A we introduce a notion of an almost constant local system on $\Bun_H$. Note that if $E$ is an irreducible almost constant local system on $\Bun_H^b$ for some $b\in\ZZ/2\ZZ$ then $E$ is of rank one and order at most two. The following will be proved in Section~7.3.2.
 
\begin{Pp} 
\label{Pp_K^d_coincide}
The irreducible subquotients $\cK^d_{\cU}$ of $\tilde\cK_{\cU}$ over $\cU_H^b$ all coincide for $d\!\!\mod 2=b$. The resulting irreducible subquotient is denoted $\cK_{\cU,b}$. If $F$ is a different irreducible subquotient of $\tilde \cK_{\cU}$ over $\cU^b_H$ then $\bar F\otimes_k \bar k$ is a direct sum of (shifted) almost constant local systems on $\Bun_H^b$.
\end{Pp}

\begin{Def} 
\label{Def_true_cK_H}
The perverse sheaf $\cK_H\in \P(\Bun_H)$ is defined as the the intermediate extension of $\cK_{\cU,0}\oplus\cK_{\cU,1}$ under $\cU_H\hook{}\Bun_H$. The perverse sheaf $\cK_H$ is irreducible over each connected component of $\Bun_H$. 
\end{Def}
  
 Proposition~\ref{Pp_K^d_coincide} immediately implies the following.
  
\begin{Cor}
\label{Cor_properties_of_cK_H}
 For each $d\in Z(e, P)$ the perverse sheaf $\nu_P^*(\cK_H)\otimes(\Qlb[1](\frac{1}{2}))^{\dimrel(\nu_P)}$ over $^e\Bun_P^d$ contains $K^d_{P,\psi}$ as an irreducible subquotient. More precisely, for each $d\in Z(e,P)$ there is an
isomorphism
$$
\Four^{-1}_{\cY_P,\psi}\nu^*_P(\cK_H)\otimes(\Qlb[1](\frac{1}{2}))^{\dimrel(\nu_P)}\,\iso\, \IC(\cZ_P)
$$
over $^e\ocY{}^d_P$.
\end{Cor} 
  
\medskip\noindent  
7.3.2  The fact that $k$ is finite will be used in the proof of the following key lemma. 
  
\begin{Lm} 
\label{Lm_really_key_lemma}
Let $\bar F$ be an irreducible perverse sheaf on $\Bun_H^b$ for some $b\in\ZZ/2\ZZ$. Let $I$ be an infinite bounded from above set of integers. Assume given for each $d\in I$ a perverse sheaf $\cF^d\in \P(^e\Bun_n^d)$ and an isomorphism (\ref{iso_descent_on_Bun_P^d}) over $^e\Bun_P^d$. Assume that if $d\in I$ then $\nu_P^*(\bar F)$ is nonzero over $^e\Bun_P^d$. Then each irreducible subquotient of $\bar F\otimes_{k} \bar k$ is a (shifted) almost constant local system on $\Bun_H^b$.
\end{Lm}
  
\begin{Rem} We do not require in Lemma~\ref{Lm_really_key_lemma} that $\cF^d$ are irreducible. We can not garantee this, as we don't know if the geometric fibres of $\nu_P: {^e\Bun_P^d}\to\Bun_H$ are connected (for generic fibres cf. Proposition~\ref{Pp_fibres_of_nu_P_connected}).
\end{Rem}
  
\medskip 
  
\begin{Prf}\select{of Proposition~\ref{Pp_K^d_coincide}} \   The perverse sheaf $\tilde\cK_{\cU}\mid_{\cU_H^b}$ admits at least one irreducible subquotient which is not an almost constant (shifted) local system. Let $F$ be such an irreducible subquotient. Then by Lemma~\ref{Lm_really_key_lemma}, the set $A=\{d\in I_F\mid F\ne \cK^d_{\cU}\}$ is finite. Let $F'$ be an irreducible subquotient of $\tilde\cK_{\cU}\mid_{\cU_H^b}$ not equal to $F$. Then for any $d\in I_F-A$ we get $F'\ne \cK^d_{\cU}$. So, by Lemmas~\ref{Lm_really_key_lemma} and \ref{Lm_two_cases}, each irreducible subquotient of
$\bar F'\otimes_k \bar k$ is a (shifted) almost constant local system on $\Bun_H^b\otimes_k \bar k$. This implies that $F'\ne \cK^d_{\cU}$ for all $d\in Z(e,P)$. Thus, $F=\cK^d_{\cU}$ for all $d\in Z(e,P)$.
\end{Prf}

\bigskip\noindent
7.4 \select{Proof of Lemma~\ref{Lm_really_key_lemma}}

\medskip\noindent
For $d_1, d_2\in\ZZ$ of the same parity write $\cX^{d_1, d_2}\subset \Bun_P^{d_1}\times_{\Bun_H}\Bun_P^{d_2}$ for the open substack given by the property that the two $P$-structures on $V\in\Bun_H$ are transversal at the generic point of $X$. So, $\cX^{d_1, d_2}$ classifies two exact sequences $0\to\wedge^2 U_i\to ?\to \cO_X$ giving rise to $0\to U_i\to V\to U_i^*\to 0$ such that the composition $U_1\to V\to U_2^*$ is an inclusion of coherent $\cO_X$-modules, and the isomorphisms 
$$
\det V\,\iso\, (\det U_i)\otimes \det U_i^*\,\iso\,\cO
$$ 
coincide for $i=1,2$.  
 
 For a point of $\cX^{d_1, d_2}$ we get a diagram $U_1\oplus U_2\subset V\subset U_2^*\oplus U_1^*$. The projections $V/(U_1\oplus U_2)\to U_2^*/U_1$ and 
$V/(U_1\oplus U_2)\to U_1^*/U_2$ are isomorphisms, so there is an isomorphism 
$$
\phi: U_1^*/U_2\,\iso\, U_2^*/U_1
$$
of torsion sheaves on $X$ such that 
$
V/(U_1\oplus U_2)=\{(v, \phi(v))\in (U_1^*/U_2)\oplus (U_2^*/U_1)\mid v\in U_1^*/U_2\}
$.  
Moreover, $\phi$ is anti-symmetric in the sense that for any $v_1, v_2\in U_1^*/U_2$ one has 
\begin{equation}
\label{anti-symmetric_description}
\<v_1, \phi(v_2)\>+\<\phi(v_1), v_2\>\in \cO_X
\end{equation}
Here $\<\cdot,\cdot\>$ is the natural pairing.

\begin{Rem} Write $\cO_x$ for the completed local ring of $X$ at $x\in X$, let $t_x\in\cO_x$ be a uniformizer. Assume that 
$a_1\ge\ldots\ge a_m>0$ and 
$$
\phi: \cO_x/t_x^{a_1}\oplus\ldots\oplus \cO_x/t_x^{a_m}\to t_x^{-a_1}\cO_x/\cO_x\oplus\ldots\oplus t_x^{-a_m}\cO_x/\cO_x
$$
is a $\cO_X$-linear map given by a matrix $b=(b_{ij})$. Then (\ref{anti-symmetric_description}) holds iff $b_{ij}\in t_x^{-\min\{a_i, a_j\}}\cO_x/\cO_x$ and for all $i,j$ one has $b_{ij}+b_{ji}=0$. Since the characteristic of $k$ is not 2, this implies in particular $b_{ii}=0$.
\end{Rem}

\medskip

 Let $\wt\cX^{d_1, d_2}\subset \cX^{d_1,d_2}$ be the open substack given by the property that there is an effective reduced divisor $D\ge 0$ on $X$ such that $\div(U_1^*/U_2)=2D$. For a point of $\wt\cX^{d_1, d_2}$ there is an isomorphism $U_1^*/U_2\,\iso\, \cO_D\oplus\cO_D$. Here $\cO_D$ is the structure sheaf of $D$. We have a diagram of smooth projections
$$
\Bun_n^{d_1} \;\getsup{q_1} \; \wt\cX^{d_1, d_2}\;\toup{q_2}\;\Bun_n^{d_2},
$$
where $q_i$ sends the above point to $U_i$. 

 Write $^e\wt\cX^{d_1, d_2}\subset \wt\cX^{d_1,d_2}$ for the preimage of $^e\Bun_n^{d_1}\times {^e\Bun_n^{d_2}}$ under $q_1\times q_2$. Consider the diagram of projections
$$
^e\Bun_n^{d_1} \;\getsup{^eq_1} \; {^e\wt\cX^{d_1, d_2}}\;\toup{^eq_2}\;\Bun_n^{d_2}
$$ 
By our assumptions, for $d_1,d_2\in I$ there are isomorphisms $\sigma: {^eq_1^*}\cF^{d_1}\,\iso\, {^eq_2^*}\cF^{d_2}$ of shifted perverse sheaves over $^e\wt\cX^{d_1, d_2}$. 

 Write $\tilde\cF^d$ for the intermediate extension of $\cF^d$ under $^e\Bun_n^d\hook{} \Bun_n^d$. The stack $\wt\cX^{d_1, d_2}$ is irreducible. So,
if $^e\Bun_n^{d_i}$ is not empty for $i=1,2$ then $^e\wt\cX^{d_1,d_2}$ is dense in $\wt\cX^{d_1,d_2}$. Thus, the isomorphisms $\sigma$ extend (by the intermediate extension) to isomorphisms 
$$
\tilde\sigma: q_1^*\tilde\cF^{d_1}\,\iso\, q_2^*\tilde\cF^{d_2}
$$ 
of shifted perverse sheaves over $\wt\cX^{d_1,d_2}$. 
For $U_2\in \Bun_n^{d_2}(k)$ write 
$$
\wt\cX^{d_1, d_2}(U_2)=\wt\cX^{d_1,d_2}\times_{\Bun_n^{d_2}}\Spec k,
$$
where we used the map $U_2: \Spec k\to \Bun_n^{d_2}$ to define the fibred product. 
 
 Given $d_1$ and $U_1\in \Bun_n^{d_1}$, there is $d_2\in I$ sufficiently small and $U_2\in\Bun_n^{d_2}(k)$ such that the projection $q_1: \wt\cX^{d_1, d_2}(U_2)\to\Bun_n^{d_1}$ is smooth over a Zariski open neighbourhood of $U_1$. Now the isomorphism $\tilde\sigma$ shows that $\tilde\cF^{d_1}$ is smooth in a neighbourhood of $U_1$. Since $U_1$ was arbitrary, $\tilde\cF^d$ is a shifted local system on $\Bun_n^d$. 
The union of the images of $\nu_P: {^e\Bun_P^d}\to\Bun_H^b$ equals $\Bun_H^b$, so $\bar F$ is also a shifted local system over $\Bun_H^a$.

 Now Conjecture~\ref{Con_pi_1_of_BunG} would imply that $\cF$ is an almost constant local system. Conjecture~\ref{Con_pi_1_of_BunG} not being known, we give another argument that applies for the finite ground field $k$. 

 For $k$-points $U_i\in\Bun_n^{d_i}$ with $d_i\in I$ say that $U_1\prec U_2$ if there is a $k$-point $\eta\in \wt\cX^{d_1,d_2}$ such that $q_i(\eta)=U_i$ for $i=1,2$. Write $\sim$ for the equivalence relation generated by $\prec$. If two $k$-points $U_1,U_2\in\Bun_n^d$ are equivalent in this sense then for the maps $\kappa_i: \Spec k\toup{U_i}\Bun_n^d$ the isomorphisms $\tilde\sigma$ yield $\kappa_1^*\tilde\cF^d\;\iso\; \kappa_2^*\tilde\cF^d$. 

\begin{Lm}
\label{Lm_euqivalence_relation_Bun_n}
Assume $n\ge 3$. Let $d\in I$ and $U_i\in \Bun_n^d(k)$ for $i=1,2$. Then $U_1\sim U_2$ if and only if there is $\cE\in\Bun_1^0(k)$ with $\cE^2\,\iso\, (\det U_1)\otimes(\det U_2)^{-1}$.
\end{Lm}

 Let $d\in I$ and $L_1\in \Bun_1^d(k)$. Let $\Bun_{n, L_1}^d$ be the stack classifying $U_1\in\Bun_n^d$, $\cE\in\Bun_1^0$ and an isomorphism $\det U_1\,\iso\, \cE^2\otimes L_1$. 
 
 By Lemma~\ref{Lm_euqivalence_relation_Bun_n}, the $*$-restrictions of $\cF^d$ to all $k$-points of $\Bun_{n, L_1}^d$ are isomorphic to each other. In particular, the function trace of Frobenius
$\tr(\cF^d, k): \Bun_{n, L_1}^d(k)\to\Qlb$ is constant. Since the same hold for any finite extension of $k$, we conclude by (\cite{Lam}, Theorem~1.1.2) that $\cF^d$ is the inverse image of a local system on $\Spec k$.
 
 Let $\GSpin_{2n}$ be the quotient of $\Gm\times \Spin_{2n}$ be the diagonally embedded subgroup $A\,\iso\, \ZZ/2\ZZ$, here $H\,\iso\, \Spin_{2n}/A$. In terms of Appendix~A, we have taken $T=\Gm$ and $T_1=\Gm/A\,\iso\,\Gm$. Pick $\bar b\in\pi_1(\GSpin_{2n})$ over $b\in\pi_1(H)$, let $c\in \pi_1(T_1)$ be the image of $\bar b$. Pick a $T_1$-torsor $\cF_{T_1}$ in $\Bun_{T_1}^c(k)$. We get the stack $\Bun^{\bar b}_{\GSpin_{2n}, \cF_{T_1}}$ defined as in Appendix~A and the morphism
$$
f: \Bun^{\bar b}_{\GSpin_{2n}, \cF_{T_1}}\to\Bun_H^b
$$ 
Let $\bar P\subset \GSpin _{2n}$ be the preimage of $P$ under the natural map $\GSpin_{2n}\to H$.

 Set $\Bun_{\bar P, \cF_{T_1}}=\Bun_{\bar P}\times_{\Bun_{T_1}}\Spec k$, where we used the map $\cF_{T_1}: \Spec k\to\Bun_{T_1}$ to define the fibred product. There is a commutative diagram for a suitable $\bar d\in \pi_1(\bar P)$
$$
\begin{array}{ccccc}
\Bun_{n, L_1}^d & \gets & \Bun^{\bar d}_{\bar P, \cF_{T_1}} & \toup{\nu_{\bar P}} & \Bun^{\bar b}_{\GSpin_{2n}, \cF_{T_1}}\\
\downarrow && \downarrow\lefteqn{\scriptstyle f_P} && \downarrow\lefteqn{\scriptstyle f} \\
\Bun_n^d & \getsup{\nu_n} & \Bun_P^d & \toup{\nu_P} & \Bun_H^b
\end{array}
$$

 Let $^e\Bun^{\bar d}_{\bar P, \cF_{T_1}}$ be the preimage of $^e\Bun_P^d$ under $f_P$. We see that $\nu_{\bar P}^*f^*\bar F$ is the inverse image of a local system on $\Spec k$. By Proposition~\ref{Pp_appendix_B} in Appendix~B, for $d\in I$ small enough the generic fibre of $\nu_{\bar P}: {^e\Bun^{\bar d}_{\bar P, \cF_{T_1}}}\to \Bun^{\bar b}_{\GSpin_{2n}, \cF_{T_1}}$ is geometrically irreducible. So, $f^*\bar F$ is the inverse image of some local system over $\Spec k$.
Lemma~\ref{Lm_really_key_lemma} is reduced to Lemma~\ref{Lm_euqivalence_relation_Bun_n}. \QED

\medskip

Recall the following notion. Let $\lambda$ be a coweight of $\GL_n$ and $F_x$ the field of fractions of $\cO_x$, $x\in X$. If $L,L'$ are two free $\cO_x$-modules of rank $n$ with an isomorphism of generic fibres $\beta: L\otimes_{\cO_x} F_x\,\iso\, L'\otimes_{\cO_x}F_x$, we say that $L$ is in the position $\lambda$ with respect to $L'$ if there is a trivialization $\sigma: L'\,\iso\, \cO_x^n$ such that the image of 
$L\hook{} L\otimes F_x\toup{\beta} L'\otimes F_x\toup{\sigma} F_x^n$ equals $t_x^{\lambda}\cO_x^n$. 

\medskip

\begin{Prf}\select{of Lemma~\ref{Lm_euqivalence_relation_Bun_n}}

\medskip\noindent
Assume that $\cE^2\,\iso\, (\det U_1)\otimes(\det U_2)^{-1}$. We must prove that $U_1\sim U_2$.

 First, we may assume $\det U_1\,\iso\, \det U_2$. Indeed, by Bertini theorems (\cite{P}), there are reduced effective divisors $D^+, D^-$ on $X$ defined over $k$ such that for $D=D^+-D^-$ one has $\cE\,\iso\, \cO(D)$. Pick any $U_2\subset U^*$ and $U_3\subset U^*$ such that $U^*/U_2\,\iso\, \cO_{D^+}\oplus\cO_{D^+}$ and $U^*/U_3\,\iso\, \cO_{D^-}\oplus\cO_{D^-}$. We may assume $D^+, D^-$ sufficiently large so that $(\deg U)\in I$. Then $U_2\sim U_3$ and $\det U_1\,\iso\, \det U_3$. We are reduced to the case $\det U_1\,\iso\, \det U_2$.

  Pick $x\in X$ and an isomorphism $\gamma: U_1\,\iso\, U_2\mid_{X-x}$. One can find a sequence of $k$-points $U_3,\ldots, U_r\in \Bun_n^d$ and isomorphisms $\gamma_i: U_i\,\iso\, U_{i+1}\mid_{X-x}$ for $i=2,\ldots, r-1$ with $U_r=U_1$ such that $U_{i+1}$ is in the position $(1,0,\ldots, 0, -1)$ with respect to $U_i$ at $x$. 
  
  We are reduced to the case of an isomorphism $\gamma: U_1\,\iso\, U_2\mid_{X-x}$ such that $U_2$ is in the position $(1,0,\ldots, 0, -1)$ with respect to $U_1$ at $x$. This means that there is a base $\{e_1,\ldots, e_n\}$ of $U_1$ in a neighbourhood of $x$ such that $\{t_xe_1, e_2,\ldots, e_{n-1}, t^{-1}_xe_n\}$ is a base of $U_2$ in a neighbourhood of $x$. Here $t_x\in \cO_x$ is a uniformizer. Let $U'\in\Bun_n$ be the modification of $U_1$ whose local base in a neighbourhood of $x$ is $\{e_1,\ldots, e_{n-2}, t_x^{-1}e_{n-1}, t_x^{-1}e_n\}$. If $-d-2\in I$ then $U_1\sim U'^*\sim U_2$. Otherwise, replace $U'$ by a bigger suitable upper modification $U''$ at some points different from $x$ such that $U_1\sim U''^*\sim U_2$. We are done.
\end{Prf}

\medskip
\begin{Rem} 
\label{Rem_about_H^i_for_i_not_zero}
Let $i\ne 0$, let $F$ be an irreducible subquotient of $^p\H^i(_a\tilde K)\mid_{\cU_H}$. Write $\bar F$ for the intermediate extension of $F$ under $\cU_H\hook{}\Bun_H$. Then each irreducible subquotient of $\bar F\otimes_k \bar k$ is an almost constant local system.

 Indeed, if $F$ is not constant then, as in Lemma~\ref{Lm_already_appears}, we see that $F$ appears as an irreducible subquotient of $^p\H^i(_{a'}\tilde K)\mid_{\cU_H}$ for all $a'\le a$. This together with Corollary~\ref{Cor_passage_to_pH^0_first} implies that $\Four_{\cY_P, \psi}^{-1}\nu_P^*\bar F$ vanishes over the stack $^e\ocY_P$. Our claim follows now from Lemma~\ref{Lm_really_key_lemma}. Thus, the whole complex $_a\tilde K\mid_{\cU_H}$ is built up from $\cK_H$ and almost constant local systems.
\end{Rem} 

\medskip\noindent
7.5 Assume $k$ algebraically closed. 
Our purpose now is to establish more properties of the sheaf $\cK_H$. From Proposition~\ref{Pp_appendix_B} of Appendix~B one easily derives the following.
 
\begin{Pp} 
\label{Pp_fibres_of_nu_P_connected}
There is $N_0\in \ZZ$ such that for all $d\le N_0$
the generic fibre of $\nu_P: {^e\Bun_P^d}\to \Bun_H^{d\!\mod 2}$ is geometrically irreducible and non empty. \QED
\end{Pp} 

\begin{Prf}\select{of Theorem~\ref{Th1}}

\medskip\noindent
By Corollary~\ref{Cor_properties_of_cK_H}, for each $d\in Z(e,P)$ there exists a semi-simple perverse sheaf $\cM^d$ on $^e\Bun_n^d$ and an isomorphism over $^e\Bun_P^d$ 
\begin{equation}
\label{iso_ss_important}
\left(\nu_P^*(\cK_H)\otimes(\Qlb[1](\frac{1}{2}))^{\dimrel(\nu_P)}\right)^{ss}\,\iso\, K^d_{P,\psi}\oplus (\nu_n^*\cM^d\otimes (\Qlb[1](\frac{1}{2}))^{\dimrel(\nu_n)})
\end{equation}
Here the upper index $ss$ stands for the semisimplification of the corresponding perverse sheaf.
Now Proposition~\ref{Pp_Euler_char_pointwise} shows that there exists a function $E_{\cM}: \Bun_H(k)\to\ZZ$ such that for each $d\in Z(e,P)$ and $\eta\in {^e\Bun_P^d(k)}$ over $V\in\Bun_H(k)$ one has  
$$
\chi(\cM^d\mid_{\nu_n(\eta)})=(-1)^{\dimrel(\nu_n^d)+\dimrel(\nu_P^d)}E_{\cM}(V)
$$
 
 Assume that $E_{\cM}$ is not identically zero on $\Bun_H^b$ for some $b\in\ZZ/2\ZZ$. Pick $d_2\in Z(e,P)$ with $d_2\! \mod 2=b$ and $U_2\in {^e\Bun_n^d}$ such that $\chi(\cM^{d_2}\mid_{U_2})\ne 0$. Argue as in the proof of Lemma~\ref{Lm_really_key_lemma}. Given $d_1\in Z(e,P)$ with $d_1\!\mod 2=b$, consider the stack $\wt\cX^{d_1,d_2}(U_2)$ introduced in Section~7.4. Let $^e\wt\cX^{d_1,d_2}(U_2)$ be the preimage of $^e\wt\cX^{d_1,d_2}$ in $\wt\cX^{d_1,d_2}(U_2)$. If $d_1$ is sufficiently small, the projection $q_1: {^e\wt\cX^{d_1,d_2}(U_2)}\to{^e\Bun_n^{d_1}}$ is dominant, so that $\chi(\cM^{d_1}\mid_{U_1})=\chi(\cM^{d_2}\mid_{U_2})$ for $U_1$ lying in some nonempty open substack of $^e\Bun_n^{d_1}$. Since $^e\Bun_P^{d_1}\to\Bun_H^b$ is dominant, we conclude that $E_{\cM}$ does not vanish over some nonempty open substack of $\Bun_H^b$. This implies that $\cK_H$ does not vanish at the generic point of $\Bun_H^b$. Then applying (\cite{G}, Lemma~4.8) together with Proposition~\ref{Pp_fibres_of_nu_P_connected}, we learn that $\nu_P^*(\cK_H)[\dimrel(\nu_P)]$ is an irreducible perverse sheaf on $^e\Bun_P^{d_1}$, so $\cM^{d_1}$ must vanish. This contradiction shows that $E_{\cM}$ is identically zero.

 Since $\cM^d$ is a perverse sheaf, this in turn implies that $\cM^d=0$ for all $d\in Z(e,P)$. So, for each $d\in Z(e,P)$ the perverse sheaf $\nu_P^*(\cK_H)[\dimrel(\nu_P)]$ is  irreducible over $^e\Bun_P^d$. 
\end{Prf}

\medskip\noindent
7.6 \select{Proof of Theorem~\ref{Th_2}}

\medskip\noindent
\Step 1 Set $G=G_1$ for brevity. Let $_aE=\Qlb[\dim\Bun_G]$ over $_a\Bun_G$. Recall that $_a\tilde K=F_H(_aE)$, where 
$$
F_H: \D^-(\Bun_G)_!\to \D^{\prec}(\Bun_H)
$$ 
is given by (\cite{L2}, Definition 2). Recall that for $a<\min\{2g-2,0\}$ we have the locally closed substack $_{un, a}\Bun_G\subset \Bun_G$ introduced in Section~7.1. Set 
$
_aR=\Qlb[\dim\Bun_G]
$ 
over $_{un, a}\Bun_G$. Write $W$ for the standard representation of $\check{H}$. Let $W_1$ denote the standard representation of $\check{G}\,\iso\, \SO_3$ and $W_0=\oplus_{i=2-n}^{n-2}\Qlb[2i]$. By (\cite{L2}, Theorem 3), one has
\begin{equation}
\label{iso_from_paper_L1}
_x\H^{\la}_H(W, {_a\tilde K})\,\iso\,  F_H(_x\H^{\la}_{G}(W_0\oplus W_1, {_aE}))\,\iso\,
F_H(_x\H^{\la}_{G}(W_1, {_aE}))\oplus (W_0\otimes {_a\tilde K})
\end{equation}

 Recall that $_a\Bun_{G}\subset {_{a-1}\Bun_{G}}$. If $a+1<\min\{2g-2, 0\}$ then $_{a-1}\Bun_{G}$ admits the stratification
$$
_{a-1}\Bun_{G}={_{a+1}\Bun_{G}}\sqcup {_{un, a}\Bun_{G}}\sqcup  {_{un, a+1}\Bun_{G}},
$$ 
and the substack $_{un, a}\Bun_{G}$ is closed in $_{a-1}\Bun_{G}$. Set 
$\bar W_1=(\Qlb[-2]\oplus \Qlb\oplus \Qlb[2])$. 
The following Lemma is straightforward.

\begin{Lm} 
\label{Lm_Hecke_action_for_SL_2}
Let $a+1<\min\{2g-2, 0\}$. The complex $_x\H^{\la}_G(W_1, {_aE})$ is the extension by zero from $_{a-1}\Bun_{G}$.\\
1) The $*$-restriction of $_x\H^{\la}_G(W_1, {_aE})$ to $_{a+1}\Bun_{G}$ identifies with 
$$
\bar W_1\otimes (_{a+1}E),
$$
where $\bar W_1=(\Qlb[-2]\oplus \Qlb\oplus \Qlb[2])$. \\
2) The $*$-restriction of $_x\H^{\la}_G(W_1, {_aE})$ to $_{un, a}\Bun_{G}$ is 
$_aR[-2]$.\\
3) The $*$-restriction of $_x\H^{\la}_G(W_1, {_aE})$ to $_{un, a+1}\Bun_{G}$ fits into an exact triangle
$$
_aR[-2]\to {_x\H^{\la}_G(W_1, {_aE})}\to {_aR}
$$
\QED
\end{Lm}

It suffices to prove that there is a complex $L\in \D(\Bun_H)$, which is a finite direct sum of shifted almost constant local systems on $\Bun_H$, and an isomorphism in $\D(\Bun_H)$
\begin{equation}
\label{iso_Hecke_for_W_main}
\H^{\la}_H(W, \cK_H)\,\iso\, (\bar W_1+W_0)\otimes \cK_H +L
\end{equation}

 Let $b$ be any integer small enough, so that $\cK_H\mid_{_b\cU_H}\ne 0$. We will show that there is such $L$ (depending eventually on $b$) and an isomorphism (\ref{iso_Hecke_for_W_main}) over $_b\cU_H$. Since $b$ is arbitrary, $L$ is independent of $b$, and this would  conclude the proof.
 
\medskip 
 
\Step 2  Write $\D_1\subset \D(_b\cU_H)$ for the full triangulated subcategory generated by objects of $\D(_b\cU_H)$ which are restrictions from $\Bun_H$ of the almost constant local systems. 
  
  Let $\alpha$ be the highest weight of $W$. For $\cF\in\D(\Bun_H)$ the complex $_x\H^{\la}_H(W, \cF)\mid_{_b\cU_H}$ is completely determined by $\cF\mid_{_{b-1}\cU_H}$. Indeed, if $V,V'\in\Bun_H$ and $V\,\iso\, V'\mid_{X-x}$ such that $V$ is in the position $\alpha$ with respect to $V'$ then $V\in {_b\cU_H}$ implies $V'\in {_{b-1}\cU_H}$.
  
  Pick $N\ge 0$ such that for any perverse sheaf $\cA$ on $\Bun_H$ the complex $_x\H^{\la}_H(W, \cA)$ over $\Bun_H$ is placed in perverse degrees $[-N, N]$ (actually, one may take $N=\dim \Gr_H^{\alpha}$).
  
  Pick $a$ small enough compared to $b$ and satisfying the assumption of Lemma~\ref{Lm_Hecke_action_for_SL_2}. Then for $a'<a$ the cone of the natural map $_a\tilde K\to {_{a'}\tilde K}$ over $_{b-1}\cU_H$ is a succesive extension of constant complexes. 
  
  By Lemmas~\ref{Lm_Hecke_action_for_SL_2} and \ref{Lm_restricting_Aut_GH}, $F_H(_x\H^{\la}_{G}(W_1, {_aE}))$ over $_b\cU_H$ has a finite filtration in the derived category, one of the graded pieces is $\bar W_1\otimes (_{a+1}\tilde K)$, and the others are constant complexes.

  Write $_{\tau \ge ?}$ for the truncation functor with respect to the perverse t-structure. Apply $_{\tau\ge -N}$ for the isomorphism (\ref{iso_from_paper_L1}) over $_b\cU_H$. 
   
  From Proposition~\ref{Pp_appendix_A} we conclude that $_{\tau\ge -N}(_x\H^{\la}_H(W, {_a\tilde K}))$ over $_b\cU_H$ admits a finite filtration in the derived category, one of whose graded pieces is $_x\H^{\la}_H(W, \cK_H)$, and all the others are shifted almost constant local systems. Similarly, 
$$
_{\tau\ge -N}\left( F_H(_x\H^{\la}_{G}(W_1, {_aE}))\oplus (W_0\otimes {_a\tilde K})\right)
$$
over $_b\cU_H$ admits a finite filtration in the derived category, one of whose graded pieces is $(\bar W_1+W_0)\otimes \cK_H$ and all the others are shifted almost constant local systems. This implies already that $(\bar W_1+W_0)\otimes \cK_H$ appears as a direct summand in $_x\H^{\la}_H(W, \cK_H)$. More precisely, by decomposition theorem (\cite{BBD}), there is a complex $L\in \D(\Bun_H)$, which is a direct sum of shifted irreducible perverse sheaves, and an isomorphism in $\D(\Bun_H)$
$$
_x\H^{\la}_H(W, \cK_H)\,\iso\, (\bar W_1+W_0)\otimes \cK_H\oplus L
$$
We also see from the above that $_x\H^{\la}_H(W, \cK_H)\,\iso\, (\bar W_1+W_0)\otimes \cK_H$ is the quotient category of $\D(_b\cU_H)$ by $D_1$. So, $L\in D_1$. Since $D_1$ is closed under taking the direct summands, we conclude that each irreducible perverse sheaf appearing in $L$ lies in $D_1$, hence extends to $\Bun_H$ as an almost constant local system. Theorem~\ref{Th_2} is proved.

\bigskip

\centerline{\scshape 8. The perverse sheaf $\cK_H$ via Eisenstein series}

\bigskip\noindent
8.1 Recall the map $\nu_Q: \Bun_Q\to\Bun_H$ defined in Section~2.3.2. Write $\Bunb_Q$ for the stack classifying 
$V\in \Bun_H$ with an isotropic subsheaf $L\subset V$, where $L\in \Bun_1$. Let $\bar\nu_Q: \Bunb_Q\to\Bun_H$ be the projection sending this point to $V$. Write $\Bunb^m_Q\subset \Bunb_Q$ for the substack given by $\deg L=m$. The restriction $\bar\nu_Q^m: \Bunb^m_Q\to\Bun_H$ of $\bar\nu_Q$ is proper. Set 
$$
\cS^m=(\bar\nu_Q^m)_!\Qlb[\dim \Bun^m_Q]
$$
This complex differs from the usual definition of geometric Eisenstein series (\cite{BG}), as we used the constant sheaf instead of $\IC(\Bunb_Q)$ on the non smooth stack $\Bunb_Q$. 

 In this section we propose one more conjectural construction of the perverse sheaf $\cK_H$ as a `residue' of the sequence $\cS^m$ as $m$ goes to minus infinity. 
Set
$$
\bar\cS^m=\Four_{\cY_P, \psi}^{-1}(\nu_P^*\cS^m)\otimes(\Qlb[1](\frac{1}{2}))^{\dimrel(\nu_P)}\in \D(\cY_P)
$$

 Recall that $G_1$ introduced in Section~2.2 is the group scheme of automorphisms of $M_0=\cO_X\oplus\Omega$ acting trivially on $\det M_0$. Let $B_1\subset G_1$ be the Borel subgroup preserving $\cO_X$. Write $\Bun^m_{B_1}$ for the connected component of $\Bun_{B_1}$ classifying exact sequences 
\begin{equation}
\label{seq_L*_by_Omega_otimes_L}
0\to \Omega\otimes L\to M\to L^*\to 0
\end{equation}
with $L\in \Bun_1^m$. Let $\nu_{B_1}: \Bun_{B_1}^m\to \Bun_{G_1}$ be the map sending (\ref{seq_L*_by_Omega_otimes_L}) to $M$. Recall that $\cZ_{P,0}$ is the stack classifying $(U,M,s)$, where $U\in\Bun_n$, $M\in\Bun_{G_1}$ and $s: U\to M$ is a surjection. 

\begin{Lm} 
\label{Lm_support_barcS_m}
1) For each $m\in\ZZ$ the complex $\bar\cS^m$ is the extension by zero under the closed immersion $\cZ_P\hook{}\cY_P$.\\
2) The restriction of $\bar\cS^m$ to the open substack $\cZ_{P,0}^d\subset \cZ_P$ identifies canonically with 
$$
(\id\times\nu_{B_1})_!(\Qlb[1](\frac{1}{2}))^{-2(d+m)+n^2(g-1)}
$$
for the map $\id\times\nu_{B_1}: \cZ_{P,0}\times_{\Bun_{G_1}}\Bun_{B_1}^m\to \cZ_{P,0}$.
\end{Lm}
 
 Our proof of Lemma~\ref{Lm_support_barcS_m} uses a general construction presented separately in Section~8.2 for the convenience of the reader.
 
\medskip\noindent
8.2 {\scshape A stack associated to a complex} \  Consider a complex $\cM=(A\toup{d} B\to C)$ of locally free $\cO_X$-modules of finite ranks placed in cohomological degrees $0,1,2$. The maps in this complex are morphisms of coherent sheaves (not necessarily morphisms of vector bundles).

Let $\cX_{\cM}$ be the stack classifying an $A$-torsor $\cF_A$ on $X$, $s\in \H^0(X, B_{\cF_A})$ whose image in $\H^0(X, C)$ vanishes. Here $a\in A$ acts on $B$ sending $b\in B$ to $b+d(a)$, and $B_{\cF_A}$ is the quotient of $B\times \cF_A$ by $A$ acting diagonally. 

\begin{Lm} $\cX_{\cM}$ is naturally isomorphic to the stack quotient of $\H^1(X, \cM)$ by the trivial action of $\H^0(X,\cM)$.
\end{Lm} 
\begin{Prf}
Let $B'$ be the kernel of $B\to C$ and $\cM'=(A\toup{d} B')$ placed in degrees $0,1$. Then $\cX_{\cM}\,\iso\, \cX_{\cM'}$ naturally. Since $\H^i(X, \cM)\,\iso\, \H^i(X,\cM')$ for $i\le 1$, we may and do assume $C=0$. 

 The category of $A$-torsors on $X$ is equivalent to the category of exact sequences $0\to A\to E\to \cO_X\to 0$ on $X$, the datum of $s$ then becomes a datum of $\alpha: E\to B$ such that the composition $A\to E\toup{\alpha} B$ equals $d$. Thus, $\cX_{\cM}$ is the stack classifying diagrams on $X$
$$
\begin{array}{ccccccc}
0\to & B' & \toup{\id} & B' &\to & 0 & \to 0\\
& \uparrow\lefteqn{\scriptstyle d} && \uparrow\lefteqn{\scriptstyle\alpha} && \uparrow\\  
0\to & A & \to & E & \to & \cO_X & \to 0
\end{array}
$$
So, a point of $\cX_{\cM}$ gives rise to a distinguished triangle 
$\cM\to \cS\to \cO_X$ on $X$, where $\cS$ is the complex $(E\toup{\alpha} B')$ placed in degrees $0,1$. This triangle yields a morphism $\H^0(X,\cO_X)\to \H^1(X, \cM)$, hence a morphism of stacks $\gamma: \cX_{\cM}\to \H^1(X, \cM)$. The group $\H^0(X,\cM)$ acts on $\cX_{\cM}$ naturally by 2-automorphisms, so $\gamma$ extends to a morphism $\cX_{\cM}\to \H^1(X, \cM)/\H^0(X,\cM)$. One checks that this is an isomorphism.
\end{Prf}
  
\medskip  
  
  Example 1. Assume that $C=0$ and $d: A\to B$ is generically surjective. Then $\H^2(X, \cM)=0$, and $\cX_{\cM}$ is the stack classifying an exact sequence $0\to A\to ?\to \cO_X\to 0$ on $X$ together with a splitting of its push-forward via $d: A\to B$.
  
  Example 2. Let $U$ be a rank $n$ vector bundle on $X$ and $t: L\hook{} U^*$ be a subsheaf. Define the complex $\cM=(\wedge^2 U\toup{d_0} \HOM(L, U)\toup{d_1} \HOM(\Sym^2 L,\cO_X))$ as follows.
The map $d_0$ sends $y: U^*\to U$ such that $y^*=-y$ to the composition $L\toup{t} U^*\toup{y} U$. The map $d_1$ sends $z: L\to U$ to $\<z,t\>+\<t,z\>$. More precisely, here $\<z,t\>+\<t,z\>$ sends a local section $w_1w_2\in \Sym^2 L$ (with $w_i\in L$) to $\<z(w_1), t(w_2)\>+\<t(w_1), z(w_2)\>\in \cO_X$. The category of $\wedge^2 U$-torsors on $X$ is naturally equivalent to the category of exact sequences (\ref{seq_cO_by_wedge2U}) on $X$. Write $V$ for the image of (\ref{seq_cO_by_wedge2U}) under $\nu_P$, it is
included into an exact sequence (\ref{seq_U^*_by_U_gives_V}). Then $B_{\cF_A}$ is the sheaf of liftings $\tilde t: L\to V$ of the morphism $t: L\to U^*$. The condition that the image of $\tilde t$ in $\H^0(X, \Sym^2 L^*)$ vanishes means that the image of $\tilde t$ is isotropic. 

Thus, $\cX_{\cM}$ is the stack classifying an exact sequence
(\ref{seq_cO_by_wedge2U}) on $X$, and for the corresponding $V\in\Bun_H$ a commutative diagram
$$
\begin{array}{ccccc}
0\to U\to &V & \to &U^* &\to 0\\
&& \nwarrow\lefteqn{\scriptstyle \tilde t} & \uparrow{\lefteqn{\scriptstyle t}}\\
&&& L,
\end{array} 
$$
where the image of $\tilde t$ is isotropic. Write $U_1$ for the kernel of $t^*: U\to L^*$. If $L$ is of rank one then the kernel of $d_1$ equals $\HOM(L, U_1)$.
 
\medskip\noindent
8.3 \select{Proof of Lemma~\ref{Lm_support_barcS_m}} \  \  
1) The stack $\Bun_P\times_{\Bun_H}\Bunb_Q^m$ classifies an exact sequence (\ref{seq_cO_by_wedge2U}) on $X$
giving rise to an exact sequence 
\begin{equation}
\label{seq_U^*_by_U_gives_V}
0\to U\to V\to U^*\to 0
\end{equation}
on $X$ with $V\in\Bun_H$, and an isotropic subsheaf $L\subset V$ with $L\in\Bun_1^m$. Denote by $\cX_1\subset \Bun_P\times_{\Bun_H}\Bunb_Q^m$ the closed substack given by the condition that $L\subset U$, write $\cX_0$ for the complement of $\cX_1$ in $\Bun_P\times_{\Bun_H}\Bunb_Q$. Write $\cP_n$ for the stack classifying $U\in\Bun_n$ with a subsheaf $t: L\hook{} U^*$, where $L\in\Bun_1^m$. 

Clearly, the contribution of $\cX_1$ to the the complex $\bar\cS^m$ is the extension by zero under the zero section $\Bun_n\to \cY_P$. Consider the diagram
$$
\AA^1\getsup{\ev} \cY_P\times_{\Bun_n}\Bun_P \getsup{\id\times q_1} \cY_P\times_{\Bun_n}\cX_0\toup{q} \cY_P,
$$
where $\ev$ is the natural pairing between $v: \wedge^2 U\to\Omega$ and the exact sequence (\ref{seq_cO_by_wedge2U}), here $(U,v)\in\cY_P$. We have denoted by $q_1: \cX_0\to \Bun_P$ and $q$ the projections. We will show that 
$$
q_!((\id\times q_1)^*\ev^*\cL_{\psi})
$$ 
is the extension by zero from $\cZ_P$. Let $f_{\cP}: \cX_0\to \cP_n$ be the map sending a collection (\ref{seq_cO_by_wedge2U}) and $L\subset V$ to the composition $L\to V\to U^*$. Then $q$ is the composition
$$
\cY_P\times_{\Bun_n} \cX_0\toup{\id\times f_{\cP}} \cY_P\times_{\Bun_n} \cP_n\toup{\pr} \cY_P
$$ 
Consider a $k$-point $\eta$ of $\cP_n$ given by $t: L\hook{} U^*$. Write $U_1$ for the kernel of $t^*: U\to L^*$. As in (Section~8.2, example~2), we get a complex $\cM=(\wedge^2 U\toup{d} L^*\otimes U_1)$ placed in degrees $0,1$. The fibre $\cX_{\cM}$ of $f_{\cP}$ over $\eta$ identifies with the stack quotient of $\H^1(X,\cM)$ by $\H^0(X, \cM)$. 

 Since $d$ is generically surjective, $\H^2(X,\cM)=0$. The distinguished triangle $\cM\to \wedge^2 U\to L^*\otimes U_1$ on $X$
yields an exact sequence 
$$
\H^1(X,\cM)\to \H^1(X, \wedge^2 U)\to \H^1(X, L^*\otimes U_1)\to 0
$$ 
Thus, integrating $(\id\times q_1)^*\ev^*\cL_{\psi}$ over $\cX_{\cM}$, one gets zero unless $v\in \H^1(X, L^*\otimes U_1)^*$. So, the restriction of $v: U\to U^*\otimes\Omega$ to $U_1$ must factor through $L\otimes\Omega$, in particular $v: U\to U^*\otimes\Omega$ is of generic rank at most 2. So, 
\begin{equation}
\label{complex_for_Lm_about_cP_n}
(\id\times f_{\cP})_!(\id\times q_1)^*\ev^*\cL_{\psi}
\end{equation}
is the extension by zero under $\cZ_P\times_{\Bun_n}\cP_n\hook{} \cY_P\times_{\Bun_n}\cP_n$. 
Part 1) follows. 

\medskip\noindent
2) Let $^0\cP_n\subset \cP_n$ be the open substack given by the property that $v: L\hook{} U^*$ is a subbundle. Let us show that the restriction of (\ref{complex_for_Lm_about_cP_n}) to the open substack $\cZ_{P,0}\times_{\Bun_n}\cP_n$ is the extension by zero under $\cZ_{P,0}\times_{\Bun_n}{^0\cP_n}\hook{} \cZ_{P,0}\times_{\Bun_n}\cP_n$. Indeed, consider a $k$-point of $\cZ_{P,0}\times_{\Bun_n}\cP_n$ given by $s: U\to M$ and $t: L\hook{} U^*$. Assume that the $*$-fibre of (\ref{complex_for_Lm_about_cP_n}) at this point does not vanish. Let $U_1$ be the kernel of $t^*: U\to L^*$. We have seen in 1) that $v\in \H^0(X, \Omega\otimes U_1^*\otimes L)$. Let $D$ be an effective divisor on $X$ such that $t: L(D)\hook{} U^*$ is a subbundle. Then $v$ writes as a composition 
$$
\wedge^2 U\to U_1\otimes L^*(-D)\hook{} U_1\otimes L^*\to\Omega
$$
Since $v: \wedge^2 U\to\Omega$ is surjective, $D=0$. 

 Write $U_2$ for the kernel of $s: U\to M$. Since $v$ vanishes on $\wedge^2 U_1$, we also get $U_2\subset U_1$, and the exact sequence $0\to U_1/U_2\to M\to L^*\to 0$ is a point of $\Bun_{B_1}^m$. We have a closed immersion 
$$
i_0: \cZ_{P,0}\times_{\Bun_{G_1}}\Bun^m_{B_1}\hook{} \cZ_{P,0}\times_{\Bun_n}{^0\cP_n}
$$ 
given by the condition that $t: L\hook{} U^*$ factors through $s^*: M^*\to U^*$. We conclude that 
the $*$-restriction of (\ref{complex_for_Lm_about_cP_n})
to $\cZ_{P,0}\times_{\Bun_n}{^0\cP_n}$ identifies with 
$(i_0)_!\Qlb$ up to a shift and a twist.

 To calculate the shift note that $\Bun_Q^m$ is smooth of dimension 
$$
m(2-2n)+(2n^2-3n+2)(g-1)
$$ 
and $\dim\Bun_H=(2n^2-n)(g-1)$. Further, $\dim\Bun_P^d=(1-n)d+\frac{3n^2-n}{2}(g-1)$. For a point of $^0\cP_n$ as above, $\cM\,\iso\, \wedge^2 U_1$, so 
$$
\dim\cX_{\cM}=-\chi(\wedge^2 U_1)=(2-n)(d+m)+\frac{(n-1)(n-2)}{2}(g-1),
$$ 
where $d=\deg U$. Lemma~\ref{Lm_support_barcS_m} follows.

\medskip\noindent
8.4 Note that $\cZ^d_{P,0}$ is smooth of dimension $(n^2+3)(g-1)-2d$, and $\Bun_{B_1}^m$ is smooth of dimension $-2m$.   

\begin{Lm} 
\label{Lm_fibres_nu_B_1}
1) For $g=0$ (resp. for $g\ge 1$) assume that $m\le 1$ (resp., $m\le 2-2g$). Then 
\begin{equation}
\label{map_nu_B_1_m}
\nu_{B_1}^m:\Bun^m_{B_1}\to\Bun_{G_1}
\end{equation}
is generically smooth. If $g=0$ then the generic fibre of (\ref{map_nu_B_1_m}) is irreducible.\\
2) If $g\ge 1$ and $m< 3-3g$ then the generic fibre of (\ref{map_nu_B_1_m}) is irreducible.
\end{Lm}
\begin{Prf}
1) is elementary. 2) Recall the stack $_a\Bun_{G_1}$ introduced in Section~2.3.4. Under our assumption the stack $_{m+4g-4}\Bun_{G_1}$ is nonempty. Indeed, this follows from the semistability of generic $M\in \Bun_{G_1}$. Let $\nu: \cX\to {_{m+4g-4}\Bun_{G_1}}$ be the stack classifying a point $M\in {_{m+4g-4}\Bun_{G_1}}$, $L\in\Bun_1^m$ and a section $s: L\otimes\Omega\to M$. The projection $\cX\to {_{m+4g-4}\Bun_{G_1}}\times\Bun_1^m$ forgetting $s$ is a vector bundle of strictly positive rank. So, the generic fibre $\cX_{\tau}$ of the composition $\cX\to {_{m+4g-4}\Bun_{G_1}}\times\Bun_1^m\to {_{m+4g-4}\Bun_{G_1}}$ is irreducible. The generic fibre of (\ref{map_nu_B_1_m}) is open in $\cX_{\tau}$, so it is also irreducible.
\end{Prf}
 
\medskip

 Combining Lemmas~\ref{Lm_fibres_nu_B_1} and \ref{Lm_support_barcS_m} one gets the following.
\begin{Cor} 
\label{Cor_bar_cS_contains_IC(Z_P)}
Assume that $^e\cZ^d_{P,0}$ is not empty.
For $g=0$ assume  $m\le 1$, for $g\ge 1$ assume $m\le 2-2g$ (resp., $m<3-3g$). Then the perverse sheaf $^p\H^{3-3g-2m}(\bar\cS^m)$ over $^e\cY^d_P$ contains $\IC(\cZ_P)$ (resp., contains $\IC(\cZ_P)$ with multiplicity one). \QED
\end{Cor}

\begin{Rem} 
\label{Rem_Galois_covering_with_Galois_group_W}
i) The following is well-known (a similar claim with moduli stacks replaced by coarse moduli spaces is
proved in \cite{Las}). Assume $g=1$. Let  $G$ be a semisimple connected group, $T\subset G$ its maximal torus, $W$ the Weyl group of $(G,T)$. Then there is an open substack $\cW\subset \Bun_G^0$ over which the natural map $\nu_T^0: \Bun_T^0\to\Bun_G^0$ is a Galois covering with Galois group $W$. Here the action of $W$ on $\Bun_T^0$ is the one induced by the standard $W$-action on $T$. Given an irreducible representation $\sigma$ of $W$,  denote by 
$\cL_{\sigma}$ a perverse sheaf, which is the intermediate extension under $\cW\hook{}\Bun_G^0$ of the isotypic component of $(\nu_T^0)_!\Qlb\mid_{\cW}$ corresponding to $\sigma$. Since $\Bun_T^0$ is irreducible, each $\cL_{\sigma}$ is an irreducible perverse sheaf.

\medskip\noindent
ii) Using i) one can strengthen Corollary~\ref{Cor_bar_cS_contains_IC(Z_P)} in the case $g=1$ as follows. If $^e\cZ^d_P$ is not empty then the perverse sheaf $^p\H^0(\bar\cS^0)$ over $^e\cY^d_P$ contains $\IC(\cZ_P)$ with multiplicity one. Indeed, for $g=1$ and $m=0$ the map (\ref{map_nu_B_1_m})  over a suitable open substack of $\Bun_{G_1}$ is a Galois covering with Galois group $\ZZ/2\ZZ$.
\end{Rem} 

\medskip\noindent
8.5 From Corollary~\ref{Cor_bar_cS_contains_IC(Z_P)} and Lemma~\ref{Lm_general_ab_categories} one derives the following.

\begin{Cor} 
\label{Cor_sheaves_cS^m_d}
For $g=0$ assume $m\le 1$. For $g=1$ assume $m\le 0$. For $g>1$ assume $m<3-3g$. 
If $^e\cZ^d_{P,0}$ is not empty then $^p\H^{3-3g-2m}(\cS^m)$ contains a unique irreducible subquotient $\cS^m_d$ with the following property. The perverse sheaf
$$
\nu_P^*\cS^m_d\otimes(\Qlb[1](\frac{1}{2}))^{\dimrel(\nu_P)}
$$
over $^e\Bun_P^d$ contains $K^d_{P,\psi}$ as an irreducible subquotient. \QED
\end{Cor}

\begin{Rem} 
\label{Rem_construction_cK_H_via_Eis_for_Q}
We expect that each perverse sheaf $\cS^m_d$ from Corollary~\ref{Cor_sheaves_cS^m_d} is isomorphic to $\cK_H$ over $\Bun_H^{d\! \mod 2}$. Though we did not check this claim completely (except in the cases $g=0$ and $g=1$ considered in Sections~8.7 and 8.8), a partial evidence for this is collected in Section~8.6 for the convenience of the reader. 
\end{Rem}

\medskip\noindent
8.6 {\scshape Partial evidence for Remark~\ref{Rem_construction_cK_H_via_Eis_for_Q}} \   Write $F_H: \D^-(\Bun_{G_1})_!\to \D^{\prec}(\Bun_H)$ for the theta-lifting functor introduced in (\cite{L2}, Definition~2). 
For the map (\ref{map_nu_B_1_m}) set 
$$
\cF^m_{G_1}=(\nu_{B_1}^m)_!(\Qlb[1](\frac{1}{2}))^{
3-3g-4m}
$$ 

 Let $\cW_Q^m$ be the stack classifying $L\in \Bun_1^m$, $V\in\Bun_Q$ and an isotropic section $s: L\to V$. Denote by $\nu_{\cW}^m: \cW_Q^m\to \Bun_H$ the map sending the above collection to $V$. 
 
\begin{Lm} 
\label{Lm_calculating_F_H_cF_G1}
There is an isomorphism over $\Bun_H$
$$
F_H(\cF^m_{G_1})\,\iso\,  (\nu_{\cW}^m)_!(\Qlb[1](\frac{1}{2}))^r,
$$ 
where $r=-2nm+\dim\Bun_H+(2n+1)(1-g)=3-3g-2m+\dim\Bun_Q^m$.
\end{Lm}
\begin{Prf} For the map $\nu_{B_1}\times\id: \Bun_{B_1}\times\Bun_H\to \Bun_{G_1}\times\Bun_H$ the complex 
$(\nu_{B_1}\times\id)^*\Aut_{G_1, H}$ is as follows. Let $\cW_{B_1,H}$ be the stack classifying $V\in\Bun_H$, a point (\ref{seq_L*_by_Omega_otimes_L}) of $\Bun_{B_1}^m$, and any section $s: L\to V$. 

 For a point of $\cW_{B_1, H}$ write $\bar s$ for the composition $L^2\to \Sym^2 V\to\cO_X$. Let $\ev_{\cW}: \cW_{B_1,H}\to\A^1$ be the map sending 
the above collection to the pairing of (\ref{seq_L*_by_Omega_otimes_L}) with $\bar s$. Let $p_{\cW}: \cW_{B_1,H}\to \Bun_{B_1}\times \Bun_H$ be the projection forgetting $s$. By (\cite{L3}, Proposition~1), there is an isomorphism
$$
(\nu_{B_1}\times\id)^*\Aut_{G_1, H}\otimes(\Qlb[1](\frac{1}{2}))^{\dimrel(\nu_{B_1})}
\,\iso\, p_{\cW !}\ev_{\cW}^*\cL_{\psi}\otimes(\Qlb[1](\frac{1}{2}))^b,
$$
where $b$ is a function of a connected component of $\cW_{B_1, H}$ whose value at a point $((\ref{seq_L*_by_Omega_otimes_L}), L\toup{s} V)$ equals $\dim\Bun_{B_1}^m+\dim\Bun_H+\chi(L^*\otimes V)$, here $m=\deg(L)$. 

 Write $\cW_{B_1,H}^m\subset\cW_{B_1,H}$ for the substack given by the property $\deg L=m$. Let $\ov{\cW}_Q^m$ be the stack classifying $L\in\Bun^m_1, V\in\Bun_H$ and any section $s: L\to V$. By definition, $F_H(\cF^m_{G_1})$ is the direct image with compact support under the projection $\cW_{B_1,H}^m\to\Bun_H$. The latter decomposes as $\cW_{B_1,H}^m\toup{\xi} \ov{\cW}_Q^m\to \Bun_H$. The direct image $\xi_!\ev_{\cW}^*\cL_{\psi}$ is the extension by zero from the closed substack $\cW^m_Q\hook{}\ov{\cW}^m_Q$. Our assertion follows.
\end{Prf} 

\medskip

 Note that $\Bunb^m_Q\subset \cW^m_Q$ is the open substack given by the condition that $s: L\to V$ does not vanish. So, Lemma~\ref{Lm_calculating_F_H_cF_G1} yields a natural map over $\Bun_H$
\begin{equation}
\label{map_throwing_zero_section}
\cS^m\otimes(\Qlb[1](\frac{1}{2}))^{3-3g-2m}\to F_H(\cF^m_{G_1})
\end{equation}
whose cone is a constant complex. 

 Recall that the complex $F_H(\IC(\Bun_{G_1}))\,\iso\, q_{H !}\Aut_{G_1,H}$ does not literally make sense, here $q_H: \Bun_{G_1}\times\Bun_H\to\Bun_H$ is the projection. However, our perverse sheaf $\cK_H$ appears in $^p\H^0$ of a suitable truncation of the latter complex. 

 Assume that $m$ satisfies the conditions of Corollary~\ref{Cor_sheaves_cS^m_d} then $^p\H^0(\cF^m_{G_1})$ contains $\IC(\Bun_{G_1})$ with multiplicity one. So, for this $m$ the perverse sheaf $^p\H^0(F_H(\cF^m_{G_1}))$ should contain $\cK_H$. Now (\ref{map_throwing_zero_section}) shows that $\cK_H$ should appear in $^p\H^{3-3g-2m}(\cS^m)$. By Corollary~\ref{Cor_sheaves_cS^m_d}, $\cK_H$ can appear as an irreducible subquoient of $^p\H^{3-3g-2m}(\cS^m)$ with multiplicity at most one. 
 
\medskip\noindent
8.7 {\scshape Case $g=0$} 

\medskip\noindent
Our purpose is to prove Proposition~\ref{Pp_answer_for_g=0}. 
We will also calculate the sheaves $\cS^1_d$ and compare the answers (the two calculations are independent and will produce the same result). 

 We will use the Shatz stratification of $\Bun_H$ (cf. \cite{BD}, Section~2.10.4 and also \cite{BH}, \cite{Sh}). Let $T\subset B\subset H$ be a maximal torus and Borel subgroups. Let $\vartriangle$ be the corresponding set of simple roots of $B$. Write $\Lambda_H^+$ for the dominant coweights of $H$. For $g=0$ the Shatz strata are indexed by $\Lambda^+_H$. Namely, for $\lambda\in\Lambda^+_H$ let $M^{\lambda}\subset H$ be the standard Levi whose simple roots are $\check{\alpha}\in \vartriangle$ such that $\<\lambda, \check{\alpha}\>=0$. Let $P^{\lambda}$ be the standard parabolic subgroup with Levi factor $M^{\lambda}$. Write $\cF_{M^{\lambda}}$ for the push-forward of $\cO(1)$ under $\Gm\toup{\lambda} T\hook{} M^{\lambda}$. Let $Shatz^{\lambda}\subset \Bun_{P^{\lambda}}$ be the open substack classifying $\cF_{P^{\lambda}}$ such that $\cF_{P^{\lambda}}\times_{P^{\lambda}} M^{\lambda}$ is isomorphic to $\cF_{M^{\lambda}}$. 
The natural map $Shatz^{\lambda}\to\Bun_H$ is a locally closed immersion, and these substacks form the Shatz stratification. 
  
  For $b\in\ZZ/2\ZZ$ write $OSh^b$ for the open Shatz stratum in $\Bun_H^b$. Then $OSh^0=Shatz^{\lambda}$ for $\lambda=0$, and $OSh^1=Shatz^{\lambda}$ for $\lambda=(1,0,\ldots, 0)$. 
  
  Note that $\dim\Bun_Q^1=\dim\Bun_H=n-2n^2$. The stack $\Bunb_Q^1$ classifies $V\in\Bun_H$ with an isotropic subsheaf $L\subset V$ such that there is an isomorphism $L\,\iso\, \cO(1)$. The open stratum $OSh^0$ is not in the image of $\bar\nu_Q^1: \Bunb^1_Q\to\Bun_H$. The map $\bar\nu_Q^1$ is an isomorphism over $OSh^1$. So, for each $b\in\ZZ/2\ZZ$ the perverse sheaf $^p\H^1(\cS^1)$ vanishes over $OSh^b$. 
  
\begin{Lm} For each $b\in\ZZ/2\ZZ$ the stack $\Bun_H^b-OSh^b$ is irreducible, its open Shatz stratum is $Shatz^{\lambda}$, where $\lambda=(1,1,0,\ldots,0)$ (resp., $\lambda=(1,1,1,0,\ldots, 0)$) for $b=0$ (resp., $b=1$).
The perverse sheaf $^p\H^1(\cS^1)$ vanishes over $OSh^b$, and over the subregular Shatz stratum $Shatz^{\lambda}$ there is an isomorphism 
\begin{equation}
\label{iso_over_subregular_stratum}
^p\H^1(\cS^1)\,\iso\, \IC(Shatz^{\lambda})
\end{equation}
\end{Lm}
\begin{Prf}  
1) The image of the proper map $\bar\nu_Q^1$ in $\Bun_H^0$ equals $\Bun_H^0-OSh^0$. By (\cite{BG}, Proposition~1.3.8), $\Bun_Q^1$ is dense in $\Bunb^1_Q$, so $\Bunb^1_Q$ is irreducible. This implies that $\Bun_H^0-OSh^0$ is irreducible. The open Shatz stratum in $\Bun_H^0-OSh^0$ is $Shatz^{\lambda}$ for $\lambda=(1,1,0,\ldots, 0)$, the subregular Shatz stratum. 
For $V\in Shatz^{\lambda}$ there is an isomorphism 
$$
V\,\iso\, \cO(1)\oplus\cO(1)\oplus \cO^{2n-4}\oplus \cO(-1)\oplus\cO(-1)
$$
So, the fibre of $\bar\nu^1_Q$ over $V$ identifies with $\PP^1$. The codimension of $Shatz^{\lambda}$ in $\Bun_H^0$ is one, so we get an isomorphism over $Shatz^{\lambda}$
$$
\cS^1\,\iso\, \IC(Shatz^{\lambda})[1]\oplus \IC(Shatz^{\lambda})[-1],
$$
and the desired isomorphism (\ref{iso_over_subregular_stratum}) over 
$Shatz^{\lambda}$.
  
\noindent
2) Recall the parabolic subgroup $R\subset H$ defined in  Section~2.3.3. Note that $R/[R,R]\,\iso\,\Gm$. The Levi quotient of $R$ identifies with $\GL_2\times H_{n-2}$. Write $\check{\Lambda}_{H,R}^+$ for the semigroup of $H$-dominant weights which are orthogonal to all the simple coroots of $\GL_2\times H_{n-2}$.

Let $\Bunb_R$ be the stack classifying a $R/[R,R]$-torsor $\cF_{R/[R,R]}$ on $X$, an $H$-torsor $\cF_H$ on $X$, and for each $\check{\lambda}\in \check{\Lambda}_{H,R}^+$ a map $\kappa^{\check{\lambda}}: \cL^{\check{\lambda}}_{\cF_{R/[R,R]}}\hook{} \cV^{\check{\lambda}}_{\cF_H}$ such that the Pl\" ucker relations hold as in (\cite{BG}, Section~1.3.2). Here $\cV^{\check{\lambda}}$ is the corresponding Weyl module (as in \cite{BG}, Section~0.4.1). We may simply think of $\Bunb_R$ as the stack classifying $L\in\Bun_1$, $V\in\Bun_H$ and a section $\kappa: L\hook{} \wedge^2 V$ such that the Pl\" ucker relations hold. Write $\Bunb_R^2$ for the substack of $\Bunb_R$ given by the properties $\deg L=2$ and $V\in\Bun_H^1$. The projection $\bar\nu_R: \Bunb_R^2\to\Bun_H^1$ is proper, and its image equals $\Bun_H^1- OSh^1$. 
 
  Let $\Bun_R^2\subset \Bunb^2_R$ be the open substack given by the property that $L\hook{}\wedge^2 V$ is a subbundle. As in (\cite{BG}, Proposition~1.3.8) one checks that $\Bun_R^2$ is dense in $\Bunb_R^2$. Since $\Bun_R^2$ is an irreducible component of $\Bun_R$, $\Bunb^2_R$ is irreducible, so $\Bun_H^1- OSh^1$ is also irreducible.
  
  The open Shatz stratum in $\Bun_H^1- OSh^1$ is $Shatz^{\lambda}$ for $\lambda=(1,1,1,0,\ldots,0)$. For any $V\in Shatz^{\lambda}$ the fibre of $\bar\nu_Q^1$ over $V$ identifies with $\PP^2$. The codimension of $Shatz^{\lambda}$ in $\Bun_H^1$ is 3. So, the $*$-restriction of $\cS^1$ to $Shatz^{\lambda}$ identifies with
$$
\IC(Shatz^{\lambda})[3]\oplus  \IC(Shatz^{\lambda})[1]\oplus \IC(Shatz^{\lambda})[-1] 
$$
The $*$-retsriction of $\IC(\Bun_H)$ to $Shatz^{\lambda}$  
identifies with $\IC(Shatz^{\lambda})[3]$. This yields the desired isomorphism (\ref{iso_over_subregular_stratum}). 
\end{Prf}  
  
\medskip

 One has the involution $s$ of $\Lambda^+_H$ sending 
$\lambda=(a_1,\ldots, a_n)$ to $s\lambda=(a_1,\ldots, a_{n-1}, -a_n)$. Note that $\dim Shatz^{\lambda}=\dim Shatz^{s\lambda}$, and the fibre of $\bar\nu^1_Q$ over a point of $Shatz^{\lambda}$ identifies with $\PP^{a-1}$, where $a=a_1+\ldots+a_{n-1}+\mid\! a_n\!\mid$. Indeed, for $V\in Shatz^{\lambda}$ one has $\dim\Hom(\cO(1),V)=a$, and any section $\cO(1)\to V$ is isotropic. Let $G_V$ be the group scheme of automorphisms of $V$ preserving the symmetric form. 

  Assume that $\lambda=(a_1,\ldots, a_m, 0,\ldots, 0)$ with $a_m>0$ and 
$$
\lambda=(b_1,\ldots, b_1; b_2,\ldots, b_2;\ldots ;b_k,\ldots, b_k; 0,\ldots, 0),
$$
where $b_i$ appears $r_i$ times for $i=1,\ldots,k$, and $b_1>\ldots >b_k>0$. 

\begin{Lm} 
\label{Lm_dim_G_V}
Set $a=\sum_i a_i$. Then
one has 
\begin{multline*}
\dim G_V=(n-m)(2n-2m-1)+\sum_{1\le i\le j\le k} r_ir_j(1+b_i-b_j)\\ +(2n-2m)(m+a)+\frac{m(m-1)}{2}+(m-1)a
\end{multline*}
\end{Lm}
\begin{Prf} We have 
\begin{equation}
\label{V_decomposition_g=0}
V\,\iso\, W\oplus \cO^{2n-2m}\oplus W^*
\end{equation}
with $W\,\iso\, \cO(a_1)\oplus\ldots\oplus \cO(a_m)$. 
Recall that $\dim H_n=n(2n-1)$. One gets for the Levi part 
$$
\dim H_{n-m}+\sum_{1\le i\le j\le k} r_ir_j(1+b_i-b_j)
$$ 
We have to add for the unipotent part $\dim\Hom(V', W)+\dim\H^0(X, \wedge^2 W)$. One has
$$
\dim\Hom(\cO^{2n-2m}, W)=(2n-2m)(r_1(b_1+1)+\ldots+r_k(b_k+1))
$$ 
One also has 
$$
\dim\H^0(X, \wedge^2 W)=\sum_{1\le i<j\le m} (a_i+a_j+1)=\frac{m(m-1)}{2}+ (m-1)(\sum_i a_i)
$$
Note that $r_1b_1+\ldots+r_kb_k=a_1+\ldots+a_m$ and $\sum r_i=m$. 
\end{Prf}

\begin{Lm} 
\label{Lm_perverse_deg_less_zero_for_cS1}
1) The $*$-restriction of $\cS^1$ to any Shatz stratum (except the open ones and the subregular ones) is placed in perverse degrees $\le 0$.\\
2) Let $b\in\ZZ/2\ZZ$. For all the Shatz strata in $\Bun_H^b-OSh^b$ one has 
\begin{equation}
\label{key_inequality}
2a-3\le\codim Shatz^{\lambda}=\dim\Bun_H-\dim Shatz^{\lambda},
\end{equation}
The inequality is strict unless $Shatz^{\lambda}$ is the subregular Shatz stratum. Here for $\lambda=(a_1,\ldots, a_n)$ we set $a=\sum a_i$. 
\end{Lm}
\begin{Prf} 
1) Follows immediately from 2).\\
2) Use the notations of Lemma~\ref{Lm_dim_G_V}. 
For $V\in Shatz^{\lambda}$ given by (\ref{V_decomposition_g=0}) we have
$\dim\End(W)\ge m^2$. Indeed, if $i\ne j$ then $\dim\Hom(\cO(a_i), \cO(a_j))\oplus \Hom(\cO(a_j),\cO(a_i))\ge 2$.

 By Lemma~\ref{Lm_dim_G_V}, it suffices to show that
\begin{equation}
\label{eq_first_to_prove}
-4nm+m(2m+1)+m^2+(2n-2m)(m+a)+\frac{m(m-1)}{2}+(m-1)a\ge 2a-3,
\end{equation} 
where $a=\sum a_i$, and the equality holds only in the cases indicated above. Now (\ref{eq_first_to_prove}) rewrites as
\begin{equation}
\label{eq_second_to_prove}
2a(2n-3-m)\ge -3m^2+m(4n-1)-6
\end{equation}
We always have $a\ge m$ and the equality is strict unless $\lambda=(1,\ldots,1)$. Using the inequality $a\ge m$, we are reduced to show that
\begin{equation}
\label{eq_third_to_prove}
2m(2n-3-m)\ge -3m^2+m(4n-1)-6
\end{equation}
The latter inequality rewrites as $m^2-5m+6\ge 0$. The quadratic function $x^2-5x+6$ takes its minimal value $-1/4$ at $x=5/2$. So, if $m\in\ZZ$ then $m^2-5m+6\ge 0$ and the equality takes place exactly for $m=2$ and $m=3$. 

 The cases $m=2$ and $m=3$ under the condition $a=m$ 
correspond exactly to the subregular coweights $\lambda$. For them (\ref{key_inequality}) is an equality, otherwise the inequality (\ref{key_inequality}) is strict. 
\end{Prf} 

\medskip

\begin{Lm} 
\label{Lm_answer_for_cS^1_d_for_g=0}
Let $d$ be as in Corollary~\ref{Cor_sheaves_cS^m_d}. Then $\cS^1_d$ is isomorphic to the $\IC$-sheaf of the subregular Shatz stratum over $\Bun_H^{d\mod 2}$.
\end{Lm}
\begin{Prf}
\Step 1 Recall that for each $a\in\ZZ$ one has the open substack $_a\cU_H\subset \Bun_H$ defined in Section~7.1. We claim that the isomorphism (\ref{iso_over_subregular_stratum}) actually holds over $_{-2}\cU_H$. Indeed, the preimage of $_{-2}\cU_H$ under $\bar\nu_Q^1: \Bunb_Q^1\to\Bun_H$ is contained in $\Bun_Q^1$, which is smooth. So, $\cS^1$ is self-dual over $_{-2}\cU_H$. Lemma~\ref{Lm_perverse_deg_less_zero_for_cS1} now implies that for any $b\in\ZZ/2\ZZ$, the perverse sheaf $^p\H^1(\cS^1)$ over $_{-2}\cU_H\cap\Bun_H^b$ is the intermediate extension from the subregular Shatz stratum. 

\smallskip

\Step 2  Assume first that $d=-n$ or $1-n$.
Let $\cW_n^{-n}$ (resp., $\cW_n^{1-n}$) be the stack classifying vector bundles $U\in\Bun_n$ isomorphic to $\cO(-1)^n$ (resp., to $\cO(-1)^{n-1}\oplus\cO$). Then $\cW_n^d\subset {^e\Bun_n^d}$ is a substack. Since $n\ge 4$, any $U\in \cW_n^d$ admits a quotient vector bundle isomorphic to $\cO(-1)\oplus\cO(-1)$, so the preimage of $\cW_n^d$ in $^e\cZ_{P,0}^d$ is not empty. Write $\cW\Bun_P^d$ for the preimage of $\cW_n^d$ under $\Bun_P^d\to\Bun_n^d$. Then $K_{P,\psi}^d$ does not vanish over $\cW\Bun_P^d$. Since the image of $\cW\Bun_P^d$ in $\Bun_H^d$ is contained in $_{-2}\cU_H^{d\!\mod 2}$, $K_{P,\psi}^d$ does not vanish over $^e\Bun_P^d\cap \nu_P^{-1}(_{-2}\cU_H^{d\!\mod 2})$. Now, for example by Theorem~1, the same holds for any $d$ such that $^e\cZ^d_{P,0}$ is not empty.
Our assertion follows from Corollary~\ref{Cor_sheaves_cS^m_d} and Step 1.
\end{Prf}
   
\medskip

\begin{Prf}\select{of Proposition~\ref{Pp_answer_for_g=0}} 

\medskip\noindent
For $b\in\ZZ/2\ZZ$ let $O\!\Bun_H^b\subset \Bun_H^b$ be the open substack equal to the union of $OSh^b$ and the subregular Shatz stratum in $\Bun_H^b$. 
We have already seen in Lemma~\ref{Lm_answer_for_cS^1_d_for_g=0} that for each $b\in\ZZ/2\ZZ$, $\cK_H$ does not vanish on $_{-2}\cU_H^b$. Now by Corollary~\ref{Cor_modulo_constant_complexes} and Proposition~\ref{Pp_K^d_coincide}, for each $b\in\ZZ/2\ZZ$ the perverse sheaf $^p\H^0(_{-2}\tilde K)\mid_{_{-2}\cU_H^b}$ contains the unique irreducible subquotient $\cK_H\mid_{_{-2}\cU_H^b}$. 

 Now it suffices to show that for each $b\in\ZZ/2\ZZ$ the perverse sheaf $^p\H^0(_{-2}\tilde K)$ over $O\!\Bun_H^b$ identifies with $\IC(Shatz^{\lambda})$, where $Shatz^{\lambda}$ is the subregular stratum.
 
  To do so, note that $_{-2}\Bun_{G_1}\subset\Bun_{G_1}$ is the open substack classifying trivial $G_1$-torsors, it is isomorphic to the classifying stack $B(G_1)$. Recall the map $_{-2}q: {_{-2}\Bun_{G_1}}\times\Bun_H\to\Bun_H$ from Section~2.3.4. 
  
  Let first $b=0$ and $\lambda=(1,1,0,\ldots,0)$. Over ${_{-2}\Bun_{G_1}}\times OSh^0$ the complex $\Aut_{G_1,H}$ identifies with $\Qlb[\dim\Bun_H-3]$, and the $*$-restriction of $\Aut_{G_1,H}$ to ${_{-2}\Bun_{G_1}}\times Shatz^{\lambda}$ identifies with $\Qlb[\dim\Bun_H-7]$ by (\cite{L1}, Theorem~1). Since the codimension of $Shatz^{\lambda}$ in $\Bun_H^0$ is one and 
$$
\RG_c(B(G_1),\Qlb)\,\iso\, \oplus_{k=0}^{\infty}\Qlb[6+4k]\;\iso\; \Qlb[6]\oplus \Qlb[10]\oplus\ldots,
$$  
our assertion for $b=0$ follows.

 Let now $b=1$ and $\lambda=(1,1,1,0,\ldots, 0)$. Over ${_{-2}\Bun_{G_1}}\times OSh^1$ the complex $\Aut_{G_1,H}$ identifies with $\Qlb[\dim\Bun_H-5]$. The $*$-restriction of $\Aut_{G_1,H}$ to ${_{-2}\Bun_{G_1}}\times Shatz^{\lambda}$ identifies with $\Qlb[\dim\Bun_H-9]$. Since the codimension of $Shatz^{\lambda}$ in $\Bun_H^1$ is 3, our assertion similarly follows for $b=1$. 
\end{Prf}   
   
\medskip

\begin{Rem} Actually, one may show that for $g=0$ and any $d\in Z(e,P)$ the fibre of $\nu_P: {^e\Bun_P^d}\to\Bun_H$ over a point of the subregular Shatz stratum of $\Bun_H^{d\mod 2}$ is irreducible. So, in this case the isomorphism (\ref{iso_Th1_on_eBunP}) of Theorem~\ref{Th1} determines $\cK_H$ up to a unique isomorphism.
\end{Rem} 
   
\medskip\noindent
8.8 {\scshape Case $g=1$}

\medskip\noindent
8.8.1 Let $T\subset H$ be the standard maximal torus, write $W$ for the Weyl group of $(H,T)$. Sometimes we write $W=W(H_n)$ to express the dependence on $n$. 
The stack $\Bun_T^0$ classifies $U_1,\ldots, U_n\in\Bun_1^0$, we have denoted by $U_i$ the push-forward of the $T$-torsor under the weight $(0,\ldots, 0,1,0,\ldots, 0)$, where $1$ appears on $i$-th place.
 
 The natural map $\nu_T^0: \Bun_T^0\to\Bun_H^0$ sends this point to $V=\sum_i(U_i\oplus U_i^*)$ with the induced symmetric form $\Sym^2 V\to\cO_X$ and a trivialization $\det V\,\iso\, \cO_X$. 
 
 Write $^0\Bun_T^0\subset \Bun_T^0$ for the open substack given by the properties: $U_i\otimes U_j$ is nontrivial for all $i,j$, and $U_i$ is not isomorphic to $U_j$ for $i\ne j$. Let $\cW^0_H\subset \Bun_H^0$ be the image of  $^0\Bun_T^0$ under $\nu_T^0$. The restriction $^0\Bun_T^0\to \cW^0_H$ of $\nu_T^0$ is a Galois covering with Galois group $W$ (cf. Remark~\ref{Rem_Galois_covering_with_Galois_group_W}). 
   
  Recall the stack $\cW^1_{H_2}$ introduced in Section~2.4.2. For $n\ge 3$ consider the map 
\begin{equation}
\label{etale_map_f^1}  
f^1: \Bun^1_{H_2}\times \Bun^0_{H_{n-2}}\to \Bun_H^1
\end{equation}
sending $(V, V')$ to $V\oplus V'$, the symmetric form being the orthogonal sum of the forms for $V, V'$. The restriction of $f^1$ to the open substack $\cW^1_{H_2}\times \cW^0_{H_{n-2}}$ is \'etale, and the image of this open substack under $f^1$ will be denoted $\cW^1_H$. Write $\cW_H$ for the disjoint union of $\cW^0_H$ and $\cW^1_H$.
   
  Let $W'\subset W$ be the stabilizor of the coweight $(1,0,\dots, 0)$ in $W$. The induced representation $ind_{W'}^W(\Qlb)$ from the trivial representaition of $W'$ to $W$ decomposes as 
$$
ind_{W'}^W(\Qlb)\,\iso\,\Qlb\oplus \sigma_n\oplus \sigma'_n,
$$
where $\sigma_n, \sigma'_n$ are irreducible, $\sigma_n$ is the representation defined in Section~2.4.2, and $\Qlb$ is the trivial representation. Note that $\dim\sigma_n=n-1$ and $\dim\sigma'_n=n$. 
  
   According to Corollary~\ref{Cor_sheaves_cS^m_d}, we will look for $\cK_H$ inside the perverse sheaf $^p\H^0(\cS^0)$. The map $\bar\nu_Q^0: \Bunb^0_Q\to \Bun_H$ over $\cW_H$ is an \'etale covering. It is easy to see that 
\begin{equation}
\label{iso_cS^0_over_cW^0_case_g=1}
\cS^0\mid_{\cW^0_H}\,\iso\, \Qlb\oplus \cL_{\sigma_n}\oplus \cL_{\sigma'_n}
\end{equation}
naturally over $\cW^0_H$. Similarly, over $\cW^1_{H_2}\times \cW^0_{H_{n-2}}$ one has a natural isomorphism
\begin{equation}
\label{iso_after_f^1}
(f^1)^*\cS^0\,\iso\, \pr_2^*(\Qlb\oplus  \cL_{\sigma_{n-2}}\oplus \cL_{\sigma'_{n-2}}),
\end{equation}
where $\pr_2: \cW^1_{H_2}\times \cW^0_{H_{n-2}}\to \cW^0_{H_{n-2}}$ is the projection.
   
\medskip\noindent   
8.8.2  Recall for $a\in\ZZ$ the open substack $_a\cU_H$ introduced in Section~7.1. Note that $\cW_H\subset {_{-1}\cU_H}$. Recall the complex $_a\tilde K$ given by (\ref{complex_a_tildeK_on_BunH}). We will analyse $_{-1}\tilde K$ over $\cW_H$. This will be sufficient, because, by Corollary~\ref{Cor_modulo_constant_complexes}, for any $a\le -1$ the cone of the natural map $_a\tilde K\to {_{a-1}\tilde K}$ over $_{-1}\cU_H$ is a constant complex.  So, for any $a\le -1$ any non constant irreducible subquotient of $^p\H^0(_a\tilde K)\mid_{\cW_H}$ already appears in 
$^p\H^0(_{-1}\tilde K)\mid_{\cW_H}$.    

 Recall the stack $_a\Bun_{G_1}$ defined in Section~2.3.4. Note that $_{-1}\Bun_{G_1}\subset\Bun_{G_1}$ coincides with the open substack of semistable $G_1$-torsors. 

 Note that $\GL_n\subset H$ is the standard Levi subgroup containing $T$. Let $\nu_{T,n}: \Bun_T\to \Bun_n$ be the extension of scalars map with respect to $T\hook{} \GL_n$. It sends $(U_1,\ldots, U_n)$ to $U_1\oplus\ldots\oplus U_n$. Let $\cW\Bun_n^0\subset \Bun_n^0$ be the image of $^0\Bun_T^0$ under $\nu_{T,n}$. The restriction 
$$
\nu_{T,n}: {^0\Bun_T^0}\to \cW\Bun_n^0
$$ 
is a Galois covering with Galois group $S_n$, the Weyl group of $(T, \GL_n)$. For an irreducible representation $\tau$ of $S_n$ write $\cL_{\GL_n, \tau}$ for the isotypic component of $(\nu_{T,n})_!\Qlb\mid_{\cW\Bun_n^0}$ corresponding to $\tau$, this is an irreducible perverse sheaf on $\cW\Bun_n^0$. 

 Let $\cW\Bun_P^0$ be the preimage of $\cW\Bun_n^0$ under $\Bun_P^0\to\Bun_n^0$. The natural map $\cW\Bun_P^0\to \cW\Bun_n^0$ is an isomorphism, so we identify these two stacks. Clearly, $\cW\Bun_n^0\subset {^e\Bun_n^0}$. 
 
\begin{Lm} 
\label{Lm_about_cWBun_P}
There is a morphism in $\D^-(\cW\Bun_P^0)$
$$
\Qlb\oplus \cL_{\GL_n, \sigma_n}\to \nu_P^*(_{-1}\tilde K)
$$
whose cone is a constant complex.
\end{Lm}
\begin{Prf}
Recall the map $\pi_P: \cS_P\to\cY_P$ introduced in Section~2.3.1. Write $\cW\cS_P\subset \cS_P$ for the open substack classifying $(s: U\to M)\in \cS_P$ such that $M\in {_{-1}\Bun_{G_1}}$, $U\in \cW\Bun_n^0$. Let $^1\cW\cS_P\subset \cW\cS_P$ be the closed substack given by $s=0$, and $^0\cW\cS_P$ be its complement in $\cW\cS_P$. 
 
 If $U\in \cW\Bun_n^0$ then $\Hom(\wedge^2 U,\Omega)=0$, so the projection $\cY_P\to\Bun_n$ becomes an isomorphism over $\cW\Bun_n^0$. The restriction of $\pi_P$ to $\cW\cS_P$ becomes a morphism 
$\pi_{\cW, P}: \cW\cS_P\to \cW\Bun_n^0$ sending $(U,M,s)$ to $U$. This is our definition of $\pi_{\cW,P}$.  

 Proposition~\ref{Pp_explicit_from_Wald_periods} implies an isomorphism 
$$
\nu_P^*(_{-1}\tilde K)\,\iso\, (\pi_{\cW,P})_!\Qlb
$$ 
over $\cW\Bun_n^0$. The contribution of $^1\cW\cS_P$ to the latter direct image is a constant complex.
Finally, $^0\cW\cS_P$ identifies with the stack classifying $U\in \cW\Bun_n^0$ and a surjection $U\to U_1$ on $X$, where $U_1\in\Bun_1^0$. The map $\nu_{T,n}$ decomposes as
$$
^0\Bun_T^0\to {^0\cW\cS_P}\toup{^0\pi_{\cW,P}} \cW\Bun_n^0,
$$
where $^0\pi_{\cW,P}$ is the restriction of $\pi_{\cW,P}$.
This yields an isomorphism $(^0\pi_{\cW,P})_!\Qlb\,\iso\, \Qlb\oplus \cL_{\GL_n, \sigma_n}$. We are done.
\end{Prf} 
 
\medskip

 Lemma~\ref{Lm_about_cWBun_P} combined with Proposition~\ref{Pp_K^d_coincide} implies that $\cK_H$ does not vanish at the generic point of $\Bun_H^0$. 
Note that $\nu_P^*\cL_{\sigma_n}\,\iso\, \cL_{\GL_n, \sigma_n}$ canonically. Now (\ref{iso_cS^0_over_cW^0_case_g=1}) together with Corollary~\ref{Cor_sheaves_cS^m_d} show that all the perverse sheaves $\cS^0_d$ coincide with the subquotient  $\cL_{\sigma_n}$ of $^p\H^0(\cS^0)$. 
Now Proposition~\ref{Pp_fibres_of_nu_P_connected} together with Corollary~\ref{Cor_sheaves_cS^m_d} imply that $\cK_H\,\iso\,\cL_{\sigma_n}$ over $\Bun_H^0$. The first part of Proposition~\ref{Pp_case_g=1} is proved. 
 
\medskip\noindent 
8.8.3 {\scshape End of the proof of Proposition~\ref{Pp_case_g=1}} \  
Let $^{in}\Bun_n^1\subset \Bun_n^1$ be the open substack of indecomposable vector bundles. Recall that the map $^{in}\Bun_n^1\to\Bun_1^1$ sending $U$ to $\det U$ is an isomorphism (cf. \cite{Po}). Denote by $L\mapsto W_n(L)$ the inverse of this map.
Recall that if $U\in{^{in}\Bun_n^1}$ then $U$ is stable, in particular $\End(U)=k$.

 Let $^{in}\Bun_n^{-1}\subset\Bun_n^{-1}$ be the open substack of indecomposable vector bundles. The map $U\mapsto U^*$ yields an isomorphism $^{in}\Bun_n^1\,\iso\; {^{in}\Bun_n^{-1}}$.  
Let $^{in}\Bun_P^{-1}$ be the preimage of $^{in}\Bun_n^{-1}$ under $\nu_P: \Bun_P^{-1}\to\Bun_n^{-1}$. To finish the proof of Proposition~\ref{Pp_case_g=1}, we will analyze the perverse sheaf $K_{P,\psi}^{-1}$ over $^{in}\Bun_P^{-1}$. 

  First, let us remind some well-known properties of indecomposable vector bundles. They are either proved in or easily obtained from the results of \cite{Po}. Let $L\in\Bun_1^1$ and $W_n=W_n(L)$. One has $\H^0(X, W_n^*)=0$. If $\cA\in \Bun_1^1$ then $\dim\Hom(\cA, W_n)=1$, the image of a nonzero map $\cA\to W_n$ is a subbundle, and $W_n/\cA$ is indecomposable. By induction, $W_n$ admits a canonical flag of subbundles 
$$
0=W^0\subset W^1\subset\ldots\subset W^{n-1}\subset W_n
$$
such that $W^i/W^{i-1}$ is non canonically isomorphic to $\cO_X$, and $W_n/W^{n-1}$ is non canonically isomorphic to $L$. There is an exact sequence $0\to \cO_X\to W_n(L)\to W_{n-1}(L)\to 0$ on $X$. It gives rise to an exact sequence 
\begin{equation}
\label{seq_for_wedge2_W_n}
0\to W_{n-1}(L)\to \wedge^2 W_n(L)\to \wedge^2 W_n(L)\to 0
\end{equation}
This sequence allows to show by induction that $\H^0(X, \wedge^2 W_n^*)=0$ and $\dim\H^0(X, \wedge^2 W_n)=n-1$. This implies $^{in}\Bun_n^{-1}\subset {^e\Bun_n^{-1}}$.

 Further, any subsheaf of $W_n$ of degree $\ge 1$ coincides with $W_n$. If $\cA\in\Bun_1^0$ then $W_n(L)\otimes\cA$ is also indecomposable, so $W_n(L)\otimes\cA\;\iso\; W_n(L\otimes \cA^n)$. 
 
 Write $Y_2(L)$ for the scheme classifying subbundles $E\subset W_n(L)$ of rank 2 such that there exists an isomorphism $\det E\,\iso\, \cO_X$, but it is not fixed. In this definition one may equally require that $E$ is a subsheaf, then it is actually a subbundle. 
 
 Let $Y_1(L)$ be the scheme classifying subbundles of rank one and degree zero in $W_n(L)$. If $n\ge 2$ then $Y_1(L)\,\iso\, \uBun_1^0$ naturally, here $\uBun_1^0$ denotes the Picard scheme of line bundles of degree zero on $X$. The latter map sends $(\cA\subset W_n)$ to the isomorphism class of $\cA$.
  
  Let $Y_{1,1}(L)$  be the scheme classifying flags $E_1\subset E\subset W_n(L)$ in $W_n(L)$, where $E\in Y_2(L)$, and $E_1$ is a subbundle of rank one and degree zero. 
  
  Assume $n\ge 3$. The map $Y_{1,1}(L)\to Y_1(L)\,\iso\, \uBun_1^0$ sending $(E_1\subset E)$ to $E_1$ is an isomorphism. This follows from $\dim\Hom(E_1^*, E/E_1)=1$.
  
  The map $\pi_{1,1}: Y_{1,1}(L)\to Y_2(L)$ sending $(E_1\subset E)$ to $E$ is surjective and proper. If $E\subset W_n$ is a subbundle of rank 2 and determinant $\cO_X$, then $E$ is semistable. Pick a line subbundle $E_1\subset E$ with $\deg E_1=0$. If $E_1^2\,\iso\, \cO_X$ then the exact sequence $0\to E_1\to E\to E_1^*\to 0$ does not split, and $\pi_{1,1}^{-1}(E)$ is a point.
If $E_1^2$ is not trivial then the latter exact sequence splits, and $\pi_{1,1}^{-1}(E)$ consists of 2 points. This shows that $Y_2(L)$ is the quotient of $\uBun_1^0$ by the automorphism $\cA\mapsto \cA^{-1}$. 
In turn, this yields an isomorphism $Y_2(L)\,\iso\, \PP^1$. 
  
  Write $\cL_{Y_2}$ for the line bundle on $Y_2(L)$ with fibre $\H^0(X, \wedge^2 E)$ at $E$. Let $a_L\in\ZZ$ be such that $\cL_{Y_2}\,\iso\, \cO(a_L)$ as a line bundle on $\PP^1$. 
  
  Write $\tilde Y_2(L)$ for the scheme of $v\in \H^0(X, \wedge^2 W_n(L))$ such that the image of $v: W_n^*\to W_n$ is of generic rank at most 2. Then $\tilde Y_2(L) -\{0\}$ is the total space of $\cL_{Y_2}$ with zero section removed. 
  
\begin{Rem} Let $W$ be a finite-dimensional $k$-vector space. The exteriour product $(\wedge^2 W)\otimes(\wedge^2 W)\to \wedge^4 W$ is symmetric, so yields a map $\Sym^2(\wedge^2 W)\to\wedge^4 W$. An element $\omega\in\wedge^2 W$ is decomposable iff $\omega\wedge\omega=0$ in $\wedge^4 W$. So, $\tilde Y_2(L)$ is the scheme of sections $v\in\H^0(X, \wedge^2 W_n)$ such that $v\wedge v=0$ in $\H^0(X, \wedge^4 W_n)$.
\end{Rem} 

\begin{Lm} If $n=4$ then there is a nondegenerate quadratic form $q: \H^0(X, \wedge^2 W_n)\to k$ such that $\tilde Y_2(L)$ is given by the equation $q(v)=0$ for $v\in \H^0(X, \wedge^2 W_n)$. In this case $a_L=-2$.
\end{Lm}
\begin{Prf} Define $q$ as the composition 
$$
\H^0(X, \wedge^2 W_n)\to \H^0(W, \Sym^2(\wedge^2 W_n))\to\H^0(X, \det W_n),
$$
where the first map sends $v$ to $v\otimes v$. One has $\dim\H^0(X, \det W_n)=1$, so we may view $q$ as a quadratic form with values in $k$. 
  
  Let us first show that the kernel of $q$ is at most 1-dimensional. Pick a subsheaf $\cA_1\oplus \cA_1^*\oplus\cA_2\oplus\cA_2^*\subset W_4$ such that $\cA_i\in\Bun_1^0$, and all the 4 line bundles $\cA_1,\cA_1^*, \cA_2,\cA_2^*$ are pairwise non isomorphic. Let $v_i\in \H^0(\wedge^2(\cA_i\oplus\cA_i^*))$ be a nonzero section. Then $q$ is nondegenerate on the subspace of $\H^0(X, \wedge^2 W_n)$ generated by $v_1, v_2$. So, the kernel of $q$ is at most 1-dimensional. Since $\tilde Y_2(L)-\{0\}$ is smooth, it follows that $q$ is non degenerate. The last assertion is now easy to check.
\end{Prf}  
  
\medskip    
    
    The following lemma is standard, its proof is left to a reader.
    
\begin{Lm} 
\label{Lm_Finkelberg}
Let $V$ be a 3-dimensional $k$-vector space with a nondegenerate quadratic form $q: V\to k$. Let $Y\subset V$ be the closed subscheme given by $q=0$. Let $b: V\,\iso\, V^*$ be the symmetric bilinear form corresponding to $q$. Let $\cV\subset V^*$ be the open subscheme, the complement to the image of $Y$ by $b$. 
There is a unique up to isomorphism rank one and order two local system $\cE_{\cV}$ on $\cV$, which is $\Gm\times \SO(V,q)$-equivariant. Let $\cI_q$ be the intermediate extension of $\cE_{\cV}[3]$ to $V^*$. Then $\Four_{\psi}(\IC(Y))\,\iso\, \cI_q$. 
\end{Lm}    
    
    We need the following generalization of Lemma~\ref{Lm_Finkelberg}. 

\begin{Lm} 
\label{Lm_rank_Four_calculation}
Let $W$ be a 2-dimensional vector space and $d\ge 1$. Let $Y\subset \Sym^d W$ be the closed subscheme, the image of the finite map $W\to \Sym^d W$ given by $w\mapsto w^d$. Then $\Four_{\psi}(\IC(Y))$ generically over $\Sym^d W^*$ is a (shifted) local system of rank $d-1$.
\end{Lm}
\begin{Prf}
For a $f\in \Sym^d W^*$ let $\beta_f: Y\to \A^1$ be the composition $Y\hook{}\Sym^d W\toup{f} \A^1$. Clearly, $\IC(Y)\,\iso\, \Qlb[2]$. For a sufficiently general $f$ calculate the Euler characteristic of $\beta_f^*\cL_{\psi}$. To this end, calculate first $h_!\beta_f^*\cL_{\psi}$, where $h: Y-\{0\}\to \PP(W)$ is the natural $\Gm$-torsor. Here $\PP(W)$ is the projective space of lines in $W$. The details are easy and left to a reader.
\end{Prf}
    
\begin{Lm} 
\label{Lm_degree_cL_Y2}
For any $L\in\Bun_1^1$ one has $a_L=2-n$, that is, $\cL_{Y_2}\,\iso\, \cO(2-n)$ on $\PP^1$.
\end{Lm}
\begin{Prf}
We will show that $\pi_{1,1}^*\cL_{Y_2}$ is of degree $4-2n$. Consider the line bundle $\cL_1$ on $\Bun_1^0$ whose fibre at $E_1\in\Bun_1^0$ is $\Hom(E_1, W_n(L))$. Let $x\in X$ be such that $L\,\iso\, \cO(x)$. Let $r_x: X\to \Bun_1^0$ be the map sending $y$ to $\cO(x-y)$. We claim that $r_x^*\cL_1\,\iso\, \cO((1-n)x)$. This is proved by induction. For $n=1$ one has canonically $\Hom(\cO(x-y), \cO(x))\,\iso\, k$, so the line bunlde $\cL_1$ is trivialized in this case. For $n>1$ consider the exact sequence $0\to \Hom(E_1, \cO_X)\to\Hom(E_1, W_n(L))\toup{\xi} \Hom(E_1, W_{n-1}(L))$. The map $\xi$ between the corresponding line bundles on $\Bun_1^0$ is regular and vanishes only at $E_1\,\iso\, \cO$ with multiplicity one. So,  $r_x^*\cL_1\,\iso\, \cO((1-n)x)$.

 Consider the line bunlde $\cL_2$ on $\Bun_1^0$ with fibre $\Hom(E_1^*, W_n)$ at $E_1\in \Bun_1^0$. One shows by induction that $r_x^*\cL_2\,\iso\, \cO(-(n+1)x)$. Indeed, for $n=1$ we have a regular map $\H^0(X, \cO(2x-y))\to \cO(2x)/\cO(x)$ which is nonzero for $y\ne x$, and has a zero at $y=x$ of order 2. The induction step is as above. 
  
 Write $\ocL_2$ for the total space of $\cL_2$ with zero section removed. Let $\cL_3$ be the line bundle on $\ocL_2$ whose fibre at $(E_1, s)$ is $\Hom(E_1^*, W_n/\Im(s))$, here $s: E_1\hook{} W_n$. Of coarse, $\cL_2$ descends with respect to the projection $\ocL_1\to \Bun_1^0$ sending $(E_1,s)$ to $E_1$, so we view it as a line bundle on $\Bun_1^0$. We have an exact sequence $0\to\Hom(E_1^*, E_1)\to\Hom(E_1^*, W_n)\toup{\nu} \Hom(E_1^*, W_n/\Im(s))$ for any inclusion $s: E_1\hook{} W_n$. The map $\nu$ between line bundles over $\Bun_1^0$ is regular and vanishes exactly at 4 points, so $\deg(\cL_3)=3-n$.
 
 Finally, we obtain $\deg \pi_{1,1}^*\cL_{Y_2}=\deg \cL_1+\deg \cL_3=4-2n$.
\end{Prf}

 \medskip
 
Since $\dim\H^0(X, \wedge^2 W_n)=n-1$, Lemma~\ref{Lm_degree_cL_Y2} implies that there is a 2-dimensional space $E$ and an isomorphism $\Sym^{n-2} E\,\iso\, \H^0(X, \wedge^2 W_n)$ such that $\tilde Y_2(L)$ is the image of the map 
$$
E\toup{\xi} \Sym^{n-2}E\,\iso\, \H^0(X, \wedge^2 W_n),
$$
where $\xi(e)=e^{n-2}$. Now by Lemma~\ref{Lm_rank_Four_calculation}, the Fourier trasform $\Four_{\psi}(\IC(\tilde Y_2(L)))$ is generically a shifted local system of rank $n-3$. This shows that $K^{-1}_{P,\psi}\mid_{\Bun_H^1}$ is generically a (shifted) local system of rank $n-3$. 

 Now the isomorphism (\ref{iso_after_f^1}) combined with Corollary~\ref{Cor_sheaves_cS^m_d} shows that $(f^1)^*\cS^0_{-1}\,\iso\, \Qlb\boxtimes \cL_{\sigma_{n-2}}$. Now applying Proposition~\ref{Pp_fibres_of_nu_P_connected}, we see that $\cK_H\mid_{\Bun_H^1}$ is uniquely determined by the isomorphism (\ref{iso_Th1_on_eBunP}) of Theorem~\ref{Th1}. The fact that for $n=4$, $\cK_H\mid_{\Bun_H^1}$ is generically a (shifted) local system of order two was already established in Lemma~\ref{Lm_Finkelberg}. Proposition~\ref{Pp_case_g=1} is proved.
 
\bigskip

\centerline{\scshape 9. Generalizations for other simple groups}    
    
\bigskip\noindent    
9.1 Inspired by our construction of $\cK_H$ and the results on the Fourier coefficients of minimal representations (\cite{GRS}, Theorem~5.2), we propose the following generalizations of Theorem~\ref{Th1}.

 Let $G$ be a connected simple algebraic group (not necessarily simply connected). Let $P\subset G$ be a maximal parabolic subgroup with an abelian unipotent radical $U\subset P$. Let $M\subset P$ be a Levi subgroup of $P$. Let $P^-$ be the opposite parabolic subgroup with respect to some maximal torus $T\subset M$. Write $U^-$ for the unipotent radical of $P^-$. 
 
  The maximal parabolic subgroups with an abelian unipotent radical have been classified in (\cite{RRS}, list of possible cases in Remark~2.3). So, $G$ is of type $A_n$, $B_n$, $C_n$, $D_n$, $E_6$ or $E_7$. 
  
  We may view $U$ as a linear representation of $M$, write $U^*$ for the dual representation. By \select{loc.cit.}, the group $M$ has finitely many orbits on $U^-\,\iso\, U^*$. On the set of $M$-orbits in $U^*$ one has an order, namely $O_1\le O_2$ iff $O_1 $ is contained in the closure of $O_2$. By (\cite{RRS}, Proposition~2.15), this order is linear, and there is a unique $M$-orbit $Z\subset U^*$ such that the closure $\bar Z$ of $Z$ is $Z\cup\{0\}$. 
    
   Let $\Bun_G$ be the stack of $G$-torsors on $X$, and similarly for $\Bun_P, \Bun_M$. The stack $\Bun_P$ classifies $ \cF_M\in\Bun_M$ and an exact sequence $0\to U_{\cF_M}\to ?\to \cO_X\to 0$ on $X$. Here $U_{\cF_M}=(U\times \cF_M)/M$ is the vector bundle on $X$ obtained out of $U$ by twisting with $\cF_M$. 
 
 Let $\cY_P$ be the stack classifying $\cF_M\in\Bun_M$ and a section $v: U_{\cF_M}\to\Omega$. Then $\Bun_P$ and $\cY_P$ are dual generalized vector bundles over $\Bun_M$, so one has the corresponding Fourier transform functor  $\Four_{\cY_P,\psi}: \D^{\prec}(\cY_P)\to \D^{\prec}(\Bun_P)$. 
 
 The $\Gm$-action on $U^*$ by scalar multiplications commutes with the $M$-action. Let $\cZ_P\subset \cY_P$ be the closed substack classifying $(\cF_M, v)$ such that $v$ is a section of $\bar Z_{\cF_M, \Omega}$. Here $\bar Z_{\cF_M, \Omega}$ is the closed subscheme of the total space of $U^*_{\cF_M}\otimes\Omega$ obtained as the corresponding twisting of $\bar Z$. Let $\cZ_{P,0}\subset\cZ_P$ be the open substack given by the property that $v$ is a section of $Z_{\cF_M,\Omega}\subset \bar Z_{\cF_M, \Omega}$.
 
  Let $^e\Bun_M\subset\Bun_M$ be the open substack given by $\H^0(X, U_{\cF_M})=\H^0(X, \Omega\otimes U_{\cF_M})=0$. Write $^e\cY_P$, $^e\Bun_P$ for the preimage of $^e\Bun_M$ in the corresponding stack. The natural map $\nu_P: {^e\Bun_P}\to\Bun_G$ is smooth. 
  
 If $G$ is of type $C_n$, assume $G$ simply-connected. If $G$ is of type $B_n$ or $C_n$ write $W$ for the standard representation of $G$. Write $\cA$ for the line bundle on $\Bun_G$ with fibre $\det\RG(X, W_{\cF_G})$ at $\cF_G\in\Bun_G$. Let $\Bunt_G$ be the $\mu_2$-gerb over $\Bun_G$ of square roots of $\cA$. It classifies $\cF_G\in\Bun_G$, a 1-dimensional vector space $\cB$ and an isomorphism $\cB^2\,\iso\,\det\RG(X, W_{\cF_G})$. 
  
   We have a diagram
$$
\begin{array}{ccccccc}
^e\cY_P &&&& ^e\Bun_P\\
& \searrow && \swarrow && \searrow\lefteqn{\scriptstyle \nu_P}\\
&& ^e\Bun_M &&&& \Bun_G,
\end{array}
$$ 
where $^e\cY_P$ and $^e\Bun_P$ are dual vector bundles over $^e\Bun_M$. For $G$ of type $C_n$ the map $\nu_P$ lifts naturally to a map $\tilde\nu_P: \Bun_P\to\Bunt_G$.

\begin{Con} 
\label{Con_general_descent}
Assume that $G$ is not of type $B_n$ or $C_n$. There is a perverse sheaf $\cK_G$ on $\Bun_G$ with the following property. If $d\in\pi_1(M)$ and $^e\cZ_{P,0}^d$ is not empty then there exists an isomorphism over $^e\Bun_P^d$
$$
\nu_P^*\cK_G\otimes(\Qlb[1](\frac{1}{2}))^{\dimrel(\nu_P)}\,\iso\, \Four_{\cY_P, \psi}(\IC(\cZ_P))
$$
\end{Con}
\begin{Rem} i) For $G=\Sp_{2n}$ Conjecture~\ref{Con_general_descent} should be corrected as follows. In this case $\cK_G$ is a direct summand in the theta-sheaf on $\Bunt_G$ introduced in (\cite{L1}, Definition~1), and $\nu_P$ should be replaced by $\tilde\nu_P$. With this correction Conjecture~\ref{Con_general_descent} holds for $G=\Sp_{2n}$ (\cite{L1}, Proposition~7).\\
ii) We don't know if Conjecture~\ref{Con_general_descent} should be true for type $B_n$ (even with $\Bun_G$ eventually replaced by $\Bunt_G$). Recall that at the level of functions, a metaplectic $\mu_2$-covering of $\SO_7$ admits a minimal representation (\cite{Ro, To}), whence for $n\ge 4$ it is known that the minimal representation does not exist for $\SO_{2n+1}$ (or its metaplectic coverings).
\end{Rem}

\medskip\noindent
9.2 If $G$ is of type $E_6$ or $E_7$, the perverse sheaf $\cK_H$ from Conjecture~\ref{Con_general_descent} should be the geometric analog of the corresponding minimal representation. More precisely, it should satisfy the Hecke property corresponding to the subregular unipotent orbit in the Langlands dual group $\check{G}$ (precisely as in    
Conjecture~\ref{Con_1}).

\bigskip

\centerline{\scshape Appendix~A. Almost constant local systems on $\Bun_G$}

\bigskip\noindent
A.1 Assume the ground field $k$ algebraically closed. Let $G$ be a semi-simple algebraic group, $G^{sc}$ its simply-connected covering, let $A$ be the finite abelian group defined by the exact sequence $1\to A\toup{i} G^{sc}\to G\to 1$. 

 Pick a connected torus $T$ and an injective homomorphism $\phi: A\hook{} T$, set $T_1=T/A$. Let $G_1=(T\times G^{sc})/A$, where the map $A\to T\times G^{sc}$ is $(\phi, i)$. We get exact sequences $1\to A\to G_1\to T_1\times G\to 1$ and $1\to T\to  G_1\to G\to 1$ over $\Spec k$. 

 Given $b\in \pi_1(G)$, pick any $\bar b\in \pi_1(G_1)$ over $b$ and let $(c, b)\in \pi_1(T_1)\times\pi_1(G)$ be the image of $\bar b$ under $\pi_1(G_1)\to \pi_1(T_1\times G)$. 
Pick any $\cF_{T_1}\in \Bun_{T_1}^c$. Write $\Bun^{\bar b}_{G_1, \cF_{T_1}}$ for the stack $\Bun^{\bar b}_{G_1}\times_{\Bun_{T_1}}\Spec k$, where we used the map $\cF_{T_1}: \Spec k\to \Bun_{T_1}$ to define the fibred product. The projection 
\begin{equation}
\label{map_cF_{T_1}}
f: \Bun^{\bar b}_{G_1, \cF_{T_1}}\to \Bun_G^b
\end{equation}
is smooth and surjective. It is not representable, the group $A$ act on each fibre of $f$ by 2-automorphisms. But after getting rid of this 2-action, the map $f$ becomes a Galois covering of $\Bun_G^b$ with Galois group $\H^1(X, A)$. 

\begin{Def} Say that a local system $K$ on $\Bun^b_G$ is \select{almost constant} if $f^*K$ is constant, that is,
there is $m\ge 1$ and an isomorphism $f^*K\,\iso\, \Qlb^m$.
\end{Def}
 
 Since $T$ is contained in the center of $G_1$, we have a natural action map $\Bun_T\times\Bun_{G_1}\to\Bun_{G_1}$. 
 
 If $\bar b'\in \pi_1(G_1)$ is another element over $b$ let $(c', b)\in \pi_1(T_1)\times\pi_1(G)$ be the image of $\bar b'$. 
Pick any $\cF'_{T_1}\in\Bun^{c'}_{T_1}$. 
Let $\bar\cF$ be the $T_1$-torsor on $X$ of isomorphisms $Isom(\cF_{T_1}, \cF'_{T_1})$. In other words, $\Bun_{T_1}$ is naturally a group stack and $\cF'_{T_1}\,\iso\, \cF_{T_1}\otimes \bar\cF$ for the corresponding product $\otimes$ in $\Bun_{T_1}$. If $\bar\cF\in \Bun_{T_1}^{\bar c}$ then $\bar c$ is in the subgroup $\pi_1(T)\hook{} \pi_1(T_1)$, so one may pick a lifting of $\bar\cF$ to a $T$-torsor $\cF$ on $X$. Then $\cF$ gives rise to a commutative diagram
$$
\begin{array}{ccc}
\Bun^{\bar b}_{G_1, \cF_{T_1}} & \toup{f} & \Bun^b_G\\
\downarrow\lefteqn{\!\wr} & \nearrow\lefteqn{\scriptstyle f'}\\
\Bun^{\bar b'}_{G_1, \cF'_{T_1}},
\end{array}
$$
where the vertical arrow is the action of $\cF$ on $\Bun_{G_1}$. This shows that the notion of an almost constant  local system does not depend on our choices of $\bar b$ and $\cF_{T_1}$.

 One checks that this notion does not depend on a choice of $(\phi, T)$. Note that for $b=0$ an irreducible perverse sheaf $K\in \P(\Bun_G^0)$ is almost constant if its restriction to $\Bun_{G^{sc}}$ is constant. To see this, take $\bar b=0$ and $\cF_{T_1}$ to be a trivial $T_1$-torsor then (\ref{map_cF_{T_1}}) idenitifes with the projection $\Bun_{G^{sc}}\to\Bun_G^0$.
 
 Since $\H^1(X, A)$ is abelian, any almost constant irreducible local system on $\Bun_G^b$ is of rank one and finite order.
 
 The following conjecture was communicated to us by Drinfeld. We have not found a reference for its formulation or a proof\footnote{Conjecture~\ref{Con_pi_1_of_BunG} would follow form the $\ell$-adic version of the results of \cite{G2}.}.
 
\begin{Con}  
\label{Con_pi_1_of_BunG}
Any smooth $\Qlb$-sheaf on $\Bun_{G^{sc}}$ is constant.
\end{Con}

 In view of this conjecture any local system on $\Bun_G$ should be almost constant. 
 
 Consider an Arthur parameter $(\alpha,\sigma):\pi_1(X)\times \SL_2\to \check{G}$, where $\alpha: \pi_1(X)\to Z(\check{G})$ is a homomorphism with values in the center of $\check{G}$, and $\sigma$ corresponds to the principal nilpotent. We have canonically $\Hom(\H^1(X, A), \mu_{\infty})\,\iso\, \H^1(X, Z(\check{G}))$, so $\alpha$ can be seen as a character $\alpha: \H^1(X, A)\to \mu_{\infty}$. We associate to $\alpha$ the local system on $\Bun_G$ whose restriction to $\Bun_G^b$ is the isotypic component in $f_!\Qlb$ for the map (\ref{map_cF_{T_1}}) on which $\H^1(X, A)$ acts by $\alpha$. We expect this local system to be the automorphic sheaf corresponding to the above Arthur parameter.
 
 We will use only the following weaker result.
 
\begin{Pp} 
\label{Pp_appendix_A}
Let $W\in \Rep(\check{G})$, $K$ be an almost constant local system on $\Bun_G$, and $x\in X$. There is $r>0$, almost constant local systems $K_i$ on $\Bun_G$ and $d_i\in\ZZ$ such that
$_x\H^{\la}_G(W,K)\,\iso\, \oplus_{i=1}^r K_i[d_i]$.
\end{Pp}
\begin{Prf}
Pick any $\bar b, \bar b'\in \pi_1(G_1)$. Denote by $b,b'\in \pi_1(G)$ and $c, c'\in \pi_1(T_1)$ the images of $\bar b, \bar b'$ respectively. Pick a $T_1$-torsor $\cF_{T_1}\in \Bun_{T_1}^c$ and set $\cF'_{T_1}=\cF_{T_1}((c'-c)x)$. This makes sense, because $c'-c$ is a coweight of $T_1$. 

 Write $_x\cH_G(b,b')$ for the Hecke stack classifying $\cF_G\in \Bun_G^b$, $\cF'_G\in\Bun_G^{b'}$, $\beta: \cF_G\,\iso\,\cF'_G\mid_{X-x}$.  We have the diagram
$$
\Bun_G^b\;\getsup{h^{\la}}\; {_x\cH_G}(b,b')\;\toup{h^{\ra}}\;\Bun_G^{b'},
$$  
where $h^{\la}$ (resp., $h^{\ra}$) sends the above point to $\cF_G$ (resp., to $\cF'_G$). 
Similarly, we have the stack $_x\cH_{G_1}(\bar b,\bar b')$ included into a diagram of projections
$$
\Bun_{G_1}^{\bar b}\;\getsup{h^{\la}}\; {_x\cH_{G_1}(\bar b,\bar b')}\;\toup{h^{\ra}}\;\Bun_{G_1}^{\bar b'}
$$  
Let $_x\cH(\bar b, \bar b')$ be the stack obtained from $_x\cH_{G_1}(\bar b, \bar b')$ by the base change $(\cF_{T_1}, \cF'_{T_1}): \Spec k\to \Bun_{T_1}^{c}\times \Bun_{T_1}^{c'}$.
We get the diagram
$$
\begin{array}{ccccc}
\Bun^{\bar b}_{G_1, \cF_{T_1}} & \getsup{\bar h^{\la}} & _x\cH(\bar b, \bar b') & \toup{\bar h^{\ra}} & \Bun^{\bar b'}_{G_1, \cF'_{T_1}}\\
\downarrow\lefteqn{\scriptstyle f} && \downarrow\lefteqn{\scriptstyle h} && \downarrow\lefteqn{\scriptstyle f}\\
\Bun_G^b & \getsup{h^{\la}} & {_x\cH_G}(b,b') & \toup{h^{\ra}} & \Bun_G^{b'}
\end{array}
$$
The key observation is that both squares in this diagram are cartesian.  

 We may assume that $K$ is supported on $\Bun_G^{b'}$. Let $\cS$ be the spherical perverse sheaf on $\Gr_G$ corresponding to $W$. By definition, 
$$
\H^{\la}_G(\cS, K)\,\iso\, h^{\la}_!(\ast\cS\tboxtimes K)^r
$$
(see \cite{L2}, Section~2.2.1). We may assume that $\ast\cS\tboxtimes K$ is supported on $_x\cH_G(b,b')$. The above diagram yields an isomorphism
$$
_x\H^{\la}_G(\cS,  f_!\Qlb)\,\iso\, f_! {_x\H^{\la}_{G_1}}(\cS, \Qlb)
$$
Since ${_x\H^{\la}_{G_1}}(\cS, \Qlb)$ is a constant complex on $\Bun^{\bar b}_{G_1, \cF_{T_1}}$, our claim follows.
\end{Prf} 
 
\bigskip 
 
\centerline{\scshape Appendix B. Connectedness issues} 
  
\bigskip\noindent
B.1 Assume the ground field $k$ algebraically closed. Let $G$ be a semi-simple algebraic group. We pick a Borel subgroup $B\subset G$ and its maximal torus $T_G\subset B$, write $\check{\Lambda}$ (resp., $\Lambda$) for the weights (resp., coweights) lattice of $T_G$. Write $\check{\Lambda}^+$ for the dominant weights of $(G,T_G)$, let $w_0$ be the longuest element of the Weyl group of $(G,T_G)$. Let
$B\subset P\subset G$ be a standard parabolic subgroup, $U_P\subset P$ its unipotent radical, write $M$ for the corresponding standard Levi subgroup of $P$. 

 For $\lambda\in\Lambda$ write $\Bun_{T_G}^{\lambda}$ for the connected component of $\Bun_{T_G}$ classifying $\cF\in\Bun_{T_G}$ such that for each $\check{\lambda}\in \check{\Lambda}$ one has $\deg\cL^{\check{\lambda}}_{\cF}=\<\lambda, \check{\lambda}\>$.
For $d\in\pi_1(M)$ write $\Bun_M^{d}$ for the connected component of $\Bun_M$ containing the image of $\Bun_{T_G}^{\lambda}$ for any $\lambda\in\Lambda$ over $d$. Write $\Bun_P^d$ for the preimage of $\Bun_M^d$ under the natural map $\Bun_P\to \Bun_M$. Write $\gg$ (resp., $\gp$) for the Lie algebras of $G$ (resp., of $P$). Let $^0\Bun_P^d\subset\Bun_P^d$ be the open substack classifying $\cF\in\Bun_P^d$ such that for any irreducible $P$-submodule $V$ of $\gg/\gp$ one has $\H^1(X, V_{\cF})=0$. The natural map $\nu_P: {^0\Bun_P^d}\to\Bun_G$ is smooth.

\begin{Pp} 
\label{Pp_appendix_B}
Let $d\in \pi_1(M)$, write $b\in \pi_1(G)$ for the image of $d$ in $\pi_1(G)$. Assume that there is a lifting of $d$ to an anti-dominant coweight $\lambda\in\Lambda$ such that 
for each negative root $\check{\alpha}$ of $(G, T_G)$ one has $2g-2<\<\check{\alpha}, \lambda\>$. Then the generic fibre of $\nu_P: {^0\Bun_P^d}\to\Bun_G^b$ is geometrically connected.
So, there is a non empty open substack of $\Bun_G^b$ such that each fibre of the latter map over a point of this substack is geometrically connected.
\end{Pp}
\begin{Prf}
Pick $T$, $\phi: A\to T$ and define $T_1$, $G_1$ as in Section~A.1. Let $P_1$ (resp., $B_1$, $M_1$) be the preimage of $P$ (resp., of $B$, $M$) under $G_1\to G$. The diagram is cartesian 
$$
\begin{array}{ccccc}
\Bun_{B_1} & \to &
\Bun_{P_1}&\toup{\nu_{P_1}} & \Bun_{G_1}\\
\downarrow && \downarrow && \downarrow\\
\Bun_B & \to &\Bun_P & \toup{\nu_P} & \Bun_G
\end{array}
$$
Let $\Lambda_1$ be the coweights lattice of $B_1$, the projection $\Lambda_1\to\Lambda$ is surjective. Pick $\bar \lambda\in \Lambda_1$ over $\lambda$, write $\bar d$ (resp., 
$\bar b$) for the image of $\bar \lambda$ in $\pi_1(M_1)$ 
(resp., in $\pi_1(G_1)$). It suffices to prove that the generic fibre of $\nu_{P_1}: {^0\Bun^{\bar d}_{P_1}}\to \Bun^{\bar b}_{G_1}$ is geometrically connected. The second assertion would also follow using Lemma~\ref{Lm_EGA4}.
  
  Note that $\bar \lambda$ is anti-dominant for $G_1$, and for each negative root $\check{\alpha}$ of $G_1$, we have $2g-2<\<\check{\alpha},\lambda\>$. This implies that the map $\nu_{B_1}: \Bun_{B_1}^{\bar\lambda}\to \Bun_{G_1}^{\bar b}$ is smooth, and similarly for $\Bun_{B_1}\to\Bun_{P_1}$. So, it suffices to show that the generic fibre of $\nu_{B_1}: \Bun_{B_1}^{\bar\lambda}\to \Bun_{G_1}^{\bar b}$ is connected. 
By Lemma~\ref{Lm_appendix_B_second_lemma} below, it suffices to show that $\Bun_{B_1}^{\bar\lambda}\times_{\Bun_{G_1}^{\bar b}}\Bun_{B_1}^{\bar\lambda}$ is connected. 

 The stack $\Bun_{B_1}^{\bar\lambda}$ is smooth of dimensioon $(g-1)\dim B_1-\<\bar\lambda, 2\check{\rho}\>$, and $\Bun_G$ is smooth of dimension $(g-1)\dim G$. Here $\check{\rho}$ is the half sum of positive roots of $(G_1, B_1)$. 
So, the stack $\Bun_{B_1}^{\bar\lambda}\times_{\Bun_{G_1}^{\bar b}}\Bun_{B_1}^{\bar\lambda}$ is smooth of pure dimension $(g-1)\dim T_{G_1}-2\<\bar\lambda, 2\check{\rho}\>$, here $T_{G_1}$ is the preimage of $T_G$ under $G_1\to G$. Let 
$$
j: \cY\hook{} \Bun_{B_1}^{\bar\lambda}\times_{\Bun_{G_1}^{\bar b}}\Bun_{B_1}^{\bar\lambda}
$$ 
be the open substack given by the property that the two $B_1$-structures on a given $G_1$-torsor are transversal at the generic point of $X$. One checks that the complement $\cY'$ of $\cY$ is of dimension $<(g-1)\dim T_{G_1}-2\<\bar\lambda, 2\check{\rho}\>$. Thus, it suffices to prove that $\cY$ is connected.
To do so, consider the map $q: \cY\to \Bun_{T_{G_1}}^{\bar\lambda}$ defined as the composition 
$$
\cY\toup{j} \Bun_{B_1}^{\bar\lambda}\times_{\Bun_{G_1}}^{\bar b}\Bun_{B_1}^{\bar\lambda}\toup{\pr_2}\Bun_{B_1}^{\bar\lambda}\to \Bun_{T_{G_1}}^{\bar\lambda},
$$
here $\pr_2$ is the projection on the second factor. Since $\Bun_{T_{G_1}}^{\bar\lambda}$ is smooth and irreducible and $q$ is smooth of constant relative dimension, it siffices to show that each fibre of $q$ is connected.

 The fibre of $q$ over any $\cF\in \Bun^{\bar\lambda}_{T_{G_1}}$ is isomorphic to the twisted versions of the Zastava space
$\ocZ{}^{\theta}_{\gg_1, \bb_1}(X)$ in the notation of (\cite{BFG}, Section~2.12) corresponding to the parameter $\theta=w_0(\bar\lambda)-\bar\lambda$. The twist is due to the fact that the trivial $T_{G_1}$-torsor used in the definitrion of the Zastava is replaced by an arbitrary point of $\Bun^{\bar\lambda}_{T_{G_1}}$. Note that $\theta$ is a sum of positive coroots. Now (\cite{BFG}, Proposition~2.25) combined with (\select{loc.cit.}, Propositions~2.19 and 2.21) imply that each fibre of $q$ is connected. We have also used the fact that the derived group $[G_1, G_1]$ of $G_1$ is simply-connected, as the results of \cite{BFG} require this assumption.  
\end{Prf}  
  
\medskip

 The following is proved in (\cite{Groth}, Proposition~9.7.8, p. 82). 

\begin{Lm} 
\label{Lm_EGA4}
Let $S$ be an irreducible scheme with generic point $\eta$. Let $f: Y\to S$ be a finitely presented morphism of schemes. Assume that the fibre $f^{-1}(\eta)$ is geometrically irreducible (resp., geometrically connected). Then there is a non empty open subscheme $U\subset S$ such that 
for any $s\in U$ the fibre $f^{-1}(s)$ is geometrically irreducible (resp., geometrically connected). \QED
\end{Lm}
 
\begin{Lm} 
\label{Lm_appendix_B_second_lemma}
Let $S$ and $Y$ be smooth irreducible $k$-schemes, let $f: Y\to S$ be a smooth morphism (of some constant relative dimension). Assume that $Y\times_S Y$ is connected. Then there is a non empty open subscheme $U\subset S$ such that for each $s\in U$ the scheme $f^{-1}(s)$ is geometrically connected.
\end{Lm}
\begin{Prf} Write $\eta$ (resp., $\zeta$) for the generic point of $S$ (resp., of $Y$). Set $Y_{\eta}=f^{-1}(\eta)$. Since $Y\times_S Y$ is connected and smooth over $k$, it is irreducible. Let $\nu$ be the generic point of $Y\times_S Y$.

 Consider the projection $\pr_2: Y\times_S Y\to Y$. The generic fibre of $\pr_2$ over $\zeta$ is $Y_{\eta}\times_{\eta}\zeta$. Since $\nu$ is dense in $Y_{\eta}\times_{\eta}\zeta$, the scheme $Y_{\eta}\times_{\eta}\zeta$ is irreducible. Clearly, it is actually geometrically irreducible. Now apply Lemma~\ref{Lm_EGA4} to $\pr_2$.
\end{Prf}

\bigskip

\centerline{\scshape Appendix C. Case of characteristic zero}

\bigskip\noindent
C.1 Assume the base field $k$ algebraically closed of characteristic zero. Work with $\cD$-modules instead of \'etale $\Qlb$-sheaves. In this case one uses the homogeneous Fourier transform, so we omit $\psi$ in some notations, for example the $D$-module $K_{P,\psi}$ introduced in Section~2.3.1 is now denoted $K_P$. 

Theorem~\ref{Th1} holds in the characteristic zero case. In this appendix we briefly explain the changes to be made in our proof of Theorem~\ref{Th1} for $\cD$-modules.
  
  Write $\gh$ for the Lie algebra of $H$. The cotangent bundle $T^*\Bun_{H}$ is the stack classifying $(V,\sigma)$, where $V\in \Bun_H$ and $\sigma\in \H^0(X, \gh_{V}^{*}\otimes \Omega)\,\iso\, \Hom(\wedge^2 V, \Omega)$ (cf. \cite{BD}). Write $\bar Z\subset \gh^*$ for the closure of the minimal nilpotent orbit $Z$ in $\gh^*$, the complement of $Z$ in $\bar Z$ is the origin in $\gh^*$. 
  
  Let $\cC\subset T^{*}\Bun_{H}$ be the substack classifying $(V,\sigma)$ such that $\sigma$ is a section of $\bar Z_{V,\Omega}\to X$. Here $\bar Z_{V,\Omega}\subset  \gh_{V}^{*}\otimes \Omega$ is the closed subscheme obtained as the corresponding twist of $\bar Z$. Then $\cC$ contains the zero section of $T^*\Bun_{H}$, write $\cC'$ for the complement of this zero section in $\cC$. Then $\cC'$ admits a stratification by locally closed substacks $\cC^m$, $m\ge 0$. Here $\cC^m$ is the stack classifying $V\in \Bun_H$ with an isotropic subbundle $V_2\subset V$ of rank two, $D\in X^{(m)}$, and an isomorphism $\Omega(-D)\,\iso\, \det(V/V_{-2})$. Here we have denoted by $V_{-2}\subset V$ the orthogonal complement to $V_2$ in $V$. The stack $\cC^m$ is smooth of dimension $\dim \Bun_{H}-(2n-4)m$. Since $\cC^0$ is contained in the global nilpotent cone,
from (\cite{BD}, Theorem~2.10.2) one derives that $\cC^0$ is a lagrangian substack of $T^{*}\Bun_{H}$.

   The formulation of Proposition~\ref{Pp_K^d_coincide} is changed as follows.
   
\begin{Pp} 
\label{Pp_appedix_C_first}
The irreducible subquotients $\cK^d_{\cU}$ of $\tilde\cK_{\cU}$ over $\cU^b_H$ all coincide for $d\mod 2=b$. The resulting irreducible subquotient is denoted $\cK_{\cU,b}$. If $F$ is a different irreducible subquotient of $\tilde\cK_{\cU}$ over $\cU^b_H$ then $\bar F$ is a (shifted) local system over the whole of $\Bun_H^b$. 
\end{Pp}

 To prove Proposition~\ref{Pp_appedix_C_first}, keep only the following part of Lemma~\ref{Lm_really_key_lemma} (its proof holds without changes in characteristic zero case).
 
\begin{Lm}  
\label{Lm_appendix_C_first}
Let $\bar F$ be an irreducible $\cD$-module on $\Bun_H^b$ for some $b\in\ZZ/2\ZZ$. Let $I$ be an infinite bounded from above set of integers. Assume given for each $d\in I$ a $\cD$-module $\cF^d$ on $^e\Bun_n^d$ and an isomorphism (\ref{iso_descent_on_Bun_P^d}) over $^e\Bun_P^d$. Assume that if $d\in I$ then $\nu_P^*(\bar F)$ is nonzero over $^e\Bun_P^d$. Then $\bar F$ is a (shifted) local system on $\Bun_H^b$. \QED
\end{Lm}

 If $E\to S$ is a vector bundle, there is a canonical symplectomorphism $\iota:T^*(E)\to T^*(E^*)$ between the cotangent bundles to the corresponding total spaces. If $E\,\iso\, S\times E_0$ is a trivialization of $E$ over $S$, here $E_0$ is a vector space, then $T^*E\,\iso\,T^{*}S\times E_0
\times E_0^*$ and  $T^*E^*\,\iso\,T^*S\times E_0^*
\times E_0$ natually. In this case $\iota (a,x,y)=(a,y,-x)$ for $a\in T^*S$, $x\in E_0$ and  $y\in E^*_0$. If is remarkable that this symplectomorphism does not depend on the trivialization of $E$. 
For a $\Gm$-equivariant $\cD$-module $M$ on $E$ the characteristic variety of its Fourier transform $\Four(M)$ on $E^*$ is the image under $\iota$ of the characteristic variety of $M$ (cf. \cite{KSh}, Theorem~5.5.5). 
 
  Recall the stack $^e\cZ_{P,0}$ introduced in Section~2.3.1, it classifies $U\in {^e\Bun_n}$, $M\in\Bun_{G_1}$ and a surjection $U\to M$. Let $\cT\subset {^e\cZ_{P,0}}$ be the open substack given by $\H^0(X, U\otimes M)=0$. Recall that $^e\Bun_P$ and $^e\cY_P$ are dual vector bundles over $^e\Bun_n$, denote by $\iota: T^*(^e\cY_P)\to T^*(^e\Bun_P)$ the symplectomorphism as above. 
  
   The substack $\cT\subset {^e\cY_P}$ is locally closed, and we denote by $N^*_{\cT}(^e\cY_P)$ the conormal bundle of $\cT$ in $^e\cY_P$. 
      
\begin{Lm} 
\label{Lm_appendix_C_second}
Let $b\in\ZZ/2\ZZ$, let $\cF$ be an irreducible $\cD$-module on $\Bun_H^b$. Assume that $d\in Z(e,P)$ with $d\mod 2=b$, and over $^e\Bun_P^d$ the $\cD$-module $\nu_P^*(\cF)[\dimrel(\nu_P)]$ contains $K^d_P$ as an irreducible subquotient. Then $\cC^0\cap T^*(\Bun_H^b)$ is contained in the characteristic variety of $\cF$. 
\end{Lm}
\noindent
\select{Sketch of the proof} \  One checks that the substack $\iota(N^*_{\cT}(^e\cY_P))\subset T^*(^e\Bun_P)$ is contained in 
$$
\cC^0\times_{\Bun_H}{^e\Bun_P}\subset T^*\Bun_H\times_{\Bun_H}{^e\Bun_P}\subset  T^*(^e\Bun_P)
$$
This implies that $\cC^0\times_{\Bun_H}{^e\Bun_P^d}$ is contained in the characteristic variety of $\nu_P^*(\cF)[\dimrel(\nu_P)]$. 
\QED 
  
  \medskip
  
\noindent  
\select{Sketch of the proof of Proposition~\ref{Pp_appedix_C_first}} Combining Lemmas~\ref{Lm_appendix_C_first} and \ref{Lm_appendix_C_second} one gets the following. For each $d\in Z(e,P)$ the characteristic variety of $\cK^d_{\cU}$ contains $\cC^0\times_{\Bun_H} \cU^b_H$. 
The $\cD$-module $\tilde\cK_{\cU}$ over $\cU^b_H$ admits a unique irreducible subquotient, whose characteristic variety contains 
$\cC^0\times_{\Bun_H} \cU^b_H$,
and all its other irreducible subquotients are (shifted) local systems on $\cU^b_H$. 
\QED

\begin{Con} The characteristic variety of $\cK_H$ is $\cC^0$ to which one should possibly add the zero section of $T^*\Bun_H\to\Bun_H$.
\end{Con} 

\begin{Rem} We have also proved the following (these claims are not used in the present paper, so we don't provide a proof). \\
i) If $n\ge 4$ then the stack $\cC^0$ is smooth, and $\cC^0\cap T^*(\Bun_H^b)$ is irreducible for each $b\in\ZZ/2\ZZ$. If $n=2$ then $\cC^0$ is contained in $T^*(\Bun_H^0)$ and is irreducible. For $n=3$ the stack $\cC^0$ is not irreducible. \\
ii) Write $\bar\cC^0$ for the image of $\cC^0$ under the projection $T^*\Bun_H\to \Bun_H$. If $n=2$ then $\bar\cC^0\subset\Bun_H^0$ is of codimension one and irreducible.  If $n\ge 4$ then for each $b\in\ZZ/2\ZZ$ the substack $\bar\cC^0\cap \Bun_H^b$ is irreducible and of codimension one in $\Bun_H^b$. \\
iii) Let $g>1$ be odd. Then, for $d$ sufficiently small, Braden's condition (\cite{Brad}, Corollary~3) holds for $\IC(\cZ_P)$ over $^e\cY^d_P$, so that $K^d_P$ does not vanish at the generic point of $\Bun_P^d$. Thus, $\cK_H$ does not vanish at the generic point of $\Bun_H^b$ for each $b$ in this case. In particular, the isomorphism (\ref{iso_Th1_on_eBunP}) determines $\cK_H$ up to a unique isomorphism.
\end{Rem}

\end{document}